\documentclass[12pt]{article}

\usepackage[dvips]{color}
\usepackage[OT1]{fontenc}
\usepackage{anysize}
\usepackage[english]{babel}
\usepackage{multirow}
\usepackage{rotating}
\usepackage{booktabs}
\usepackage{indentfirst}
\usepackage{float}
\usepackage{latexsym}
\usepackage{geometry}
\usepackage{xcolor}
\usepackage{authblk}
\usepackage{pdflscape}

\usepackage{amsmath}
\usepackage{mathrsfs}
\usepackage{amsfonts}
\usepackage[authoryear]{natbib}
\usepackage{amsthm}
\usepackage{amssymb}
\usepackage{color}
\usepackage{lscape}
\usepackage{multirow}
\usepackage{rotating}
\usepackage{caption2}
\usepackage{subfigure}
\usepackage{graphicx}
\theoremstyle{plain}

\def\be{\begin{equation}}
\def\ee{\end{equation}}
\def\bea{\begin{eqnarray}}
\def\eea{\end{eqnarray}}
\def\bd{\begin{displaymath}}
\def\ed{\end{displaymath}}
\def\bda{\begin{eqnarray*}}
\def\eda{\end{eqnarray*}}
\def\bsm{\begin{small}}
\def\esm{\end{small}}

\def\t0{\theta_0}

\def\ha1{\hat \beta_1}

\def\bnt{\begin{enumerate}}
\def\ent{\end{enumerate}}

\def\bsc{\begin{scriptsize}}
\def\esc{\end{scriptsize}}

\newcommand{\p}{{\rm p}}

\newtheorem{theorem}{Theorem}

\newtheorem{lemma}{Lemma}

\newtheorem{proposition}{Proposition}
\newtheorem{ft}{Fact}[section]
\theoremstyle{definition}
\newtheorem{as}{Assumption}

\newtheorem{rek}{Remark}

\def\boxit#1{\vbox{\hrule\hbox{\vrule\kern6pt
          \vbox{\kern6pt#1\kern6pt}\kern6pt\vrule}\hrule}}

\makeatletter
\newcommand{\figcaption}{\def\@captype{figure}\caption}
\newcommand{\tabcaption}{\def\@captype{table}\caption}
\makeatother

\newcommand{\cov}{{\rm Cov}}
\newcommand{\diag}{{\rm diag}}

\newcommand{\bA}{{\mathbf A}}
\newcommand{\bB}{{\mathbf B}}

\newcommand{\bI}{{\mathbf I}}

\newcommand{\bR}{{\mathbf R}}

\newcommand{\bU}{{\mathbf U}}

\newcommand{\bW}{{\mathbf W}}
\newcommand{\bX}{{\mathbf X}}
\newcommand{\bY}{{\mathbf Y}}
\newcommand{\bZ}{{\mathbf Z}}
\newcommand{\ba}{{\mathbf a}}

\newcommand{\bfeta}  {\boldsymbol{\eta}}

\newcommand{\bOmega}{\boldsymbol{\Omega}}

\newcommand{\bSigma}{\boldsymbol{\Sigma}}

\newcommand{\bgamma}{\boldsymbol{\gamma}}

\newcommand{\btheta} {\boldsymbol{\theta}}

\newcommand{\bD}{{\mathbf D}}
\newcommand{\bzero}{{\mathbf 0}}

\newcommand{\bchi}{\boldsymbol{\chi}}

\usepackage{fancyhdr}
\renewcommand{\baselinestretch}{1.2}

\renewcommand{\baselinestretch}{\vv}

\renewcommand{\theequation}{\arabic{equation}}
\allowdisplaybreaks

\begin{document}

\renewcommand{\baselinestretch}{1.1}\normalsize

\title{\Large \bf   Optimal Covariance Matrix Estimation for  High-dimensional  Noise in  High-frequency Data}

\author[a,b]{Jinyuan Chang}
\author[b]{Qiao Hu}
\author[c]{Cheng Liu}
\author[d]{Cheng Yong Tang}
\affil[a]{\it \small Guizhou Key Laboratory of Big Data Statistical Analysis, Guizhou University of Finance and Economics, Guiyang, China}
\affil[b]{\it \small Joint Laboratory of Data Science and Business Intelligence, Southwestern University of Finance and
Economics, Chengdu, China}
\affil[c]{\it School of Economics and Management, Wuhan University, Wuhan, China}
\affil[d]{\it \small Department of Statistical Science, Temple University, Philadelphia, PA, USA}

 \date{}
\maketitle

\vspace{-0.5cm}

\begin{abstract}

We consider high-dimensional measurement errors with  high-frequency data.
Our objective is on recovering the high-dimensional cross-sectional covariance matrix of the random  errors with optimality. In this problem,  not all components of the random vector are observed at the same time and the measurement errors are latent variables, leading to  major challenges besides high data dimensionality.
We propose a  new covariance matrix estimator in this context with appropriate localization and thresholding, and then
conduct a series of comprehensive theoretical investigations of the proposed estimator.
By developing a new technical device integrating the high-frequency data feature with the conventional notion of $\alpha$-mixing, our analysis successfully accommodates the challenging serial dependence in the measurement errors.
Our theoretical analysis establishes the minimax optimal convergence rates associated with two commonly used loss functions; and we demonstrate with concrete cases  when the proposed localized estimator with thresholding achieves the minimax optimal convergence rates.
 Considering that the variances and covariances can be small in reality, we conduct a second-order theoretical analysis that further disentangles  the dominating bias in the estimator. A bias-corrected estimator is then proposed to ensure its practical finite sample performance.
We also extensively analyze our estimator in the setting with jumps, and show that its performance is reasonably robust. 	
	We illustrate the  promising empirical performance of the proposed estimator with extensive simulation studies and a real data analysis.


\end{abstract}




\baselineskip=21pt


\renewcommand{\baselinestretch}{1.2}\normalsize

\section{Introduction}

High-frequency data broadly refer to those collected at time points with very small time intervals between consecutive observations.  Exemplary scenarios with high-frequency data include   longitudinal observations with intensive repeated measurements   \citep{BolgerLau_2013},  the tick-by-tick trading data in finance \citep{Zhangetal_2005_JASA}, and functional data with dense observations \citep{Zhang2016}.   High-frequency data are commonly contaminated by some noise,   broadly termed  as the measurement errors. For measurement errors in the context of functional data analysis, we refer to the review article \cite{Wang2016a} and reference therein.  In high-frequency financial data,  as another example, the microstructure noise is well known; see the monograph \cite{Ait-Sahalia2014} for an overview.  

Despite the central interests on recovering the signals contaminated by the noise,  the properties of the noise themselves are of their own great interests.   
Recently, \cite{Jacod2017} highlighted the importance of statistical
properties of the microstructure noise and studied the estimation of its moments; see the  recent study of  \cite{LiLinton} on the limiting distributions in a broad setting.
\cite{Chang2018} investigated recovering the distribution of the noise with some frequency-domain analysis.
In a simultaneous and independent work of ours, \cite{DaXiu2021} investigated the auto-covariance of the measurement errors with a semiparametric approach that utilizing a working moving-average model.
These aforementioned studies aimed on univariate cases.
\cite{UO09} considered the  covariance estimation  and testing for measurement errors in a bivariate case.
 \cite{Christensen2013} proposed an estimator for the covariance matrix of the noise vector in high-frequency finance data.   Both \cite{UO09}  and  \cite{Christensen2013} handled fixed dimensional cases with $m$-dependent or independent measurement errors. 

We are motivated to concentrate on high-dimensional cases in this study that shed light on influential practical applications where covariances between different components of the noise could bring us  useful information in solving various problems.
For example, for functional-type observations, the covariations between the measurement errors may help identifying the common source or reasons of contaminations so that improvement can be developed in designing future investigations. For financial data, such covariations in the high-dimensional microstructure noise may help in better understanding the trading behaviors that may show substantial different pattern between equities. 
Indeed, \cite{Li2016} found that a parametric function incorporating the market information may account for a substantial contribution to the variations in the microstructure noise.  
Nevertheless,  studying the covariance between different components of the high-dimensional noise in high-frequency data remains little explored.  

Our primary interests in this study are on the validity and optimality of the covariance matrix estimation procedure for the high-dimensional noise in high-frequency data. This problem has unique challenges from multiple aspects.  First,  since the noise of interest are not directly observable, the targeted random vectors are latent. 
Second, the latency arises together with high data dimensionality and high sampling frequency, two challenging features that interrelates to each other in this investigation.
The high-dimensional noise sequence is expected to contain serial dependence, posing a major methodological and theoretical challenge.
The properties of high-dimensional covariance matrix estimation have not yet been explored in this important scenario.
Third,  the high-dimensional observations may not be synchronous, i.e. different components of the contaminated observation for the high-dimensional noise may be observed at different time points.
How these data  features  affect the statistical properties  on the validity and optimality  of the covariance matrix estimation remains unclear.

High-dimensional covariance matrix estimation is an important problem in the current state of knowledge,  and has received intensive attentions in the past decade; see, among others, \cite{BickelLevina_2008_AOS, BickelLevina_2008_thre},  \cite{LamFan_2009_AOS}, \cite{Rothman2009}, \cite{Cai2010}, \cite{Cai2011} and \cite{Cai2012a,Cai2012}.   For high-dimensional sparse covariance matrices,     the minimax optimality  of the estimations were investigated in-depth in \cite{Cai2012a,Cai2012}. We note that the existing estimation methods for high-dimensional sparse covariance matrices are developed when the underlying data of interest are fully observed;  hence  they are not applicable for the covariance matrix estimation of the noise in high-frequency data with latency and asynchronous observations.
In the literature on multivariate and high-dimensional high-frequency data analysis, existing studies mainly concern the estimations of the so-called realized covariance matrix. Specifically, the major objective is on the signal part,  attempting to eliminate the impact from the noise; see, for example,   \cite{Aitetal_2010_JASA}, \cite{Fan2012},  \cite{Taoetal_2013_AOS}, \cite{LiuTang_2014_JOE}, \cite{Lametal_2017_Biometrika},  and \cite{Xia2018}.  However, it remains little explored on the high-dimensional covariance matrix of the noises in high-frequency data, accommondating all aforementioned challenging features.

Our study makes  several  contributions to the area.  To our best knowledge, our method is  the first handling covariance matrix estimation of the serially dependent high-dimensional noises in high-frequency data.  Methodologically, to overcome the difficulties due to the latency,  asynchronicity, and serially dependent observations, we propose a new approach with appropriate localization and thresholding.
Theoretically,  to our best knowledge,  our technical device integrating high-frequency serial dependence with the $\alpha$-mixing is a new development of the current state of knowledge; and it can be more broadly applied for solving this class of problems.
Meanwhile, our theoretical analysis establishes the minimax optimal convergence rates associated with two commonly used loss functions for the covariance matrix estimations of the high-dimensional noise in high-frequency data. The minimax optimal rates in this setting are our new theoretical discoveries, and we establish cases  when the proposed estimator achieves such rates. Our result also reveals that the optimal convergence rates reflect the impact due to the asynchronous data,   which are slower than those with synchronous data.
The higher the level of the data asynchronicity is, the slower the convergence rates are expected.
We show that the proposed localized estimator has the same accuracy as if the high-dimensional noise are directly observed in the sense of the same convergence rates. Furthermore, our theory includes a second-order analysis revealing the dominating bias of the estimator.  We then propose a bias-corrected estimator and show that removing such a bias leads to more promising performance, especially when components in the covariance matrix are small. Our analysis also indicates that the proposed localized estimator is robust to the setting with jumps.


The rest of this paper is organized as follows. The methodology is outlined in Section \ref{s2}, followed by theoretical development in Section \ref{sec:theoryfU}. Section \ref{sec:4} presents the theory and method handling situation when the level of the noises is small. Section \ref{s5j} investigates the robustness of our method in the setting with jumps. Numerical studies with simulation and a real data analysis are presented in Section \ref{s3}. Section \ref{s5} includes some discussions. All technical proofs are given in Section \ref{proof}. Some additional numerical results are presented in the supplementary material.

\section{Methodology}
\label{s2}

\subsection{Model and data}\label{sec:model}

We introduce some notations first. For any positive integer $q$, we write $[q]:=\{1,\ldots,q\}$. For a matrix $\bB=(b_{i,j})_{s_1\times s_2}$, let $\bB^{\otimes2}=\bB\bB^{\top}$, $|\bB|_\infty=\max_{i\in[s_1],j\in[s_2]}|b_{i,j}|$, $\|\bB\|_1=\max_{j\in[s_2]}\sum_{i=1}^{s_1}|b_{i,j}|$, $\|\bB\|_\infty=\max_{i\in[s_1]}\sum_{j=1}^{s_2}|b_{i,j}|$ and $\|\bB\|_2=\lambda_{\max}^{1/2}(\bB^{\otimes2})$, where $\lambda_{\max}(\bB^{\otimes2})$ denotes the largest eigenvalue of $\bB^{\otimes2}$. Denote by $I
(\cdot)$ the indicator function. For a countable set $\mathcal{G}$, we use $|\mathcal{G}|$ to denote its cardinality. For two sequences of positive numbers $\{a_n\}$ and $\{b_n\}$, we write $a_n\lesssim b_n$ or $b_n\gtrsim a_n$ if there exist a positive constant $c$ and a large enough integer $n_0$ such that $a_n/b_n\leq c$ for all $n\geq n_0$. We write $a_n\asymp b_n$ if and only if $a_n\lesssim b_n$ and $b_n\lesssim a_n$ hold simultaneously.

The setting of our study contains the signal part -- a $p$-dimensional continuous-time process $(\bX_t)_{t\in[0,T]}$,  where, without loss of generality, $[0,T]$ is the time frame in which the high-frequency data are observed.
We begin with a setting that $\bX_t=(X_{1,t},\ldots,X_{p,t})^{\top}$ satisfies:
\begin{equation}\label{eq:model}
{\rm d} X_{i,t}=\mu_{i,t}\,{\rm d}t+\sigma_{i,t}\,{\rm d}B_{i,t}~~~\textrm{and}~~~\mathbb{E}({\rm d}B_{i,t}\cdot{\rm d}B_{j,t})=\rho_{i,j,t}\,{\rm d}t\,,
\end{equation}
where $\mu_{i,t}$ and $\sigma_{i,t}$ are progressively measurable processes, and $B_{1,t},\ldots,B_{p,t}$ are univariate standard Brownian motions. Here $\sigma_{i,t}$ and $\rho_{i,j,t}$ are, respectively, governing the volatilities and correlations, where both of them may be dynamic over time.   A theoretical study of our method in the setting with jumps will be considered in Section \ref{s5j}.

For each $i\in[p]$, we use $\mathcal {G}_i=\{t_{i,1},\ldots,t_{i,n_i}\}$ to denote the grid of time points at which we observe the noisy data of  the $i$th component process $X_{i,t}$, where $0\leq t_{i,1}<\cdots<t_{i,n_i}\leq T$.
The subject-specific set $\mathcal {G}_i$ reflects the asynchronous nature of the problem.
For the special case with synchronous data, all $\mathcal{G}_i$'s are the same. However, $\mathcal{G}_i$'s are typically different in many practical high-frequency data. Let $n$ be the number of different time points in $\cup_{i=1}^p\mathcal{G}_i$, and we denote the different time points in $\cup_{i=1}^p\mathcal{G}_i$ by $0\leq t_1<\cdots<t_n\leq T$. For any $i,j\in[p]$,  we define
\[
n_{i,j}=|\mathcal{G}_i\cap\mathcal{G}_j|\,,
\]
where $n_{i,j}$ evaluates how many time points $t_k$'s  at which we observe the noisy data of the $i$th and $j$th component processes $X_{i,t}$ and $X_{j,t}$ simultaneously. Clearly,  $n_{i,i}=n_i$ for any $i\in[p]$.

We consider  that the actual observed data are  contaminated by additive measurement errors in the sense that
\[
Y_{i,t_{i,k}}=X_{i,t_{i,k}}+U_{i,t_{i,k}}
\]
with $\mathbb{E}(U_{i,t_{i,k}})=0$ for each $i\in[p]$ and $k\in[n_i]$. The additive noise assumption is common in the literature; see \cite{Ait-Sahalia2014}. Formally, we can write
\begin{equation}\label{eq:mod}
\bY_{t_k}=\bX_{t_k}+\bU_{t_k}\,,~~~~k\in[n]\,,
\end{equation}
and assume the measurement errors $\{\bU_{t_k}\}_{k=1}^n$ are independent of the process $(\bX_t)_{t\in[0,T]}$. At each time point $t_k$, we only observe $\sum_{i=1}^pI(t_k\in\mathcal{G}_i)$ components of $\bY_{t_k}$.


Besides the cross-sectional dependence,  serial dependence is expected to be the case for $\{\bU_{t_k}\}_{k=1}^n$; our study accommodates such a feature with an innovative device.
Denote by $\mathcal{F}_{-\infty}^{s}$ and $\mathcal{F}_s^\infty$ the $\sigma$-fields generated by $\{\bU_{t_k}\}_{k\leq s}$ and $\{\bU_{t_k}\}_{k\geq s}$, respectively,  
 the $\alpha$-mixing coefficients are defined as
\begin{align}\label{eq:alphacoeff}
\alpha_n(m)=\sup_{s}\sup_{A\in\mathcal{F}_{-\infty}^s,B\in\mathcal{F}_{s+m}^\infty}|\mathbb{P}(AB)-\mathbb{P}(A)\mathbb{P}(B)|\,,~~~~m\geq1\,.
\end{align}
Then $\{\bU_{t_k}\}_{k=1}^n$ is  an {\it $\alpha$-mixing sequence} if $\alpha_n(m)\rightarrow0$ as $m\rightarrow\infty$.
The notion of $\alpha$-mixing is a conventional foundation  for  broadly characterizing the serial dependence.
%
%
%
Among others, causal ARMA processes with continuous
innovation distributions are $\alpha$-mixing with exponentially decaying
$\alpha$-mixing coefficients,  so are stationary Markov chains satisfying
certain conditions;  see Section 2.6.1 of \cite{FanYao_2003}.
Stationary GARCH models with finite second
moments and continuous innovation distributions are also $\alpha$-mixing
with exponentially decaying $\alpha$-mixing coefficients;  see Proposition 12 of \cite{CarrascoChen_2002}. Under certain conditions, vector auto-regressive (VAR) processes, multivariate ARCH processes, and multivariate GARCH processes are all $\alpha$-mixing with exponentially decaying $\alpha$-mixing coefficients;  see \cite{HP09}, \cite{BFS11} and \cite{Wong2020}.

In \eqref{eq:alphacoeff},  we highlight the necessary inclusion of  $n$, the frequency related sample size, in the $\alpha$-mixing coefficient.        The reason is that   in a high-dimensional data setting,   $p$ is commonly specified as a function of the sample size $n$.  Such an intrinsic dependence    makes characterizing the serial dependence substantially  more challenging.  To handle it in our study, we impose the following assumption on  $\alpha_n(m)$ defined in \eqref{eq:alphacoeff}.

\begin{as}\label{as:alphmix}
There exist some universal constants $C_1>1$, $C_2>0$ and $\varphi>0$ such that $\alpha_n(m)\leq C_1\exp\{-C_2(L_n^{-1}m)^{\varphi}\}$ for any $m\geq 1$, where $L_n>0$ may diverge with $n$.
\end{as}

Assumption \ref{as:alphmix} is our new dedicated device for characterizing the serial dependence of $\{\bU_{t_k}\}_{k=1}^n$ in  the context of high-frequency  high-dimensional data.
Here $L_n$   is introduced as a parameter to handle the aforementioned challenge due to the high data dimensionality, together with the conventional $m$ as  in the $\alpha$-mixing settings for analyzing time series.
As a development of its own interests, the synthetic  device in Assumption \ref{as:alphmix} successfully integrates the considerations of high-frequency and high-dimensional data, where
the usual interpretation of the $\alpha$-mixing remains:  the between-observation dependence is still getting weaker when they are further away in the serial data,  as characterized by both $L_n$ and $m$.
Intuitively,  the rationale is that $L_n$, as a standalone parameter, may diverge together with the sampling frequency and data dimensionality in a synthetic manner.  %
Such a divergence reflects the nature of this more challenging problem due to relatively limited data information, in the sense that the serial dependence in the measurement errors will become stronger as $L_n$ increases.

More specifically, Assumption \ref{as:alphmix} does not require $\{\bU_{t_k}\}$ to be strictly stationary, and it includes several commonly used models for $\{\bU_{t_k}\}$ as special cases. For an independent sequence $\{\bU_{t_k}\}$, we can select $L_n=1/2$ and $\varphi=\infty$ in Assumption \ref{as:alphmix}. For an $L_n$-dependent sequence $\{\bU_{t_k}\}$, we can select $\varphi=\infty$ in Assumption \ref{as:alphmix}. If $\{\bU_{t_k}\}$ follows VAR model, multivariate ARCH model or multivariate GARCH model with certain conditions, we can select $L_n=\varphi=1$ in Assumption \ref{as:alphmix}.  We provide a concrete example here with  a diverging $L_n$. For each $i\in[r]$, let $Z_{i,t}$ satisfy the diffusion process ${\rm d}Z_{i,t}=\tilde{\mu}_i(Z_{i,t};\btheta_i)\,{\rm d}t+\tilde{\sigma}_i(Z_{i,t};\bgamma_i)\,{\rm d}{W}_{i,t}$, where ${W}_{i,t}$ is a univariate standard Brownian motion, $\tilde{\mu}_i(\cdot;\cdot)$ and $\tilde{\sigma}_i(\cdot;\cdot)$ are two functions of $Z_{i,t}$ with some parameters $\btheta_i$ and $\bgamma_i$, respectively. Write $\bZ_t=(Z_{1,t},\ldots,Z_{r,t})^\top$ with $r$ independent processes $Z_{1,t},\ldots,Z_{r,t}$. Letting $\bU_{t_k}=\bA\bZ_{k\delta}$ for some known loading matrix $\bA\in\mathbb{R}^{p\times r}$ and some $\delta>0$, we can select $L_n=\delta^{-1}$ and $\varphi=1$ when $\tilde{\mu}_i(\cdot;\cdot)$ and $\tilde{\sigma}_i(\cdot;\cdot)$ satisfy certain conditions\footnote{For each $i\in[r]$, Lemma 4 of \cite{AitMykland_2004} indicates that $\{Z_{i,k\delta}\}_{k\geq1}$ is a $\rho$-mixing process with $\rho$-mixing coefficient $\rho_i(m)\leq \exp(-c_im\delta)$ for any integer $m\geq0$, where $c_i>0$ is a constant depending on the properties of $\tilde{\mu}_i(\cdot;\cdot)$ and $\tilde{\sigma}_i(\cdot;\cdot)$ (see Assumption 1 of \cite{AitMykland_2004}). Theorem 5.1 of \cite{Bradley2005} implies $\{\bZ_{k\delta}\}_{k\geq1}$ is also a $\rho$-mixing process with $\rho$-mixing coefficient $\rho(m)\leq\exp(-c_{\min} m\delta)$ for any integer $m\geq0$, where $c_{\min}=\min_{i\in[r]}c_i$. Since $\rho$-mixing implies $\alpha$-mixing, then $\alpha_n(m)$ defined in \eqref{eq:alphacoeff} satisfies $\alpha_n(m)\leq 4^{-1}\rho(m)\leq4^{-1}\exp(-c_{\min}m\delta)$ for any integer $m\geq0$.}, where $L_n$ will diverge with $n$ if $\delta\rightarrow0$  as $n\rightarrow\infty$.  Here $L_n$ is also allowed to  depend directly on $p$, the dimension of $\bU_{t_k}$. As an example, if each univariate sequence $\{U_{i,t_k}\}_{k=1}^n$ is $\alpha$-mixing with exponentially decaying $\alpha$-mixing coefficients, with the independent assumption imposed on the $p$ sequences $\{U_{1,t_k}\}_{k=1}^n,\ldots,\{U_{p,t_k}\}_{k=1}^n$,  Theorem 5.1 of \cite{Bradley2005} indicates that $\alpha_n(m)$ defined in \eqref{eq:alphacoeff} satisfies $\alpha_n(m)\leq p\exp(-cm)$ for some universal constant $c>0$, which implies Assumption \ref{as:alphmix} holds for $\varphi=1$ and $L_n\asymp \log p$.

 To our best knowledge, there is no alternative assumption in the literature that is capable of handling the setting of our study.
In existing studies, some serial dependence assumptions have been imposed  on the measurement errors, with a primary objective recovering its auto-covariance.
When $p=1$, \cite{Jacod2017} assumes $U_{1,t_{k}}=\gamma_{t_{k}}\chi_k$ for some nonnegative semimartingale $\gamma_t$ and a $\rho$-mixing stationary sequence $\{\chi_k\}_{k\geq1}$, where $\{\chi_k\}_{k\geq1}$ is independent of the process $\gamma_t$;  see also  the setting of \cite{LiLinton} that covers serially dependent, endogenous,
and nonstationary noises.
If $\gamma_t$ is the solution of some stochastic differential equations, $\{\gamma_{t_{k}}\}_{k\geq1}$ is also $\rho$-mixing. See, for example, Lemma 4 of \cite{AitMykland_2004}. Based on the independence between $\{\chi_k\}_{k\geq1}$ and $\{\gamma_{t_{k}}\}_{k\geq1}$, Theorem 5.2 of \cite{Bradley2005} implies the sequence $\{U_{1,t_{k}}\}_{k\geq1}$ is also $\rho$-mixing. Since $\rho$-mixing implies $\alpha$-mixing, we know $\{U_{1,t_{k}}\}_{k\geq1}$ is also $\alpha$-mixing.  \cite{Varneskov2017} relaxes the $\rho$-mixing assumption on $\{\chi_k\}_{k\geq1}$ to the weaker $\alpha$-mixing condition.  In a recent study, \cite{DaXiu2021} assume instead a working moving average structure for the measurement errors.




For \eqref{eq:mod},
we assume ${\rm Cov}(\bU_{t_k})\equiv\bSigma_u$ for each $k\in[n]$.  Our main goal in this study  is to estimate  $\bSigma_u$, the covariance matrix that contains information on the between-component relationship of the unobserved noise $\bU_{t_k}$. Clearly, $\bU_{t_k}$ is a latent vector.
To estimate its covariance matrix,  eliminating the impact due to the process $\bX_t$ is required, which means that now $\bU_{t_k}$ performs like `signal' and $\bX_{t_k}$ is `noise'.  Our strategy is to perform a dedicated localization:  focusing on observations that are in a specific neighborhood mentioned later.
 For any $i,j\in[p]$, we write $\mathcal{G}_i\cap\mathcal{G}_j=\{t_{i,j,1},\ldots,t_{i,j,n_{i,j}}\}$ with $t_{i,j,1}<\cdots<t_{i,j,n_{i,j}}$.
Let $\Delta t_{i,j,k}=t_{i,j,k+1}-t_{i,j,k}$ for
any $k\in[n_{i,j}-1]$. In this
paper,  
we consider the scenario with $T$ being fixed but $\max_{i,j\in[p]}\max_{k\in[n_{i,j}-1]}\Delta t_{i,j,k}\rightarrow0$ as $n\rightarrow\infty$. Formally, we make the following assumption:

\begin{as}\label{as:space}
(i) As $n\rightarrow\infty$, $\min_{i,j\in[p]}\min_{k\in[n_{i,j}-1]}\Delta t_{i,j,k}/\max_{i,j\in[p]}\max_{k\in[n_{i,j}-1]}\Delta t_{i,j,k}$ is uniformly bounded away from zero. (ii) As $n\rightarrow\infty$, we have each $n_{i,j}\rightarrow\infty$, and $\min_{i,j\in[p]}n_{i,j}/\max_{i,j\in[p]}n_{i,j}$ is uniformly bounded away from zero. (iii) $\min_{i,j\in[p]}(t_{i,j,n_{i,j}}-t_{i,j,1})\asymp T$.
\end{as}

The setting with Assumption \ref{as:space} is broad and general. The first part is a  standard setting for  studying high-frequency data.
The second part requires enough number of pairwise synchronous observations.   This is a reasonable practical setting; see also \cite{Aitetal_2010_JASA} for a pairwise approach for estimating the realized covariance matrix for $(\bX_t)_{t\in[0,T]}$.
Based on part (ii) of Assumption \ref{as:space}, we write
\begin{equation}\label{eq:nstar}
\min_{i,j\in[p]}n_{i,j}\asymp\max_{i,j\in[p]}n_{i,j}\asymp n_*\,,
\end{equation}
where $n_*\rightarrow\infty$ as $n\rightarrow\infty$.
As we will show in Theorems \ref{pn:1}--\ref{tm:4}, the convergence rates for the estimates of the covariance matrix $\bSigma_{u}=\cov(\bU_{t_k})$ will depend on $n_*$ instead of $n$.
In  the special case with synchronous observations, we have $n_{i,j}=n$ for any $i,j\in[n]$ and we can set $n_*=n$.  Then all our results also apply to the setting with synchronous data. Assumption \ref{as:space} is not necessary for our theoretical analysis which is just imposed for simplicity and can be removed at the expenses of lengthier proofs. Our theoretical analysis essentially only requires the assumption that $\min_{i,j\in[p]}n_{i,j}\rightarrow\infty$ and $\max_{i,j\in[p]}\max_{k\in[n_{i,j}-1]}\Delta t_{i,j,k}\rightarrow0$ as $n\rightarrow\infty$. With such assumption, both $\min_{i,j\in[p]}\min_{k\in[n_{i,j}-1]}\Delta t_{i,j,k}/\max_{i,j\in[p]}\max_{k\in[n_{i,j}-1]}\Delta t_{i,j,k}$ and $\min_{i,j\in[p]}n_{i,j}/\max_{i,j\in[p]}n_{i,j}$ can decay to zero as $n\rightarrow\infty$. We will discuss in Section \ref{s5} how this assumption affects the convergence rates for the estimates of the covariance matrix $\bSigma_{u}$.

\subsection{Covariance matrix estimation of $\bU_{t_k}$}\label{sec:methfU}

 Write $\bSigma_{u}=(\sigma_{u,i,j})_{p\times p}$ and $\bU_{t_k}=(U_{1,t_k},\ldots,U_{p,t_k})^{\top}$. Here the subscript $u$ in $\sigma_{u,i,j}$ indicates that it is a quantity associated with the noise so as to differentiate it from the volatility process $\sigma_{i,t}$ in \eqref{eq:model}. We know $2\bSigma_{u}=\cov(\bU_{t_k}-\bU_{t_\ell})+\mathbb{E}(\bU_{t_\ell}\bU_{t_k}^{\top})+\mathbb{E}(\bU_{t_k}\bU_{t_\ell}^{\top})$ for any $\ell\neq k$. By Assumption \ref{as:alphmix} and Davydov's inequality \citep{Davydov1968}, $|\mathbb{E}(\bU_{t_\ell}\bU_{t_k}^{\top})|_\infty+|\mathbb{E}(\bU_{t_k}\bU_{t_\ell}^{\top})|_\infty\lesssim \exp(-C_*L_n^{-\varphi}|k-\ell|^{\varphi})$ for some constant $C_*>0$, provided that $\max_{j\in[p]}\mathbb{E}(|U_{j,t_k}|^\gamma)$ and $\max_{j\in[p]}\mathbb{E}(|U_{j,t_\ell}|^\gamma)$ are uniformly bounded away from infinity for some universal constant $\gamma>2$. Notice that $\bU_{t_k}-\bU_{t_\ell}=(\bY_{t_k}-\bY_{t_\ell})-(\bX_{t_k}-\bX_{t_\ell})$ and each component process $X_{i,t}$ is a continuous-time and continuous-path
stochastic process.  We have considerations from two ends.  First, due to $|X_{i,t+h}-X_{i,t}|\rightarrow0$ almost surely as $h\rightarrow0$, in a small neighborhood $\cal N$  of $t_k$,  the difference between the high-frequency observations $\bY_{t_k}$ and $\bY_{t_\ell}$, for $t_\ell\in \cal N$,  can be approximately viewed as $\bU_{t_k}-\bU_{t_\ell}$. Second, to avoid excessive impact from aggregating $\mathbb{E}(\bU_{t_\ell}\bU_{t_k}^{\top})+\mathbb{E}(\bU_{t_k}\bU_{t_\ell}^{\top})$, we cannot choose $t_\ell$ and $t_k$ too close.
Putting these two considerations together, 
we propose to estimate $\sigma_{u,i,j}$ by
\begin{equation}\label{eq:estasyn}
\hat{\sigma}_{u,i, j}=\frac{1}{2n_{i, j}}\sum_{k=1}^{n_{i, j}}\frac{1}{N_{i ,j , k}}\sum_{t_{i,j,\ell}\in S_{i,j,k}}(Y_{i,t_{i, j, \ell}}-Y_{i, t_{i, j, k}})(Y_{j,t_{i, j, \ell}}-Y_{j, t_{i, j, k}})\,,
\end{equation}
where $S_{i,j,k}=\{t_{i, j, \ell}\in\mathcal{G}_i\cap\mathcal{G}_j:K\leq|\ell-k|\leq K+\Delta_K \}$ for some integers $K\geq1$ and $\Delta_K\geq0$, and $N_{i ,j , k}=|S_{i,j,k}|$.
Here the set  $S_{i,j,k}$ is designed to meet the aforementioned two considerations -- ensuring data in an appropriate range are incorporated for estimating $\bSigma_{u}$.
%
%


Our estimator \eqref{eq:estasyn} with $S_{i,j,k}$  is generally applicable.
For an independent sequence $\{\bU_{t_k}\}$, we can select $K=1$ and then $\mathbb{E}(U_{i,t_{i,j,\ell}}U_{j,t_{i,j,k}})=\mathbb{E}(U_{i,t_{i,j,k}}U_{j,t_{i,j,\ell}})=0$ for any $t_{i,j,\ell}\in S_{i,j,k}$ due to $(L_n,\varphi)=(1/2,\infty)$ in Assumption \ref{as:alphmix}. For an $L_n$-dependent sequence $\{\bU_{t_k}\}$, we can select $K>L_n$ due to $\varphi=\infty$  in Assumption \ref{as:alphmix} and then $\mathbb{E}(U_{i,t_{i,j,\ell}}U_{j,t_{i,j,k}})=\mathbb{E}(U_{i,t_{i,j,k}}U_{j,t_{i,j,\ell}})=0$ for any $t_{i,j,\ell}\in S_{i,j,k}$. For general case with $\varphi<\infty$, with selecting $K\geq L_n(C_{**}\log n_*)^{1/\varphi}$ for some sufficiently large constant $C_{**}>0$, $\max_{i,j\in[p]}\max_{k\in[n_{i,j}]}\max_{t_{i,j,\ell}\in S_{i,j,k}}|\mathbb{E}(U_{i,t_{i,j,\ell}}U_{j,t_{i,j,k}})+\mathbb{E}(U_{i,t_{i,j,k}}U_{j,t_{i,j,\ell}})|\lesssim n_*^{-C_*C_{**}}$, which is negligible in comparison to the bias from approximating $(U_{i,t_{i,j,\ell}}-U_{i,t_{i,j,k}},U_{j,t_{i,j,\ell}}-U_{j,t_{i,j,k}})$ by $(Y_{i,t_{i,j,\ell}}-Y_{i,t_{i,j,k}},Y_{j,t_{i,j,\ell}}-Y_{j,t_{i,j,k}})$. 
To simplify our presentation, we assume $\Delta_K\geq0$ is a fixed integer in this paper. Our theoretical results  can be  parallel extended to the scenario with  diverging $\Delta_K$.

For a fixed $T$,  (\ref{eq:nstar}) and Assumption~\ref{as:space} imply that \[
\min_{i,j\in[p]}\min_{k\in[n_{i,j}-1]}\Delta t_{i,j,k}\asymp\max_{i,j\in[p]}\max_{k\in[n_{i,j}-1]}\Delta t_{i,j,k}\asymp n_*^{-1}\,.
\] Let \begin{equation}\label{eq:hatsigma}
\widehat{\bSigma}_{u}=(\hat{\sigma}_{u,i,j})_{p\times p}
\end{equation} for $\hat{\sigma}_{u,i,j}$ defined as (\ref{eq:estasyn}). Theorem \ref{pn:1} in Section \ref{sec:theoryfU} shows that the elements of $\widehat{\bSigma}_{u}$ are uniformly consistent to the corresponding elements of $\bSigma_{u}$ with a suitable selection of $K$, i.e.
\[
\mathbb{E}\big(|\widehat{\bSigma}_{u}-\bSigma_{u}|_\infty\big)\lesssim(Kn_*^{-1}\log p)^{1/2}\,.
\]
Theorem \ref{tm:2} in Section \ref{sec:theoryfU} shows that $(n_*^{-1}\log p)^{1/2}$ is the minimax optimal rate in the maximum element-wise loss for the covariance matrix estimations of the high-dimensional noise $\bU_{t_k}$ in high-frequency data. If $\{\bU_{t_k}\}$ is an independent or $L_n$-dependent sequence with fixed $L_n$, we can select $K$ as a fixed integer and then the associated convergence rate of $|\widehat{\bSigma}_u-\bSigma_u|_\infty$ is minimax optimal. For general cases with $\varphi<\infty$ and fixed $L_n$, with selecting $K\asymp (\log n_*)^{1/\varphi+\epsilon}$ for some $\epsilon>0$, the convergence rate of $|\widehat{\bSigma}_u-\bSigma_u|_\infty$ is nearly optimal with an additional logarithm factor $(\log n_*)^{1/(2\varphi)+\epsilon/2}$. More importantly, as we will discuss below Remark \ref{rek3} in Section \ref{sec:theoryfU}, $(n_*^{-1}\log p)^{1/2}$ is also the minimax optimal rate in the maximum element-wise loss for the covariance matrix estimations of $\bU_{t_k}$ if we have observations of the noise, which indicates that our estimator shares some oracle property and the proposed localization actually makes the impact of the latent process $\bX_t$ be negligible.

However,  the aforementioned element-wise consistency and optimality do not imply their counterparts for the covariance matrix estimation with high-dimensional data.
That is, the estimator
$\widehat{\bSigma}_{u}$ may not be consistent to $\bSigma_{u}$ under the spectral norm $\|\cdot\|_2$ when $p\gg n$. This is a well-known phenomenon in high-dimensional covariance matrix estimation;  see, among other, \cite{BickelLevina_2008_AOS}.
For high-dimensional covariance matrix estimations,   one often resorts to some classes of the target with extra information. With the extra information, the consistency under the spectral norm and other properties associated with the covariance matrix estimations can be well established.
 In this paper, we focus on the following class -- the sparse covariance matrices considered in \cite{BickelLevina_2008_thre}:
\begin{equation}\label{eq:sparse-set}
\mathcal{H}(q,c_p,M)=\bigg\{\bSigma_u=(\sigma_{u,i,j})_{p\times p}: \sigma_{u,i,i}\leq M~\textrm{and}~\sum_{j=1}^p|\sigma_{u,i,j}|^q\leq c_p~\textrm{for all}~i\bigg\}\,,
\end{equation}
where $q\in[0,1)$ and $M>0$ are two prescribed constants, and $c_p$ may diverge with $p$. Here $c_p$ can be viewed as a parameter that characterizes the sparsity of $\bSigma_u$, i.e., if $c_p$ is smaller, then $\bSigma_u$ is more sparse. If $q=0$, we have
\[
\mathcal{H}(0,c_p,M)=\bigg\{\bSigma_u=(\sigma_{u,i,j})_{p\times p}:\sigma_{u,i,i}\leq M~\textrm{and}~\sum_{j=1}^pI(\sigma_{u,i,j}\neq 0)\leq c_p~\textrm{for all}~i\bigg\}\,,
\]
where $c_p$ evaluates the number of nonzero components in each row of $\bSigma_u$.

For $\bSigma_{u}\in\mathcal{H}(q,c_p,M)$, we  propose  the following thresholding estimator based on the element-wise estimation $\widehat{\bSigma}_{u}$ given in \eqref{eq:hatsigma}:
\begin{equation} \label{eq:est-thre}
\widehat{\bSigma}_{u}^{{\rm thre}}=\big[\hat{\sigma}_{u,i,j}{I}\big\{|\hat{\sigma}_{u,i,j}|\geq \beta(Kn_*^{-1}\log p)^{1/2}\big\}\big]_{p\times p}\,,
\end{equation}
where $\beta>0$ is a fixed constant for the thresholding level. Theorem \ref{tm:3} in Section \ref{sec:theoryfU} shows that such defined thresholding estimator $\widehat{\bSigma}_{u}^{{\rm thre}}$ is consistent to $\bSigma_{u}$ under the spectral norm with suitable selections of $K$ and $\beta$, i.e.
 \begin{align}\label{eq:convrate}
\mathbb{E}\big(\|\widehat{\bSigma}_{u}^{{\rm thre}}-\bSigma_{u}\|_2^2\big)\lesssim c_p^2(Kn_*^{-1}\log p)^{1-q}\,.
 \end{align}
 Furthermore, Theorem \ref{tm:4} in Section \ref{sec:theoryfU} indicates that $c_p(n_*^{-1}\log p)^{(1-q)/2}$ is the minimax optimal convergence rate with the spectral norm loss function for the covariance matrix estimations of the high-dimensional noise $\bU_{t_k}$ in high-frequency data, which is also the minimax optimal convergence rate in the spectral norm loss if we have observations of the noise directly. If $\{\bU_{t_k}\}$ is an independent or $L_n$-dependent sequence with fixed $L_n$, we can select $K$ as a fixed integer and then the associated convergence rate of $\|\widehat{\bSigma}_u^{{\rm thre}}-\bSigma_u\|_2$ is minimax optimal. For general cases with $\varphi<\infty$ and fixed $L_n$, with selecting $K\asymp (\log n_*)^{1/\varphi+\epsilon}$ for some $\epsilon>0$, the convergence rate of $\|\widehat{\bSigma}_u^{{\rm thre}}-\bSigma_u\|_2$ is nearly optimal with an additional logarithm factor $(\log n_*)^{(1/\varphi+\epsilon)(1-q)/2}$.

\begin{rek}
In finite samples, the thresholding estimator $\widehat{\bSigma}_u^{{\rm thre}}$ given in \eqref{eq:est-thre} may not be positive definite in general. We can first apply the singular value decomposition to $\widehat{{\bSigma}}_u^{{\rm thre}}$: $\widehat{\bSigma}_u^{{\rm thre}}=\widehat{{\bf P}}^\top{\rm diag}(\hat{\tau}_1,\ldots,\hat{\tau}_p)\widehat{{\bf P}}$, where $\hat{\tau}_1\geq\cdots\geq\hat{\tau}_p$ are the eigenvalues of $\widehat{\bSigma}_u^{{\rm thre}}$, and $\widehat{{\bf P}}$ is an orthogonal matrix. If there are $s$ negative eigenvalues, we can use $\widetilde{\bSigma}_u^{{\rm thre}}=\widehat{{\bf P}}^\top{\rm diag}(\hat{\tau}_1,\ldots,\hat{\tau}_{p-s},\hat{\tau}_{p-s+1}+\epsilon,\ldots,\hat{\tau}_p+\epsilon)\widehat{{\bf P}}$ as the estimate of $\bSigma_u$ for some $\epsilon>0$. Write $\delta_n=c_p(Kn_*^{-1}\log p)^{(1-q)/2}$ and let $\tau_1\geq\cdots\geq\tau_p>0$ be the eigenvalues of $\bSigma_u$. Since $\max_{j\in[p]}|\hat{\tau}_j-\tau_j|\leq\|\widehat{\bSigma}_u^{{\rm thre}}-{\bf \Sigma}_u\|_2=O_{\p}(\delta_n)$, if $\tau_p$ is uniformly bounded away from zero and we select $\epsilon=-\hat{\tau}_p+(Kn_*^{-1}\log p)^{1/2}$ when $\hat{\tau}_p<0$, such defined $\widetilde{\bSigma}_u^{{\rm thre}}$ is positive definite and also satisfies \eqref{eq:convrate}.
\end{rek}

\section{Theoretical analysis}\label{sec:theoryfU}
In this section,  we establish the theoretical properties of the proposed estimators.  To mimic the high-dimensional scenario, we always assume $p\geq n_*^{\kappa}$ for some universal constant $\kappa>0$ in  this paper. We also  require the following  assumptions.

\begin{as}\label{as:moment}
Write $\bU_{t_k}=(U_{1,t_k},\ldots,U_{p,t_k})^{\top}$. There exist some universal constants $C_3>1$ and $C_4>0$ such that $\mathbb{P}(|U_{i,t_k}|>u)\leq C_3\exp(-C_4u^{2})$ for any $i\in[p]$, $k\in[n]$ and $u>0$.
\end{as}

\begin{as}\label{as:drifdiff}
There exist some universal  constants $C_5>0$, $C_6>0$ and $C_7>0$ such that
(i) $\mathbb{E}(\exp[\theta\{\mu_{i,t}^2-\mathbb{E}(\mu_{i,t}^2)\}])\leq \exp(C_6\theta^2)$ and
$\mathbb{E}(\exp[\theta\{\sigma_{i,t}^2-\mathbb{E}(\sigma_{i,t}^2)\}])\leq \exp(C_6\theta^2)
$ for any $i\in[p]$, $t\in[0,T]$ and $|\theta|\leq C_5$; (ii) $\mathbb{E}(\mu_{i,t}^2)\leq C_7
$ and $\mathbb{E}(\sigma_{i,t}^2)\leq C_7$ for any $i\in[p]$ and $t\in[0,T]$.
\end{as}

\begin{as}\label{as:diff}
There exist some universal constants $\gamma>0$, $C_8>1$ and $C_9>0$ such that
$
\mathbb{P}(\sup_{0\leq t\leq T}\sigma_{i,t}>u)\leq C_8\exp(-C_9u^\gamma)
$
for any $i\in[p]$ and $u>0$.
\end{as}

All assumptions are mild for studying high-dimensional covariance matrix estimations with high-frequency data. Assumption \ref{as:moment} requires that each component of $\bU_{t_k}$ is sub-Gaussian. Following Lemma 2.2 of \cite{Petrov_1995}, we know that part (i) of Assumption \ref{as:drifdiff} holds if there exist two positive constants $C_*$ and $C_{**}$ such that $\mathbb{P}\{|\mu_{i,t}^2-\mathbb{E}(\mu_{i,t}^2)|\geq u\}\leq C_*\exp(-C_{**}u)$ and $\mathbb{P}\{|\sigma_{i,t}^2-\mathbb{E}(\sigma_{i,t}^2)|\geq u\}\leq C_*\exp(-C_{**}u)$ for any $i\in[p]$, $t\in [0,T]$ and $u>0$. Assumption \ref{as:diff} describes the behavior of the tail probability of $\sup_{0\leq t\leq T}\sigma_{i,t}$. If the spot volatility process $\sigma_{i,t}$ is uniformly bounded away from infinity over $i\in[p]$ and $t\in[0,T]$, we can select $\gamma=\infty$ in Assumption \ref{as:diff}. Then we have the following result.

\begin{theorem}\label{pn:1}
Let $\mathcal{P}_1$ denote the collections of models for $\{\bY_{t_k}\}_{k=1}^n$ such that $\bY_{t_k}=\bX_{t_k}+\bU_{t_k}$, where the noises $\{\bU_{t_k}\}_{k=1}^n$ satisfy Assumption {\rm\ref{as:moment}}, $\bX_t=(X_{1,t},\ldots,X_{p,t})^{\top}$ follows 
 model {\rm(\ref{eq:model})} with each $\mu_{i,t}$ and $\sigma_{i,t}$ satisfying Assumptions {\rm\ref{as:drifdiff}} and {\rm\ref{as:diff}}, and the grids of time points $\{\mathcal{G}_i\}_{i=1}^p$ satisfy Assumption {\rm\ref{as:space}}. Let $K> CL_n$ for some constant $C\geq 1$. Under Assumption {\rm \ref{as:alphmix}}, it holds that
$$
\sup_{\mathcal{P}_1}\mathbb{E}\big(|\widehat{\bSigma}_{u}-\bSigma_{u}|_\infty\big)\lesssim(Kn_*^{-1}\log p)^{1/2}$$
provided that $\log p=o[\min\{(n_*L_n^{-2}K)^{\varphi/(3\varphi+2)},(n_*K^{-1})^\chi\}]$ and  $K^{-\varphi}L_n^{\varphi}\log\{n_*(K\log p)^{-1}\}=o(1)$, where $n_*$ is specified in {\rm(\ref{eq:nstar})} and $\chi=\min\{\gamma/(\gamma+4),1/3\}$.
\end{theorem}

\begin{rek}\label{remark1}
Theorem \ref{pn:1} gives the convergence rate of $\max_{i,j\in[p]}|\hat{\sigma}_{u,i,j}-\sigma_{u,i,j}|$.

(i) For an independent sequence $\{\bU_{t_k}\}$, due to $(L_n,\varphi)=(1/2,\infty)$ and $K\geq1$, $\mathbb{E}(|\widehat{\bSigma}_{u}-\bSigma_{u}|_\infty)\lesssim (Kn_*^{-1}\log p)^{1/2}$ provided that $\log p=o\{(n_*K^{-1})^{\chi}\}$. For fixed $K$, $\max_{i,j\in[p]}|\hat{\sigma}_{u,i,j}-\sigma_{u,i,j}|=O_\p\{(n_*^{-1}\log p)^{1/2}\}$. 

(ii) For an $L_n$-dependent sequence $\{\bU_{t_k}\}$, due to $\varphi=\infty$ and $K>L_n$,  $\mathbb{E}(|\widehat{\bSigma}_{u}-\bSigma_{u}|_\infty)\lesssim (Kn_*^{-1}\log p)^{1/2}$ provided that $\log p=o\{(n_*K^{-1})^{\chi}\}$. For fixed $L_n$, $\max_{i,j\in[p]}|\hat{\sigma}_{u,i,j}-\sigma_{u,i,j}|=O_\p\{(n_*^{-1}\log p)^{1/2}\}$ with selecting a fixed $K>L_n$. For diverging $L_n$, $\max_{i,j\in[p]}|\hat{\sigma}_{u,i,j}-\sigma_{u,i,j}|=O_\p\{(L_nn_*^{-1}\log p)^{1/2}\}$ with selecting $K\asymp L_n$ and $K>L_n$.

(iii) For the general cases with $\varphi<\infty$ and fixed $L_n$, to make $K^{-\varphi}L_n^{\varphi}\log\{n_*(K\log p)^{-1}\}=o(1)$, we need to select $K\gg (\log n_*)^{1/\varphi}$. If we select $K\asymp (\log n_*)^{1/\varphi+\epsilon}$ for some $\epsilon>0$, $\max_{i,j\in[p]}|\hat{\sigma}_{u,i,j}-\sigma_{u,i,j}|=O_\p\{(n_*^{-1}\log p)^{1/2}(\log n_*)^{1/(2\varphi)+\epsilon/2}\}$.

\end{rek}

 Furthermore, Theorem \ref{tm:2} below shows that the convergence rate $(n_*^{-1}\log p)^{1/2}$ is minimax optimal in the maximum element-wise loss for the covariance matrix estimations of the high-dimensional noise $\bU_{t_k}$ in high-frequency data.

\begin{theorem}\label{tm:2}
Let $n/n_*\lesssim p$. Denote by $\check{\mathcal{F}}$ the class of all measurable functionals of the data. Then
$$
\inf_{\widehat{\bSigma}\in\check{\mathcal{F}}}\sup_{\mathcal{P}_1}\mathbb{E}\big(|\widehat{\bSigma}-\bSigma_{u}|_\infty\big)\gtrsim(n_*^{-1}\log p)^{1/2}\,,$$
where $\mathcal{P}_1$ is defined in Theorem {\rm \ref{pn:1}}.
\end{theorem}

\begin{rek}\label{rek3}
(i) If $\{\bU_{t_k}\}$ is an independent sequence or $L_n$-dependent sequence with fixed $L_n$, Remarks \ref{remark1}(i) and \ref{remark1}(ii) indicate that our proposed estimate $\widehat{\bSigma}_u$ is minimax optimal under the maximum element-wise loss.

(ii) For the general cases with $\varphi<\infty$ and fixed $L_n$, Remark \ref{remark1}(iii) indicates our proposed estimate $\widehat{\bSigma}_u$ is nearly minimax optimal under the maximum element-wise loss with an additional logarithm factor $(\log n_*)^{1/(2\varphi)+\epsilon/2}$.

\end{rek}

To establish the lower bound stated in Theorem \ref{tm:2}, we essentially focus on a model belonging to $\mathcal{P}_1$ with $\mu_{i,t}=0$ and $\sigma_{i,t}=0$ for any $t\in[0,T]$ and $i\in[p]$. Let $\mathcal{C}_k=\{i\in[p]:t_k\in\mathcal{G}_i\}$ for any $k\in[n]$. In this specific model, the latent process $\bX_t=\bzero$ for any $t\in[0,T]$ and thus the data we observed are $\mathcal{Z}=\{\bU_{t_1,\mathcal{C}_1},\ldots,\bU_{t_n,\mathcal{C}_n}\}$. Here $\bU_{t_k,\mathcal{C}_k}$ denotes the subvector of $\bU_{t_k}$ with components indexed by $\mathcal{C}_{k}$. Hence, $(n_*^{-1}\log p)^{1/2}$ is also the minimax optimal rate in the maximum element-wise loss for the covariance matrix estimations of $\bU_{t_k}$ with data $\mathcal{Z}=\{\bU_{t_1,\mathcal{C}_1},\ldots,\bU_{t_n,\mathcal{C}_n}\}$, which indicates that the estimator $\widehat{\bSigma}_{u}$
shares some oracle property and the proposed localization actually makes the impact of the
latent process $\bX_t$ be negligible.

Regarding the loss function  under the spectral norm $\|\cdot\|_2$ for the whole covariance matrix estimation,  Theorem \ref{tm:3} establishes the convergence rate of the thresholding estimator $\widehat{\bSigma}_{u}^{{\rm thre}}$ defined as (\ref{eq:est-thre}).

\begin{theorem}\label{tm:3}
Let $\mathcal{P}_2$ denote the collections of models for $\{\bY_{t_k}\}_{k=1}^n$ such that $\bY_{t_k}=\bX_{t_k}+\bU_{t_k}$, where the noises $\{\bU_{t_k}\}_{k=1}^n$ satisfy Assumption {\rm\ref{as:moment}} with the covariance matrix $\bSigma_{u}\in\mathcal{H}(q,c_p,M)$, $\bX_t=(X_{1,t},\ldots,X_{p,t})^{\top}$ follows model {\rm(\ref{eq:model})} with each $\mu_{i,t}$ and $\sigma_{i,t}$ satisfying Assumptions {\rm\ref{as:drifdiff}} and {\rm\ref{as:diff}}, and the grids of time points $\{\mathcal{G}_i\}_{i=1}^p$ satisfy Assumption {\rm\ref{as:space}}. Let $K> CL_n$ for some constant $C\geq 1$. Under Assumption {\rm \ref{as:alphmix}}, with sufficiently large constant $\beta>0$ in \eqref{eq:est-thre}, it holds that
$$
\sup_{\mathcal{P}_2}\mathbb{E}\big(\|\widehat{\bSigma}_{u}^{{\rm thre}}-\bSigma_{u}\|_2^2\big)\lesssim c_p^2(Kn_*^{-1}\log p)^{1-q}$$
provided that $\log p=o[\min\{(n_*L_n^{-2}K)^{\varphi/(3\varphi+2)},(n_*K^{-1})^{\chi}\}]$ and $K^{-\varphi}L_n^{\varphi}\log\{n_*(K\log p)^{-1}\}=o(1)$, where $n_*$ is specified in {\rm(\ref{eq:nstar})} and $\chi=\min\{\gamma/(\gamma+4),1/3\}$.
\end{theorem}

Our result in  the following Theorem \ref{tm:4} justifies that the convergence rate $c_p(n_*^{-1}\log p)^{(1-q)/2}$ is minimax optimal  under the spectral norm loss  function for the covariance matrix estimations of $\bU_{t_k}$ with the sparsity structure (\ref{eq:sparse-set}). Again, this rate is also the minimax optimal
 rate in the spectral norm loss for the covariance matrix estimations of $\bU_{t_k}$ with data $\mathcal{Z}=\{\bU_{t_1,\mathcal{C}_1},\ldots,\bU_{t_n,\mathcal{C}_n}\}$.

\begin{theorem}\label{tm:4}
Let $n/n_*\lesssim p$. Denote by $\check{\mathcal{F}}$ the class of all measurable functionals of the data. Then
$$
\inf_{\widehat{\bSigma}\in\check{\mathcal{F}}}\sup_{\mathcal{P}_2}\mathbb{E}\big(\|\widehat{\bSigma}-\bSigma_{u}\|_2^2\big)\gtrsim c_p^2(n_*^{-1}\log p)^{1-q}$$
provided that $c_p\lesssim n_*^{(1-q)/2}(\log p)^{-(3-q)/2}$, where $\mathcal{P}_2$ is defined in Theorem {\rm \ref{tm:3}}.
\end{theorem}

\begin{rek}
(i) If $\{\bU_{t_k}\}$ is an independent sequence or $L_n$-dependent sequence with fixed $L_n$, Theorems \ref{tm:3} and \ref{tm:4} indicate that $\widehat{\bSigma}_u^{{\rm thre}}$ defined as (\ref{eq:est-thre}) is minimax optimal under the spectral norm loss when we select $K$ as a fixed integer.

(ii) For the general cases with $\varphi<\infty$ and fixed $L_n$, if we select $K\asymp (\log n_*)^{1/\varphi+\epsilon}$ for some $\epsilon>0$, Theorems \ref{tm:3} and \ref{tm:4} indicate that $\widehat{\bSigma}_u^{{\rm thre}}$ defined as (\ref{eq:est-thre}) is nearly minimax optimal under the spectral norm with an additional logarithm factor $(\log n_*)^{(1/\varphi+\epsilon)(1-q)/2}$.

\end{rek}
In summary, we conclude that it is $n_*$ -- the effective sample size of the pairwise synchronous observations --  determining the convergence rate of the covariance matrix estimation of the noise $\bU_{t_k}$. Practically, $n_*$ is  expected to  be smaller than $n$ -- the total number of observation times.   Hence, the accuracy of the covariance matrix estimation  is affected by the level of data asynchronicity --  the more asynchronous the data are, the more difficult it is to estimate $\bSigma_{u}$. Another finding from our theoretical analysis is that although the noise $\{\bU_{t_1,\mathcal{C}_1},\ldots,\bU_{t_n,\mathcal{C}_n}\}$ are not directly observable,  the localized estimator in some scenarios has the (nearly) same accuracy as the one when the noise $\{\bU_{t_1,\mathcal{C}_1},\ldots,\bU_{t_n,\mathcal{C}_n}\}$ are observed in the sense of the (nearly) same convergence rates for estimating $\bSigma_{u}$ with high-frequency data. From the practical perspective, it can be viewed as a bless from the high-frequency data with adequate amount of data information locally, so that the statistical properties of the noise can be  accurately revealed.

\section{The effect of the smallness of the noise}\label{sec:4}
Our results in Section \ref{sec:theoryfU} assume that $n_{i,j}\to \infty$ with $T$ fixed.
Empirically, as pointed out in \cite{HL06}, the magnitude of $\sigma_{u,i,j}$ may be small; see also \cite{Christensen2014}.
To address this issue,  we study the second-order property of our estimator concerning its bias.

For $\hat{\sigma}_{u,i,j}$ defined as (\ref{eq:estasyn}), since $\{\bU_{t_k}\}_{k=1}^n$ is independent of $(\bX_t)_{0\leq t\leq T}$, we have that
\begin{align}\label{eq:highexp}
&\mathbb{E}\{\hat{\sigma}_{u,i,j}\,|\,(\bX_t)_{0\leq t\leq T}\}-\sigma_{u,i,j}\notag\\
&~~~~~~=-\underbrace{\frac{1}{2n_{i,j}}\sum_{k=1}^{n_{i,j}}\mathbb{E}\bigg( \frac{U_{i,t_{i,j,k}}}{N_{i,j,k}}\sum_{t_{i,j,\ell}\in S_{i,j,k}}U_{j,t_{i,j,\ell}}\bigg)}_{\textrm{I}'_2(i,j)}-\underbrace{\frac{1}{2n_{i,j}}\sum_{k=1}^{n_{i,j}}\mathbb{E}\bigg(\frac{U_{j,t_{i,j,k}}}{N_{i,j,k}}\sum_{t_{i,j,\ell}\in S_{i,j,k}}U_{i,t_{i,j,\ell}}\bigg)}_{\textrm{I}'_3(i,j)}\notag\\
&~~~~~~~~~~+\underbrace{\frac{1}{2n_{i,j}}\sum_{k=1}^{n_{i,j}}\frac{1}{N_{i,j,k}}\sum_{t_{i,j,\ell}\in S_{i,j,k}}(X_{i,t_{i,j,\ell}}-X_{i,t_{i,j,k}})(X_{j,t_{i,j,\ell}}-X_{j,t_{i,j,k}})}_{{\rm II}(i,j)}\,,
\end{align}
which indicates that the bias in $\hat{\sigma}_{u,i,j}$ includes three parts: ${\rm I}'_2(i,j)$, ${\rm I}'_3(i,j)$ and ${\rm II}(i,j)$. Proposition \ref{la:2} in Section \ref{proof} shows that $\max_{i,j\in[p]}|{\rm II}(i,j)|=O_{\p}\{(Kn_*^{-1}\log p)^{1/2}\}$, but $\mathbb{E}\{{\rm II}(i,j)\}\neq 0$, causing a bias of order $O(Kn_{i,j}^{-1})$ that summarized in Theorem \ref{tm:5}. Under Assumptions \ref{as:alphmix} and \ref{as:moment}, it follows from Davydov's inequality that $\max_{i,j\in[p]}|{\rm I}'_2(i,j)|+\max_{i,j\in[p]}|{\rm{I}}'_3(i,j)|\lesssim\exp(-C_*L_n^{-\varphi}K^{\varphi})$ for some universal constant $C_*>0$. If $\{\bU_{t_k}\}_{k=1}^n$ is an independent sequence or $L_n$-dependent sequence, we know $\varphi=\infty$ and then $\max_{i,j\in[p]}|{\rm I}'_2(i,j)|+\max_{i,j\in[p]}|{\rm{I}}'_3(i,j)|\lesssim\exp(-\infty)=0$ with $K>L_n$. If $\varphi<\infty$, with selecting $K=L_n(C_{**}\log n_*)^{1/\varphi}$ for some sufficiently large constant $C_{**}>0$, $\max_{i,j\in[p]}|{\rm I}'_2(i,j)|+\max_{i,j\in[p]}|{\rm{I}}'_3(i,j)|\lesssim n_*^{-C_*C_{**}}$ will be negligible in comparison to ${\rm II}(i,j)$.  

\begin{theorem}\label{tm:5}
Under Assumptions {\rm\ref{as:space}} and {\rm\ref{as:drifdiff}}, if $K=o(n_*)$,  it holds that
$$
\max_{i,j\in[p]}\bigg|\mathbb{E}\bigg\{{\rm II}(i,j)-\frac{2K+\Delta_K}{4n_{i,j}}\int_{t_{i,j,1}}^{t_{i,j,n_{i,j}}}\sigma_{i,s}\sigma_{j,s}\rho_{i,j,s}\,{\rm d}s\bigg\}\bigg|\lesssim (Kn_*^{-1})^{3/2}\,.$$
\end{theorem}

From \eqref{eq:highexp} and Theorem \ref{tm:5}, we have that
\begin{align}\label{eq:biasexp}
\mathbb{E}(\hat{\sigma}_{u,i,j})&=\sigma_{u,i,j}+\underbrace{\mathbb{E}\bigg(\frac{2K+\Delta_K}{4n_{i,j}}\int_{t_{i,j,1}}^{t_{i,j,n_{i,j}}}\sigma_{i,s}\sigma_{j,s}\rho_{i,j,s}\,{\rm d}s\bigg)}_{O(Kn_*^{-1})}+ O(K^{3/2}n_*^{-3/2})
\end{align}
provided that $K=L_n(C_{**}\log n_*)^{1/\varphi}$ for some sufficiently large constant $C_{**}>0$. Since $K=o(n_*)$, the second term on the right-hand side of \eqref{eq:biasexp} is asymptotically negligible if $\sigma_{u,i,j}$ is not vanishing; our Theorem \ref{tm:5} implies that it is the leading term in the bias.
Impact from the bias on $\hat{\sigma}_{u,i,j}$ could be empirically substantial, especially when $\sigma_{u,i,j}$ is relatively small.


For any $i,j\in[p]$, let
$
\psi_{i,j}=\int_{t_{i,j,1}}^{t_{i,j,n_{i,j}}}\sigma_{i,s}\sigma_{j,s}\rho_{i,j,s}\,{\rm d}s$.
As a remedy, we propose a bias-correction for $\hat{\sigma}_{u,i,j}$ as follows:
\begin{equation}\label{eq:biascorrect}
\hat{\sigma}_{u,i,j}^{{\rm bc}}=\hat{\sigma}_{u,i,j}-\frac{2K+\Delta_K}{4n_{i,j}}\hat{\psi}_{i,j}\,,
\end{equation}
where $\hat{\sigma}_{u,i,j}$ is given in \eqref{eq:estasyn}, and $\hat{\psi}_{i,j}$ is an estimate of $\psi_{i,j}$.
Since $\psi_{i,j}$ is an integrated covariance, it can be estimated by existing approaches, for example,  the polarization method
 \citep{Aitetal_2010_JASA}, the two time scales approach \citep{Zhang2011}, the pre-averaging method \citep{Jacod2009,Christensen2010}, and the quasi-maximum likelihood approach \citep{LiuTang_2014_JOE}. Section \ref{sec:bias-correction1} gives details for calculating $\hat{\psi}_{i,j}$ by the two time scales approach. Based on $\hat{\sigma}_{u,i,j}^{{\rm bc}}$ given in (\ref{eq:biascorrect}), we can obtain $\widehat{\bSigma}_u^{{\rm bc}}$ and $\widehat{\bSigma}_u^{{\rm bc, thre}}$, the bias-corrected version of $\widehat{\bSigma}_u$ defined as (\ref{eq:hatsigma}) and $\widehat{\bSigma}_{u}^{{\rm thre}}$ defined as (\ref{eq:est-thre}), respectively, by replacing $\hat{\sigma}_{u,i,j}$ by $\hat{\sigma}_{u,i,j}^{{\rm bc}}$. Theorem  \ref{corr1} indicates that $\widehat{\bSigma}_u^{{\rm bc}}$ and $\widehat{\bSigma}_u^{{\rm bc, thre}}$ share the same convergence rates of $\widehat{\bSigma}_u$ and $\widehat{\bSigma}_u^{{\rm thre}}$, respectively.

\begin{theorem}\label{corr1}
Assume $\max_{i,j\in[p]}|\hat{\psi}_{i,j}|=O_\p(\log p)$. The following two assertions are satisfied:

	{\rm (i)} Under the conditions of Theorem {\rm \ref{pn:1}},
$
 |\widehat{\bSigma}_u^{{\rm bc}}-\bSigma_u|_\infty=O_\p\{(Kn_*^{-1}\log p)^{1/2}\}$.

{\rm (ii)} Under the conditions of Theorem {\rm \ref{tm:3}}, $ \|\widehat{\bSigma}_u^{{\rm bc,thre}}-\bSigma_u\|_2=O_\p\{c_p(Kn_*^{-1}\log p)^{(1-q)/2}\}$ for any covariance matrix $\bSigma_{u}\in\mathcal{H}(q,c_p,M)$.
 \end{theorem}

\section{ Impact from jumps and the robustness of our methods} \label{s5j}

We now consider the setting with jumps in the underlying process $\bX_t$. Assume $\bX_t=(X_{1,t},\ldots,X_{p,t})^{\top}$ satisfies the following model:
\begin{equation}\label{eq:model2}
{\rm d} X_{i,t}=\mu_{i,t}\,{\rm d}t+\sigma_{i,t}\,{\rm d}B_{i,t}+J_{i,t}\,{\rm d}M_{i,t}~~~\textrm{and}~~~\mathbb{E}({\rm d}B_{i,t}\cdot{\rm d}B_{j,t})=\rho_{i,j,t}\,{\rm d}t\,,
\end{equation}
where $\mu_{i,t}$'s,  $\sigma_{i,t}$'s, $B_{i,t}$'s and $\rho_{i,j,t}$'s  are same as those in \eqref{eq:model}, ${J}_{i,t}$'s are the jump sizes, and $M_{i,t}$'s are   counting processes.
Our analysis reveals that our estimators proposed in Section \ref{sec:methfU} for $\bSigma_u$ are reasonably robust against jumps.  In our theoretical analysis, we impose Assumptions \ref{as:jump} and \ref{as:jumpsize} on the counting process $M_{i,t}$ and the jump size $J_{i,t}$, respectively.
\begin{as}\label{as:jump}
	Let $(\Delta M_{i,\cdot})_{i,j}=M_{i,t_{i,j,n_{i,j}}}-M_{i,t_{i,j,1}}$ for any $i,j\in[p]$.
	There exist $\lambda_1,\ldots,\lambda_p>0$ and some universal constants $C_{10}>0,C_{11}>0$ such that (i)
	$
	\mathbb{E}\{\exp(\theta[(\Delta M_{i,\cdot})_{i,j}-\mathbb{E}\{(\Delta M_{i,\cdot})_{i,j}\}])\}\leq \exp\{\lambda_i(t_{i,j,n_{i,j}}-t_{i,j,1})\theta^2\}
	$
	for any $|\theta|\leq C_{10}^{-1}$ and $i\in[p]$; (ii) $ \mathbb{E}\{(\Delta M_{i,\cdot})_{i,j}\}\leq C_{11}\lambda_i(t_{i,j,n_{i,j}}-t_{i,j,1})$ for any $i\in[p]$.
\end{as}

\begin{as}\label{as:jumpsize}
	There exist some universal   constants  $\iota>0,C_{12}>1$ and $C_{13}>0$ such that
	$
	\mathbb{P}(\sup_{0\leq t\leq  T}|J_{i,t}|>u)\leq C_{12}\exp(-C_{13}u^{\iota})
	$
	for any $i\in[p]$ and $u>0$.
\end{as}
If $M_{i,t}$ is a Poisson process with intensity $\lambda_i>0$, then $(\Delta M_{i,\cdot})_{i,j}$ follows the Poisson distribution with parameter $\lambda_i(t_{i,j,n_{i,j}}-t_{i,j,1})$ and Assumption \ref{as:jump} holds for $C_{10}=C_{11}=1$.
Assumption \ref{as:jumpsize} controls the tail behavior of the random jump size  $\sup_{0\leq t\leq  T}J_{i,t}$.
If the jump size  $J_{i,t}$ is  bounded from above uniformly for $i\in[p]$ and $t\in[0,T]$, we can select $\iota=\infty$ in Assumption \ref{as:jumpsize}.

Recall $S_{i,j,k}=\{t_{i, j, \ell}\in\mathcal{G}_i\cap\mathcal{G}_j:K\leq|\ell-k|\leq K+\Delta_K \}$ for some integers $K\geq1$ and $\Delta_K\geq0$, and $N_{i ,j , k}=|S_{i,j,k}|$.  For any $i,j\in[p]$, define
\begin{align}\label{eq:wij}
\varpi_{i,j}=\frac{1}{2n_{i,j}}\sum_{k=1}^{n_{i,j}}\frac{1}{N_{i,j,k}}\sum_{t_{i,j,\ell} \in S_{i,j,k}}\bigg(\int_{t_{i,j,k}}^{t_{i,j,\ell}}J_{i,s}\,{\rm d}M_{i,s}\bigg)\bigg(\int_{t_{i,j,k}}^{t_{i,j,\ell}}J_{j,s}\,{\rm d}M_{j,s}\bigg)\,.
\end{align}
Let  $\bOmega=(\varpi_{i,j})_{p\times p}$ and
\begin{align*}
\widehat{\bSigma}_{u}^{{\rm jump}}=(\hat{\sigma}_{u,i,j}^{{\rm jump}})_{p\times p}=\widehat{\bSigma}_{u}-\bOmega\,,
\end{align*}
where $\widehat{\bSigma}_{u}$ is defined as \eqref{eq:hatsigma}.
Analogous to \eqref{eq:est-thre}, we define the thresholding version of $\widehat{\bSigma}_{u}^{{\rm jump}}$ as
\begin{align}\label{eq:est-threJ}
\widehat{\bSigma}_{u}^{{\rm jump,  thre}}=\big[\hat{\sigma}_{u,i,j}^{{\rm jump}}{I}\big\{|\hat{\sigma}_{u,i,j}^{{\rm jump}}|\geq \beta(Kn_*^{-1 }\log p)^{1/2}\big\}\big]_{p\times p}\,,
\end{align}
where $\beta>0$ is a fixed constant for the thresholding level.

Our theory has two parts.
As the first part, parallel to Theorems \ref{pn:1} and \ref{tm:3}, we have the next theorem for the convergence rates of $\widehat{\bSigma}_{u}^{{\rm jump}}$ and $\widehat{\bSigma}_{u}^{{\rm jump, thre}}$.

\begin{theorem}\label{eq:the6}
Assume Assumptions {\rm \ref{as:alphmix}}--{\rm \ref{as:jumpsize}} hold. Let $\lambda_*=\max_{i\in[p]}\lambda_i$ and  $K> CL_n$ for some constant $C\geq 1$. If $\lambda_*^2(n_*^{-1}K)^{\iota/(2\iota+2)} (\log p)^{-1}=o(1)$, $\log p=o[\min\{(n_*L_n^{-2}K)^{\varphi/(3\varphi+2)},(n_*K^{-1})^{\tilde{\chi}}\}]$ and  $K^{-\varphi}L_n^{\varphi}\log\{n_*(K\log p)^{-1}\}=o(1)$, where $n_*$ is specified in {\rm(\ref{eq:nstar})} and  $\tilde{\chi}=\min\{\iota\gamma/(2\iota+2\gamma+\iota\gamma),\iota/(2\iota+2),\gamma/(\gamma+4),1/3\}$, then the following two assertions are satisfied:

{\rm(i)} Let $\mathcal{P}_3$ denote the collections of models for $\{\bY_{t_k}\}_{k=1}^n$ such that $\bY_{t_k}=\bX_{t_k}+\bU_{t_k}$, where the noises $\{\bU_{t_k}\}_{k=1}^n$ satisfy Assumption {\rm\ref{as:moment}}, $\bX_t=(X_{1,t},\ldots,X_{p,t})^{\top}$ follows 
model {\rm(\ref{eq:model2})} with each $\mu_{i,t}$, $\sigma_{i,t}$, $M_{i,t}$ and  $J_{i,t}$ satisfying Assumptions {\rm\ref{as:drifdiff}}--{\rm\ref{as:jumpsize}}, and the grids of time points $\{\mathcal{G}_i\}_{i=1}^p$ satisfy Assumption {\rm\ref{as:space}}. It holds that
$$
\sup_{\mathcal{P}_3}\mathbb{E}\big(|\widehat{\bSigma}_{u}^{{\rm jump}}-\bSigma_{u}|_\infty\big)\lesssim(Kn_*^{-1}\log p)^{1/2}\,.$$

{\rm(ii)} Let $\mathcal{P}_4$ denote the collections of models for $\{\bY_{t_k}\}_{k=1}^n$ such that $\bY_{t_k}=\bX_{t_k}+\bU_{t_k}$, where the noises $\{\bU_{t_k}\}_{k=1}^n$ satisfy Assumption {\rm\ref{as:moment}} with the covariance matrix $\bSigma_{u}\in\mathcal{H}(q,c_p,M)$, $\bX_t=(X_{1,t},\ldots,X_{p,t})^{\top}$ follows 
model {\rm(\ref{eq:model2})} with each $\mu_{i,t}$, $\sigma_{i,t}$, $M_{i,t}$ and  $J_{i,t}$ satisfying Assumptions {\rm\ref{as:drifdiff}}--{\rm\ref{as:jumpsize}}, and the grids of time points $\{\mathcal{G}_i\}_{i=1}^p$ satisfy Assumption {\rm\ref{as:space}}. With sufficiently large constant $\beta>0$ in \eqref{eq:est-threJ}, it holds that
$$
\sup_{\mathcal{P}_4}\mathbb{E}\big(\|\widehat{\bSigma}_{u}^{{\rm jump,  thre}}-\bSigma_{u}\|_2^2\big)\lesssim c_p^2(Kn_*^{-1}\log p)^{1-q}\,.$$
\end{theorem}

Theorems \ref{eq:the6}(i) and \ref{eq:the6}(ii) can be viewed as the generalization of Theorems \ref{pn:1} and \ref{tm:3}, respectively. If there are no jumps in $\bX_t$, we have  $J_{i,t}\equiv0$ for any $i\in[p]$ and $t\in[0,T]$, and then $\widehat{\bSigma}_{u}^{{\rm jump}}=\widehat{\bSigma}_{u}$ and $\widehat{\bSigma}_{u}^{{\rm jump,  thre}}=\widehat{\bSigma}_{u}^{{\rm thre}}$. In this scenario, we can set $\lambda_i=0$ for any $i\in[p]$ and $\iota=\infty$ in Assumptions \ref{as:jump} and \ref{as:jumpsize}, respectively, which implies  $\lambda_*^2(n_*^{-1}K)^{\iota/(2\iota+2)} (\log p)^{-1}=o(1)$ holds automatically and $\tilde{\chi}=\min\{\gamma/(\gamma+4),1/3\}$. Hence, the results of  Theorems \ref{eq:the6}(i) and \ref{eq:the6}(ii) in this scenario
are identical to Theorems \ref{pn:1} and \ref{tm:3}, respectively.

In the second part of our theory, we  establish the properties of $\widehat{\bSigma}_{u}$ itself.
Since $\widehat{\bSigma}_{u}^{{\rm jump}}=\widehat{\bSigma}_{u}-\bOmega$ where $\bOmega=(\varpi_{i,j})_{p\times p}$ with $\varpi_{i,j}$ defined as \eqref{eq:wij},   the gap between $\widehat{\bSigma}_{u}$ and $\widehat{\bSigma}_{u}^{{\rm jump}}$ is seen determined by  an extra `bias'.
We assume the following condition for controlling
 the tail probability of $\varpi_{i,j}$.   
\begin{as}\label{as:jump1}
Each $M_{i,t}$ has  independent increments.	Let $(\Delta M_{i,\cdot})_{i,j}^{(k)}=M_{i,t_{i,j,k+1}}-M_{i,t_{i,j,k}}$ for any $i,j\in[p]$ and $k\in[n_{i,j}-1]$.
	There exist $\lambda_1,\ldots,\lambda_p>0$ and some universal constants $C_{10}>0,C_{11}>0$ such that (i)
	$
	\mathbb{E}\{\exp(\theta[(\Delta M_{i,\cdot})_{i,j}^{(k)}-\mathbb{E}\{(\Delta M_{i,\cdot})_{i,j}^{(k)}\}])\}\leq \exp\{\lambda_i(t_{i,j,k+1}-t_{i,j,k})\theta^2\}
	$
	for any $|\theta|\leq C_{10}^{-1}$ and  $i\in[p]$; (ii) $ \mathbb{E}\{(\Delta M_{i,\cdot})_{i,j}^{(k)}\}\leq C_{11}\lambda_i(t_{i,j,k+1}-t_{i,j,k})$ for any $i\in[p]$.
\end{as}

 Assumption \ref{as:jump} holds automatically under Assumption \ref{as:jump1}. If $M_{i,t}$ is a Poisson process with intensity $\lambda_i>0$, then $(\Delta M_{i,\cdot})_{i,j}^{(k)}$ follows the Poisson distribution with parameter $\lambda_i(t_{i,j,k+1}-t_{i,j,k})$ and Assumption \ref{as:jump1} holds for $C_{10}=C_{11}=1$. Given Assumptions \ref{as:jumpsize} and \ref{as:jump1}, we have the following theorem for the tail probability of $\varpi_{i,j}$.
\begin{theorem}\label{thwij} Let $\lambda_*=\max_{i\in[p]}\lambda_i$. Under Assumptions {\rm \ref{as:jumpsize}} and {\rm \ref{as:jump1}}, if $K(\log n_*)^{(3\iota+4)/\iota}=o(n_*)$ and $n_*^{-1}K\lambda_*=o(1)$, it holds that
$$
	\max_{ i,j\in[p]}\mathbb{P}(|\varpi_{i,j}|\geq v)\lesssim \exp\{-C(n_*K^{-1})^{(2\iota+2)/(3\iota+4)}v\}+ \exp\{-C(n_*K^{-1})^{\iota/(3\iota+4)}\}$$
for any $v\gg (n_*^{-1}K)^{(2\iota+2)/(3\iota+4)}\lambda_*$, where  $\iota$ is specified in Assumption  $\ref{as:jumpsize}$. Furthermore, it holds that $\max_{ i,j\in[p]}\mathbb{E}(|\varpi_{i,j}|^m)\lesssim 1$ for any fixed positive integer $m$.
\end{theorem}

Theorem \ref{thwij} implies $|\bOmega|_\infty=O_{\p}\{(Kn_*^{-1})^{(2\iota+2)/(3\iota+4)}\log p\}=o_{\p}\{(Kn_*^{-1}\log p)^{1/2}\}$ if $\log p=o\{(n_*K^{-1})^{\iota/(3\iota+4)}\}$, which leads to the robustness of our proposed estimators 
against possible jumps in the underlying process $\bX_t$, as established in the following theorem.

\begin{theorem}\label{eq:jjj}
Let $\lambda_*=\max_{i\in[p]}\lambda_i$ and  $K> CL_n$ for some constant $C\geq 1$. Under Assumptions {\rm \ref{as:alphmix}}--{\rm \ref{as:diff}} and {\rm \ref{as:jumpsize}}--{\rm \ref{as:jump1}}, if  $\lambda_*^2(n_*^{-1}K)^{\iota/(3\iota+4)}(\log p)^{-1}=o(1)$, $\log p=o[\min\{(n_*L_n^{-2}K)^{\varphi/(3\varphi+2)},(n_*K^{-1})^{\chi^*}\}]$  and  $K^{-\varphi}L_n^{\varphi}\log\{n_*(K\log p)^{-1}\}=o(1)$, where $n_*$ is specified in {\rm(\ref{eq:nstar})} and $\chi^*=\min\{\iota\gamma/(2\iota+2\gamma+\iota\gamma),\iota/(3\iota+4),\gamma/(\gamma+4)\}$, then
 the following two assertions are satisfied:
	
	{\rm (i)}
$
|\widehat{\bSigma}_u-\bSigma_u|_\infty=O_\p\{(Kn_*^{-1}\log p)^{1/2}\}$.

{\rm (ii)} With sufficiently large constant $\beta>0$ in \eqref{eq:est-thre}, $ \|\widehat{\bSigma}_u^{{\rm thre}}-\bSigma_u\|_2=O_\p\{c_p(Kn_*^{-1}\log p)^{(1-q)/2}\}$ for any covariance matrix $\bSigma_{u}\in\mathcal{H}(q,c_p,M)$.
\end{theorem}

Specifically,  Theorem \ref{eq:jjj} implies that if the intensity parameter $\lambda_*$ of the counting processes diverges no faster  than $(n_*K^{-1})^{\iota/(6\iota+8)}(\log p)^{1/2}$, the convergence rates of $\widehat{\bSigma}_u$ and $\widehat{\bSigma}_u^{\rm thre}$ are the same as those in Theorems \ref{pn:1} and \ref{tm:3}, even there are jumps in the processes.  

\section{Numerical studies}\label{s3}
\subsection{Simulations}
\subsubsection{Data generating procedure}
We set $(\rho_{i,j})_{p\times p}=\{\diag(\bA^{\otimes2})\}^{-1/2}\bA^{\otimes2}\{\diag(\bA^{\otimes2})\}^{-1/2}$, where $\bA=(a_{i,j})_{p\times p}$ with $a_{i,j}=(-0.8)^{|i-j|}I(i\geq j)$. For $\bX_t=(X_{1,t},\ldots,X_{p,t})^{\top}$, we generated each $X_{i,t}$ from the following stochastic volatility model:
\begin{align}\label{eq:sigmai}
{\rm d}X_{i, t} &= \sigma_{i, t}\,{\rm d}B_{i, t}+J_{i,t}\,{\rm d}M_{i,t}\,,~~~~~
{\rm d}\sigma_{i, t}^2 = \kappa(\bar{\sigma}^2-\sigma_{i, t}^2)\,{\rm d}t+s\sigma_{i,t}\,{\rm d}W_{i,t}\,, ~~~~
\end{align}
where $M_{1,t},\ldots,M_{p,t}$ are independent Poisson processes with intensity $\lambda_{\rm J}$, and $B_{1,t},\ldots,B_{p,t},W_{1,t},\ldots,W_{p,t}$ are univariate standard Brownian motions such that (i) $W_{1,t},\ldots, W_{p,t}$ are $p$ independent Brownian motions, and (ii) $\mathbb{E}({\rm d}B_{i, t}\cdot {\rm d}W_{j, t})=-0.3I(i=j)\,{\rm d}t$ and $\mathbb{E}({\rm d}B_{i, t}\cdot {\rm d}B_{j, t})=\rho_{i,j}\,{\rm d}t$. We considered two settings -- with or without jumps: (i) $J_{1,t},\ldots,J_{p,t}\overset{{\rm i.i.d.}}{\sim}N(0,\sigma_{{\rm J}}^2)$ which are independent of  $M_{1, t},\ldots,M_{p,t}$, and (ii) $J_{i,t}\equiv0$ for  $i\in[p]$ and $t$. 
We set $(\kappa, s, \bar{\sigma}^2, \sigma_{{\rm J}}^2, \lambda_{{\rm J}})=(5, 0.5, 0.1, 0.015^2, 5)$, the same as that in the numerical studies of \cite{Aitetal_2013_JFE} that mimics the empirical features of  financial data \citep{Ait-Sahalia2007}.  This setting is reasonable; comparable settings are  found in existing studies  \citep{AY2009, Aitetal_2010_JASA,LiuTang_2014_JOE,Ait_Xiu2017_JASA}. The initial observations of $\sigma_{i, t}^2$ $(i\in[p])$ were generated from a Gamma distribution $\Gamma(2\kappa\bar{\sigma}^2/s^2, s^2/(2\kappa))$. In our simulation, we set $p\in\{50,100,200\}$. 

We
took $t\in[0, T]$ with $T = 1/252$;   here 1 unit of $t$ means one year, so $T=1/252$ is corresponding to a trading day.  We first  generated high-frequency data available at each second in a $6.5$-hour  period; this setting results in $60\times60\times6.5=23400$ observations. By letting $\tilde{t}_k=k/(252\times 23400)$, 
we generated  $\bY_{\tilde{t}_k}=\bX_{\tilde{t}_k}+\bU_{\tilde{t}_k}$ 
with $(\bX_t)_{t\in[0,T]}$ from \eqref{eq:sigmai},  and each element of $\bU_{\tilde{t}_k}=(U_{1,\tilde{t}
_k},\ldots,U_{p,\tilde{t}_k})^{\top}$ from a stationary GARCH(1,1) model:
\begin{align*}
U_{i,\tilde{t}_k}=\tilde{\sigma}_{i,\tilde{t}_{k}}\eta_{i,\tilde{t}_{k}}\,,~~~\tilde{\sigma}_{i ,\tilde{t}_k}^2=0.1\sigma_e^2+0.1\tilde{\sigma}_{i,\tilde{t}_{k-1}}^2\eta_{i,\tilde{t}_{k-1}}^2+0.8\tilde{\sigma}_{i,\tilde{t}_{k-1}}^2\,,
\end{align*}
where $\bfeta_{\tilde{t}_{k}}=(\eta_{1,\tilde{t}_{k}}, \ldots, \eta_{p,\tilde{t}_{k}})^{\top}$ is independently generated from $N(\bzero,\bR)$;  the  settings  of $\bR$ will be described later. 
In this model,  upon observing $\mathbb{E}(\tilde{\sigma}_{i,\tilde{t}_{k}}^2)=\sigma_e^2$ for each $i\in[p]$,  
we considered different settings for the signal-to-noise ratios $\sigma_e^2/(T^{-1}\int_0^{T}\sigma_{i,t}^2\,{\rm d}t)$ by varying  $\sigma_e^2$. Since the signal-to-noise ratio can be approximated by $252\sigma_e^2/\bar{\sigma}^2$ with $\bar{\sigma}^2$ in \eqref{eq:sigmai},  
we specified two selections of $\sigma_e^2$: (i) $\sigma_e^2=0.005^2$ such that $252\sigma_e^2/\bar{\sigma}^2=0.063$, and (ii) $\sigma_e^2=0.001^2$ such that $252\sigma_e^2/\bar{\sigma}^2=0.00252$.
In Part S1 of the supplementary material, we have also investigated the finite-sample performance of the proposed estimators when $\bU_{\tilde{t}_1},\ldots,\bU_{\tilde{t}_{23400}}\overset{{\rm i.i.d.}}{\sim}N(\bzero,\sigma_e^2\bR)$ with $\sigma_e^2=0.005^2$ and $0.001^2$. 

We studied the following three models for $\bR=(r_{i,j})_{p\times p}$ that controls the correlations:
\vspace{-11pt}
\begin{description}
\item Model 1:  $\bR$ is a banded matrix, where $r_{i,i}=1$, $r_{i+1,i}=r_{i, i+1} =0.6$, $r_{i+2,i}=r_{i,i+2}=0.3$, and $r_{i,j} =0$ for $|i-j|\geq 3$.

\vspace{-6pt}

\item Model 2: $\bR=\{\diag(\bR^*)\}^{-1/2}{\bR^*}\{\diag(\bR^*)\}^{-1/2}$, where $\bR^*=\widetilde{\bR}+\{|\lambda_{\min}(\widetilde{\bR})|+0.05\}\bI_p$, $\bI_p$ is the identity matrix of order $p$, $\lambda_{\min}(\widetilde{\bR})$ is the smallest eigenvalue of $\widetilde{\bR}$, and $\widetilde{\bR}=(\tilde{r}_{i, j})_{p\times p}$ satisfies that $\tilde{r}_{i,j}=w_{i,j}b_{i,j}$, $w_{i, j}$'s are independently generated from the uniform  distribution $U(0.4, 0.8)$, $b_{i, j}$'s are independently generated from the Bernoulli distribution with successful probability 0.04.

\vspace{-6pt}

\item Model 3:  $\bR$ is a bandable matrix with $r_{i,j}=0.6^{|i-j|}$.
\end{description}

\vspace{-11pt}

We considered both synchronous and asynchronous high-frequency data in our simulation. To model the synchronous data setting,
we took $\{\bY_{\tilde{t}_{k\Delta}}\}_{k=1}^{\lfloor23400/\Delta\rfloor}$ as the observed data where $\lfloor x \rfloor$ is the floor function;  by varying $\Delta$, we simulated data sets of different sizes:  
larger $\Delta$ means fewer  observations.  
Then the time points 
where we observed the noisy data are  $t_k=\tilde{t}_{k\Delta}$ with $k=1,\ldots,\lfloor23400/\Delta\rfloor$. In our numerical studies, we set $\Delta\in\{1,2,3\}$. To model the asynchronous data setting, for each $i\in[p]$, we applied the Poisson process sampling scheme with intensity $\lambda$ to $\{\tilde{t}_k\}_{k=1}^{23400}$ for generating $\mathcal {G}_i=\{t_{i,1},\ldots,t_{i,n_i}\}$, the grid of time points at which we actually observed $Y_{i,t}$. The Poisson process sampling schemes for different $i$'s are independent.
Based on this setting for asynchronous data, on average there were $\lfloor23400/\lambda\rfloor$ observations for each $Y_{i,t}$. We selected $\lambda\in\{1, 2,3\}$ in our simulation.    

\subsubsection{Implementation of bias-correction}\label{sec:bias-correction1}


To obtain the bias-corrected estimator $\hat{\sigma}_{u,i,j}^{{\rm bc}}$ in \eqref{eq:biascorrect}, we need to calculate $\hat{\psi}_{i,j}$, the estimate of the integrated covariance $\psi_{i,j}$. In the simulation, we applied the two time scales approach \citep{Zhang2011}  to estimate $\psi_{i,j}$. Recall $\mathcal {G}_i=\{t_{i,1},\ldots,t_{i,n_i}\}$ is the grid of time points we observed $Y_{i,t}$.  For given $i,j\in[p]$, we  first used the  refresh time procedure  \citep{Barndorff_2011_JOE}  to synchronize the data if $\mathcal{G}_i\neq\mathcal{G}_j$. More specifically,  let the first refresh time point be
$v_{1}=\max(t_{i,1}, t_{j,1})$, and then define the other refresh time points $v_{l+1}$ with $l\geq1$ as
$
v_{l+1}=\max[ \min\{ t\in \mathcal {G}_{i}:  t>v_{l}\},  \min\{ t\in \mathcal {G}_j: t>v_{l}\}]$ iteratively.
Denoted by $n_{i,j}^*$ the resulting refresh time points for $\mathcal {G}_i$  and $\mathcal {G}_j$, and write
$
{t}_{i,l}^*=\max\{ t\in \mathcal {G}_i: t\leq  v_{l}\}$ and ${t}_{j,l}^*=\max\{ t\in \mathcal {G}_j: t\leq  v_{l}\}$
for each $l\in[n^*_{i, j}]$. If $\mathcal{G}_i=\mathcal{G}_j$, the refresh time points based on above procedure are identical to $t_{i,1},\ldots,t_{i,n_i}$, and thus the associated $t_{i,l}^*=t_{j,l}^*=t_{i,l}$. For given positive integers $\delta_1$ and $\delta_2$, the two time scales estimator for $\psi_{i,j}$ is given by
\begin{equation}\label{eq:msrv}
	\hat{\psi}_{i,j}=[Y_{i}, Y_{j}]^{(\delta_1)}-\frac{\delta_2({n}_{i,j}^*-\delta_1+1)}{\delta_1({n}_{i,j}^*-\delta_2+1)}[Y_{i}\,, Y_{j}]^{(\delta_2)}\,,
\end{equation}
where
$
	[Y_{i}, Y_{j}]^{(\delta)}=\delta^{-1}\sum_{l=\delta+1}^{{n}_{i,j}^*}(Y_{i,{t}_{i,l}^*}-Y_{i,{t}^*_{i,l-\delta}})(Y_{j,{t}^*_{j,l}}-Y_{j,{t}^*_{j,l-\delta}})$ for any positive integer $\delta$. Following  \cite{AY2009}, we set $(\delta_1,\delta_2)=(25,1)$ in our simulation.

\subsubsection{Selections of $(K,\Delta_K)$ and the thresholding level}\label{sec:tuning}

To obtain $\hat{\sigma}_{u,i,j}$ defined as \eqref{eq:estasyn} in practice, we need to select the tuning parameters $K$ and $\Delta_K$. If each univariate sequence $\{U_{i,t_k}\}_{k=1}^n$ is $\alpha$-mixing with exponentially decaying $\alpha$-mixing coefficients\footnote{Such requirement can be easily satisfied in most commonly used univariate time series models. See our discussion below \eqref{eq:alphacoeff}.}, with the independent assumption imposed on the $p$ sequences $\{U_{1,t_k}\}_{k=1}^n,\ldots,\{U_{p,t_k}\}_{k=1}^n$,  Theorem 5.1 of \cite{Bradley2005} indicates that $\alpha_n(m)$ defined in \eqref{eq:alphacoeff} satisfies $\alpha_n(m)\leq p\exp(-cm)$ for some universal constant $c>0$, which provides a rough upper bound for $\alpha_n(m)$.  Hence, Assumption \ref{as:alphmix} holds for $\varphi=1$ and $L_n=\tilde{c}\log p$ for some sufficiently small constant $\tilde{c}>0$. Our theoretical results require $K^{-\varphi}L_n^{\varphi}\log\{n_*(K\log p)^{-1}\}=o(1)$ and $K>CL_n$ for some constant $C\geq 1$. To match these requirements, when we estimate ${\sigma}_{u,i,j}$,
we can select $K=\bar{c} (\log p)(\log n_{i,j})(\log\log n_{i,j})$ for some small constant $\bar{c}>0$. In our simulation, we have tried $\bar{c}\in[0.03,0.07]$ and the associated results are similar. We suggest to select $\bar{c}=0.05$ in practice. Since we use $(Y_{i,t_{i,j,\ell}}-Y_{i,t_{i,j,k}},Y_{j,t_{i,j,\ell}}-Y_{j,t_{i,j,k}})$ to approximate $(U_{i,t_{i,j,\ell}}-U_{i,t_{i,j,k}},U_{j,t_{i,j,\ell}}-U_{j,t_{i,j,k}})$ for $K\leq |\ell-k|\leq K+\Delta_K$, the bias issue caused by $X_{i,t_{i,j,\ell}}-X_{i,t_{i,j,k}}$ and
$X_{j,t_{i,j,\ell}}-X_{j,t_{i,j,k}}$ will impact the performance of our estimators. Notice that a smaller $\Delta_K$ results in a smaller bias. We need to select $\Delta_K$ as some small positive integers.  Table \ref{Error} shows that (i) the estimators with $\Delta_K\in\{1,2,3\}$ perform quite well, and (ii) the estimators with $\Delta_K=1$ work best in most cases and perform quite close to the best ones in other cases. This verifies our claim that $\Delta_K$ should be selected as some small integers. We suggest to select $\Delta_K\in\{1,2,3\}$ in practice.

 Based on $\widehat{\bSigma}_u=(\hat{\sigma}_{u,i,j})_{p\times p}$ and $\widehat{\bSigma}_u^{\rm bc}=(\hat{\sigma}_{u,i,j}^{\rm bc})_{p\times p}$,  to derive their thresholding version $\widehat{\bSigma}_u^{\rm thre}$ and $\widehat{\bSigma}_u^{{\rm bc, thre}}$, we need to determine the thresholding level. Our theoretical analysis shows that the thresholding level should have the order $(Kn_*^{-1}\log p)^{1/2}$. Notice that $\hat{\sigma}_{u,i,j}=n_{i,j}^{-1}\sum_{k=1}^{n_{i,j}}\zeta_{i,j,k}$ with $
\zeta_{i,j,k}=(2N_{i ,j , k})^{-1}\sum_{t_{i,j,\ell}\in S_{i,j,k}}(Y_{i,t_{i, j, \ell}}-Y_{i, t_{i, j, k}})(Y_{j,t_{i, j, \ell}}-Y_{j, t_{i, j, k}})$. Since ${\rm Var}(\hat{\sigma}_{u,i,j})$ has the order $Kn_*^{-1}$, the long-run variance of the sequence $\{\zeta_{i,j,k}\}_{k=1}^{n_{i,j}}$ has the order $K$. To incorporate the heterogeneity of the estimators $\hat{\sigma}_{u,i,j}$, we implemented the thresholding estimators in practice as
\begin{equation}\label{eq:praest}
\begin{split}
\widehat{\bSigma}_{u}^{\rm thre}&=(\hat{\sigma}^{\rm thre}_{u,i, j})_{p\times p}=\big[\hat{\sigma}_{u, i, j}{I}\big\{|\hat{\sigma}_{u, i, j}|\geq \beta_*(\hat{\theta}_{i,j}{n}_{i,j}^{-1}\log p)^{1/2}\big\}\big]_{p\times p}\,,\\
\widehat{\bSigma}_{u}^{\rm bc, thre}&=(\hat{\sigma}^{\rm bc, thre}_{u,i, j})_{p\times p}=\big[\hat{\sigma}_{u, i, j}^{{\rm bc}} {I}\big\{|\hat{\sigma}_{u, i, j}^{{\rm bc}}|\geq \beta_*(\hat{\theta}_{i,j}{n}_{i,j}^{-1}\log p)^{1/2}\big\}\big]_{p\times p}\,,
\end{split}
\end{equation}
where $\beta_*>0$ is a constant, and
 $\hat{\theta}_{i,j}$ is an estimate for the long-run variance of the sequence $\{\zeta_{i,j,k}\}_{k=1}^{n_{i,j}}$. Write $\bar{\zeta}_{i, j}=n_{i,j}^{-1}\sum_{k=1}^{n_{i, j}}\zeta_{i,j,k}$. We chose $\hat{\theta}_{i,j}$   in  \eqref{eq:praest} as 
\[
\hat{\theta}_{i,j}=\sum_{\ell=-n_{i,j}+1}^{n_{i,j}-1}\mathcal{K}\bigg(\frac{\ell}{h}\bigg)\widehat{H}_{i,j}(\ell)\,,
\]
where $\mathcal{K}(\cdot)$ is a symmetric kernel function, $h$ is the bandwidth, $\widehat{H}_{i,j}(\ell)=n_{i,j}^{-1}\sum_{k=\ell+1}^{n_{i,j}}(\zeta_{i,j,k}-\bar{\zeta}_{i, j})(\zeta_{i,j,k-\ell}-\bar{\zeta}_{i, j})$ for $\ell\geq0$ and $\widehat{H}_{i,j}(\ell)=n_{i,j}^{-1}\sum_{k=-\ell+1}^{n_{i,j}}(\zeta_{i,j,k+\ell}-\bar{\zeta}_{i, j})(\zeta_{i,j,k}-\bar{\zeta}_{i, j})$ otherwise. \cite{Andrews1991} suggested the quadratic spectral kernel
\[
\mathcal{K}(x)=\frac{25}{12\pi^2x^2}\bigg\{\frac{\sin(6\pi x/5)}{6\pi x/5}-\cos(6\pi x/5)\bigg\}
\]
with optimal bandwidth
$
h=1.3221\{4 n_{i,j}\hat{\vartheta}_{i,j}^2(1-\hat{\vartheta}_{i,j})^{-4}\}^{1/5}
$,
where $\hat{\vartheta}_{i,j}$ is the estimated autoregressive coefficient in the fitted AR(1) model for the sequence $\{\zeta_{i,j,k}\}_{k=1}^{n_{i,j}}$.  In our simulation, we have tried $\beta_*\in[1.75,2.25]$ and the associated results are similar. We suggest to select $\beta_*=2$ in practice.

\subsubsection{Simulation results}

For given estimator $\widetilde{\bSigma}$, we evaluated its relative estimation error
$\|\widetilde{\bSigma}-\bSigma_{u}\|_{2}/\|\bSigma_{u}\|_{2}$ in different settings. Table \ref{Error} summarizes the averages of the relative estimation errors based on 1000 repetitions.  We have several observations.
First, we find that in general, $\widehat{\bSigma}_{u}^{\rm bc, thre}$ performs quite well for all cases with satisfactorily small relative estimation errors compared with $\widehat{\bSigma}_{u}^{\rm thre}$. Further, $\widehat{\bSigma}_{u}^{\rm thre}$ performs quite well when $252\sigma_e^2/\bar{\sigma}^2=0.063$ but poorly when $252\sigma_e^2/\bar{\sigma}^2=0.00252$. This suggests that when the noise is quite small, the bias-correction is necessary. Second, as the dimension $p$ increases,  the relative estimation errors worsen a bit,  but at a very slow pace growing with $p$. This demonstrates the promising performance of the thresholding method for handling high-dimensional covariance estimations.
 Third,  as the sampling frequency becomes higher (smaller $\Delta$ or $\lambda$), the performance is improved by observing smaller relative estimation errors, reflecting the blessing to the covariance estimations with more high-frequency data. This is actually the reason why the performance of the estimator with synchronous data is better than that with asynchronous data when $\Delta$ and $\lambda$ are the same.
Fourth, we find that the differences are small among the performances with different $\Delta_K$, especially when the data are synchronous. 
Fifth, we find that the empirical performance of the proposed estimators is robust to jumps, confirming our finding in Theorem \ref{eq:jjj}. 

 In addition, for given estimator $\widetilde{\bSigma}=(\tilde{\sigma}_{i,j})_{p\times p}$, we also calculated in Tables \ref{TPRFPR_Syn} and \ref{TPRFPR_Asyn} the true positive rate (TPR) and the false positive rate (FPR) defined as
 \begin{align*}
 {\rm TPR}&=\frac{|\{ (i, j): \tilde{\sigma}_{i, j}\neq 0 ~{\rm and} ~{\sigma}_{u, i, j}\neq 0\}|}{|\{ (i, j): {\sigma}_{u, i, j}\neq 0\}|}\,,\\
 {\rm FPR}&=\frac{|\{ (i, j): \tilde{\sigma}_{i, j}\neq 0 ~{\rm and} ~{\sigma}_{u, i, j}= 0\}|}{|\{ (i, j): {\sigma}_{u, i, j}=0\}|}\,.
 \end{align*}
Since the covariance matrix considered in Model 3 has no exact zero element, we omit reporting the TPR and FPR in this case. Results in Tables \ref{TPRFPR_Syn} and \ref{TPRFPR_Asyn} show that the TPRs of our proposed estimators for all cases are equal to 1 or quite close to 1, and the FPRs for all cases are almost 0. This  indicates that our proposed thresholding method can recover the non-zero elements of the covariance matrix very accurately. From the results in Table \ref{TPRFPR_Asyn} when the data are asynchronous with $\lambda=2$ and 3,   we find that $\widehat{\bSigma}_u^{{\rm bc, thre}}$ performs a bit better than $\widehat{\bSigma}_u^{{\rm thre}}$, and both $\widehat{\bSigma}_u^{{\rm bc, thre}}$ and $\widehat{\bSigma}_u^{{\rm thre}}$ have lower TPRs when $252\sigma_e^2/\bar{\sigma}^2=0.00252$, which is reasonable as the signal-to-noise ratio in term of estimating the covariance matrix of noises is lower in this case.    For the FPRs, we find that there is no big difference between the FPRs of  $\widehat{\bSigma}_u^{{\rm bc, thre}}$ with different values for $252\sigma_e^2/\bar{\sigma}^2$. However, the FPRs of  $\widehat{\bSigma}_u^{{\rm thre}}$ when $252\sigma_e^2/\bar{\sigma}^2=0.00252$ are much higher than those    when $252\sigma_e^2/\bar{\sigma}^2=0.063$, showing the impact from weaker signal.   This again suggests that the bias-correction is very helpful, especially for handling relatively weaker signals.

  \begin{table}[htbp]
  	\tiny
  	\centering
  	\caption{Averages of the relative estimation errors ($\times100$) for the proposed estimators when jumps exist and do not exist (in parentheses) based on 1000 repetitions.}
	\medskip
  	\begin{tabular}{ccccccccccccc}
  		\hline
		\hline
  	\multicolumn{3}{c}{Synchronous Data} &  &  \multicolumn{3}{c}{Model 1}  &  \multicolumn{3}{c}{Model 2}  &  \multicolumn{3}{c}{Model 3}\\
  		${252\sigma_e^2}/{\bar{\sigma}^2}$  &   $p $    & Estimators  & $\Delta_K$ & $\Delta=3$     & $\Delta=2$     & $\Delta=1$     & $\Delta=3$     & $\Delta=2$     & $\Delta=1$     & $\Delta=3$     & $\Delta=2$     & $\Delta=1$ \\
  		\hline
  		
  		        0.063 & 50    & $\widehat{\bSigma}_{u}^{\rm  thre} $      & 1     & 4.5(4.5) & 4.0(4.0)  & 3.2(3.2) & 5.2(5.5) & 4.6(4.9) & 3.3(3.4) & 6.3(6.3) & 5.3(5.3) & 4.3(4.3) \\
          &       &       & 2     & 4.5(4.5) & 4.0(4.0)  & 3.2(3.2) & 5.3(5.6) & 4.7(4.9) & 3.3(3.5) & 6.1(6.0) & 5.2(5.1) & 4.2(4.2) \\
          &       &       & 3     & 4.6(4.6) & 4.1(4.1) & 3.2(3.3) & 5.6(5.9) & 4.8(5.1) & 3.4(3.6) & 5.9(5.9) & 5.0(5.0)  & 4.1(4.1) \\
          &       & $\widehat{\bSigma}_{u}^{\rm bc, thre} $      & 1     & 4.6(4.6) & 4.1(4.2) & 3.3(3.3) & 5.3(5.5) & 4.7(4.9) & 3.4(3.5) & 6.4(6.4) & 5.4(5.4) & 4.4(4.3) \\
          &       &       & 2     & 4.6(4.6) & 4.2(4.2) & 3.3(3.3) & 5.3(5.6) & 4.7(5.0) & 3.4(3.6) & 6.1(6.1) & 5.2(5.2) & 4.2(4.2) \\
          &       &       & 3     & 4.7(4.7) & 4.2(4.3) & 3.3(3.4) & 5.6(5.9) & 4.9(5.2) & 3.5(3.7) & 6.0(6.0)  & 5.1(5.1) & 4.1(4.1) \\
\cline{2-13}
          & 100   &  $\widehat{\bSigma}_{u}^{\rm thre} $     & 1     & 5.2(5.2) & 4.4(4.4) & 3.6(3.6) & 5.7(5.9) & 4.6(4.7) & 3.5(3.5) & 7.0(6.9) & 6.0(6.0)  & 4.8(4.7) \\
          &       &       & 2     & 5.2(5.2) & 4.4(4.4) & 3.6(3.6) & 5.8(5.9) & 4.6(4.7) & 3.5(3.6) & 6.7(6.7) & 5.8(5.8) & 4.6(4.6) \\
          &       &       & 3     & 5.3(5.3) & 4.5(4.5) & 3.6(3.7) & 6.0(6.1) & 4.8(4.9) & 3.6(3.7) & 6.6(6.5) & 5.7(5.7) & 4.6(4.5) \\
          &       & $\widehat{\bSigma}_{u}^{\rm bc, thre} $      & 1     & 5.3(5.4) & 4.5(4.5) & 3.7(3.7) & 5.8(6.0) & 4.7(4.8) & 3.6(3.7) & 7.1(7.0) & 6.1(6.1) & 4.8(4.8) \\
          &       &       & 2     & 5.4(5.4) & 4.5(4.6) & 3.7(3.7) & 5.9(6.1) & 4.8(4.9) & 3.7(3.7) & 6.9(6.8) & 5.9(5.9) & 4.7(4.6) \\
          &       &       & 3     & 5.5(5.5) & 4.6(4.7) & 3.8(3.8) & 6.1(6.3) & 5.0(5.1) & 3.8(3.9) & 6.8(6.7) & 5.8(5.8) & 4.6(4.6) \\
\cline{2-13}
          & 200   & $\widehat{\bSigma}_{u}^{\rm  thre} $      & 1     & 6.0(6.0)  & 5.0(5.0)  & 4.0(4.0)  & 5.5(5.5) & 4.4(4.5) & 3.3(3.4) & 7.5(7.6) & 6.5(6.5) & 5.2(5.1) \\
          &       &       & 2     & 5.9(5.9) & 5.0(5.0)  & 4.0(4.0)  & 5.4(5.5) & 4.4(4.5) & 3.3(3.3) & 7.3(7.3) & 6.3(6.3) & 5.0(5.0) \\
          &       &       & 3     & 5.9(5.9) & 5.0(5.0)  & 4.0(4.0)  & 5.4(5.5) & 4.5(4.5) & 3.3(3.3) & 7.2(7.2) & 6.3(6.2) & 5.0(4.9) \\
          &       & $\widehat{\bSigma}_{u}^{\rm bc, thre} $      & 1     & 6.2(6.2) & 5.2(5.2) & 4.2(4.2) & 5.7(5.7) & 4.6(4.7) & 3.6(3.6) & 7.7(7.8) & 6.6(6.6) & 5.3(5.2) \\
          &       &       & 2     & 6.3(6.3) & 5.3(5.3) & 4.2(4.2) & 5.7(5.7) & 4.7(4.7) & 3.6(3.6) & 7.5(7.6) & 6.5(6.5) & 5.1(5.1) \\
          &       &       & 3     & 6.3(6.3) & 5.4(5.4) & 4.3(4.3) & 5.8(5.8) & 4.8(4.8) & 3.6(3.7) & 7.5(7.5) & 6.4(6.4) & 5.1(5.1) \\


        \hline

    0.00252 & 50    & $\widehat{\bSigma}_{u}^{\rm  thre} $      & 1     & 23.5(18.6) & 21.0(16.6) & 9.9(8.0) & 40.6(34.8) & 36.4(30.9) & 16.5(14.2) & 19.9(16.3) & 17.1(14.1) & 9.0(7.5) \\
          &       &       & 2     & 28.2(22.1) & 24.3(19.0) & 11.5(9.1) & 49.0(41.6) & 42.1(35.7) & 19.4(16.6) & 22.7(18.6) & 19.0(15.7) & 10.0(8.3) \\
          &       &       & 3     & 33.0(25.9) & 27.5(21.6) & 13.1(10.4) & 57.4(49.0) & 47.6(40.7) & 22.4(19.0) & 25.5(20.9) & 20.9(17.3) & 10.9(9.1) \\
          &       &  $\widehat{\bSigma}_{u}^{\rm bc, thre} $     & 1     & 4.6(4.6) & 4.1(4.1) & 3.3(3.3) & 5.7(5.8) & 4.9(5.1) & 3.4(3.5) & 6.7(6.7) & 5.6(5.5) & 4.4(4.4) \\
          &       &       & 2     & 4.7(4.6) & 4.1(4.1) & 3.3(3.3) & 5.9(6.1) & 5.0(5.2) & 3.4(3.6) & 6.5(6.4) & 5.4(5.4) & 4.3(4.3) \\
          &       &       & 3     & 4.8(4.7) & 4.2(4.2) & 3.3(3.3) & 6.3(6.5) & 5.2(5.4) & 3.6(3.7) & 6.3(6.3) & 5.3(5.3) & 4.2(4.2) \\
\cline{2-13}
          & 100   &$\widehat{\bSigma}_{u}^{\rm  thre} $       & 1     & 31.9(32.1) & 20.5(20.8) & 12.6(12.8) & 46.7(48.9) & 29.8(31.4) & 18.0(19.0) & 25.5(25.3) & 17.1(17.2) & 11.0(10.9) \\
          &       &       & 2     & 36.7(36.6) & 23.8(23.9) & 14.2(14.3) & 53.7(56.0) & 34.7(36.2) & 20.5(21.5) & 28.2(28.0) & 19.0(19.1) & 11.9(11.9) \\
          &       &       & 3     & 41.5(41.3) & 27.1(27.0) & 15.9(16.0) & 60.7(63.0) & 39.5(41.0) & 23.0(24.0) & 31.1(30.7) & 21.0(21.0) & 12.9(12.8) \\
          &       &$\widehat{\bSigma}_{u}^{\rm bc, thre} $       & 1     & 5.4(5.4) & 4.5(4.5) & 3.7(3.7) & 6.4(6.5) & 4.9(4.9) & 3.7(3.7) & 7.5(7.4) & 6.3(6.3) & 4.9(4.8) \\
          &       &       & 2     & 5.4(5.4) & 4.5(4.5) & 3.7(3.7) & 6.6(6.7) & 5.0(5.1) & 3.7(3.8) & 7.3(7.2) & 6.1(6.1) & 4.8(4.7) \\
          &       &       & 3     & 5.6(5.6) & 4.6(4.7) & 3.7(3.8) & 6.9(7.0) & 5.2(5.3) & 3.8(3.9) & 7.2(7.1) & 6.0(6.0)  & 4.7(4.6) \\
\cline{2-13}
          & 200   &$\widehat{\bSigma}_{u}^{\rm  thre} $       & 1     & 40.8(40.9) & 26.5(26.6) & 15.5(15.6) & 44.7(45.4) & 29.0(29.4) & 16.8(17.1) & 31.0(30.9) & 21.0(21.1) & 12.8(12.9) \\
          &       &       & 2     & 45.4(45.5) & 29.7(29.7) & 17.2(17.2) & 49.8(50.5) & 32.5(32.9) & 18.6(18.9) & 33.8(33.6) & 22.9(22.9) & 13.8(13.8) \\
          &       &       & 3     & 50.1(50.1) & 32.9(32.8) & 18.8(18.8) & 55.0(55.6) & 35.9(36.3) & 20.4(20.6) & 36.7(36.5) & 24.8(24.8) & 14.7(14.8) \\
          &       &$\widehat{\bSigma}_{u}^{\rm bc, thre} $       & 1     & 7.5(6.2) & 5.4(5.1) & 4.2(4.1) & 8.3(6.4) & 5.7(4.8) & 3.8(3.6) & 8.2(8.2) & 6.9(6.9) & 5.3(5.3) \\
          &       &       & 2     & 8.1(6.3) & 5.7(5.2) & 4.2(4.1) & 8.9(6.5) & 6.1(4.9) & 3.9(3.6) & 8.1(8.1) & 6.7(6.7) & 5.2(5.2) \\
          &       &       & 3     & 8.8(6.5) & 6.0(5.3) & 4.3(4.2) & 9.6(6.7) & 6.6(5.0) & 4.1(3.7) & 8.0(8.0) & 6.7(6.6) & 5.2(5.1) \\

  		\hline
  	   	\multicolumn{3}{c}{Asynchronous Data} &  &  \multicolumn{3}{c}{Model 1}  &  \multicolumn{3}{c}{Model 2}  &  \multicolumn{3}{c}{Model 3}\\
		${252\sigma_e^2}/{\bar{\sigma}^2}$ &   $p $    & Estimators  & $\Delta_K$ & $\lambda=3$     & $\lambda=2$     & $\lambda=1$     & $\lambda=3$     & $\lambda=2$     & $\lambda=1$     & $\lambda=3$     & $\lambda=2$     & $\lambda=1$ \\
		\hline		
  		        0.063 & 50    &$\widehat{\bSigma}_{u}^{\rm  thre} $       & 1     & 6.3(6.3) & 5.1(5.1) & 3.7(3.7) & 9.8(10.1) & 7.3(7.6) & 4.4(4.7) & 9.8(9.7) & 7.3(7.2) & 5.5(5.5) \\
          &       &       & 2     & 6.5(6.4) & 5.2(5.2) & 3.7(3.7) & 9.8(10.3) & 7.5(7.8) & 4.6(4.8) & 9.2(9.2) & 7.0(6.9) & 5.2(5.2) \\
          &       &       & 3     & 6.7(6.7) & 5.4(5.3) & 3.8(3.8) & 10.3(10.7) & 7.8(8.2) & 4.8(5.1) & 8.9(8.8) & 6.8(6.7) & 5.1(5.1) \\
          &       &$\widehat{\bSigma}_{u}^{\rm bc, thre} $       & 1     & 6.2(6.3) & 5.2(5.2) & 3.8(3.8) & 9.8(10.2) & 7.4(7.8) & 4.5(4.7) & 9.7(9.6) & 7.2(7.2) & 5.5(5.5) \\
          &       &       & 2     & 6.4(6.4) & 5.3(5.3) & 3.8(3.8) & 9.9(10.4) & 7.6(8.0) & 4.6(4.9) & 9.1(9.1) & 6.9(6.9) & 5.3(5.3) \\
          &       &       & 3     & 6.7(6.8) & 5.4(5.5) & 3.9(3.9) & 10.4(11.0) & 8.0(8.4) & 4.9(5.2) & 8.7(8.7) & 6.7(6.7) & 5.1(5.1) \\
\cline{2-13}
          & 100   &$\widehat{\bSigma}_{u}^{\rm  thre} $       & 1     & 7.6(7.6) & 5.6(5.6) & 4.2(4.2) & 11.0(11.3) & 7.3(7.4) & 4.9(5.0) & 10.5(10.4) & 8.3(8.2) & 5.9(5.9) \\
          &       &       & 2     & 7.8(7.8) & 5.7(5.7) & 4.3(4.3) & 11.0(11.3) & 7.5(7.7) & 5.0(5.1) & 10.1(10.0) & 7.9(7.8) & 5.7(5.7) \\
          &       &       & 3     & 8.1(8.1) & 5.9(5.9) & 4.4(4.4) & 11.3(11.6) & 7.8(8.0) & 5.2(5.3) & 9.9(9.8) & 7.7(7.6) & 5.6(5.6) \\
          &       &$\widehat{\bSigma}_{u}^{\rm bc, thre} $       & 1     & 7.7(7.7) & 5.7(5.7) & 4.3(4.3) & 11.2(11.6) & 7.5(7.6) & 5.1(5.2) & 10.4(10.3) & 8.2(8.1) & 5.9(5.9) \\
          &       &       & 2     & 8.0(8.0)  & 5.8(5.9) & 4.4(4.4) & 11.3(11.7) & 7.7(7.9) & 5.2(5.3) & 10.0(9.9) & 7.9(7.8) & 5.8(5.8) \\
          &       &       & 3     & 8.3(8.3) & 6.1(6.1) & 4.6(4.6) & 11.8(12.1) & 8.1(8.3) & 5.4(5.5) & 9.8(9.7) & 7.7(7.6) & 5.7(5.7) \\
\cline{2-13}
          & 200   &$\widehat{\bSigma}_{u}^{\rm  thre} $       & 1     & 9.1(9.1) & 6.6(6.6) & 4.8(4.8) & 10.6(11.0) & 7.2(7.4) & 4.7(4.8) & 11.3(11.2) & 8.8(8.7) & 6.3(6.3) \\
          &       &       & 2     & 9.2(9.2) & 6.7(6.7) & 4.8(4.9) & 10.4(10.8) & 7.3(7.4) & 4.7(4.8) & 11.0(11.0) & 8.5(8.5) & 6.2(6.2) \\
          &       &       & 3     & 9.4(9.4) & 6.8(6.8) & 4.9(4.9) & 10.6(10.9) & 7.4(7.5) & 4.7(4.8) & 10.9(10.8) & 8.4(8.3) & 6.1(6.1) \\
          &       &$\widehat{\bSigma}_{u}^{\rm bc, thre} $       & 1     & 9.4(9.4) & 6.9(6.9) & 5.0(5.0) & 11.1(11.6) & 7.6(7.7) & 5.0(5.1) & 11.2(11.2) & 8.7(8.7) & 6.4(6.4) \\
          &       &       & 2     & 9.7(9.6) & 7.0(7.1) & 5.1(5.1) & 11.1(11.5) & 7.7(7.8) & 5.0(5.1) & 11.0(11.0) & 8.5(8.5) & 6.3(6.3) \\
          &       &       & 3     & 10.0(10.0) & 7.3(7.3) & 5.2(5.2) & 11.4(11.7) & 7.9(8.1) & 5.1(5.2) & 10.9(10.9) & 8.4(8.4) & 6.2(6.2) \\
  		\hline
    0.00252 & 50    &$\widehat{\bSigma}_{u}^{\rm  thre} $       & 1     & 61.8(49.2) & 46.7(36.8) & 18.4(14.7) & 109(93.2) & 81.8(69.6) & 31.5(27.0) & 48.2(40.4) & 36.7(30.3) & 16.1(13.3) \\
          &       &       & 2     & 78.9(62.3) & 55.5(43.9) & 22.1(17.4) & 137(118) & 97.2(83.5) & 38.2(32.4) & 57.5(47.8) & 41.5(34.5) & 18.3(15.1) \\
          &       &       & 3     & 95.8(75.9) & 64.3(51.3) & 26.0(20.3) & 165(143) & 113(97.7) & 45.0(38.1) & 67.4(55.5) & 46.4(38.8) & 20.5(16.9) \\
          &       &$\widehat{\bSigma}_{u}^{\rm bc, thre} $       & 1     & 6.9(6.7) & 5.4(5.3) & 3.8(3.8) & 15.1(14.8) & 8.1(8.1) & 4.7(4.8) & 12.7(12.2) & 8.5(8.3) & 5.8(5.7) \\
          &       &       & 2     & 7.2(6.9) & 5.4(5.4) & 3.8(3.8) & 15.8(15.5) & 8.3(8.2) & 4.8(5.0) & 12.9(12.2) & 8.4(8.2) & 5.6(5.5) \\
          &       &       & 3     & 7.6(7.0) & 5.5(5.5) & 3.9(3.9) & 16.6(16.6) & 8.6(8.5) & 5.0(5.3) & 13.2(12.3) & 8.3(8.1) & 5.4(5.4) \\
\cline{2-13}
          & 100   &$\widehat{\bSigma}_{u}^{\rm  thre} $       & 1     & 91.9(89.5) & 46.8(46.9) & 24.7(24.8) & 133(136) & 68.9(72.1) & 35.9(37.6) & 65.3(64.0) & 37.4(37.1) & 20.4(20.3) \\
          &       &       & 2     & 107(104) & 55.7(55.6) & 28.5(28.5) & 156(160) & 82.1(85.4) & 41.6(43.3) & 73.9(72.1) & 42.1(42.3) & 22.6(22.4) \\
          &       &       & 3     & 121(119) & 64.8(64.8) & 32.4(32.2) & 179(183) & 95.3(98.7) & 47.3(49.0) & 82.8(81.1) & 47.0(47.7) & 24.9(24.6) \\
          &       &$\widehat{\bSigma}_{u}^{\rm bc, thre} $       & 1     & 9.0(8.8) & 5.8(5.8) & 4.3(4.3) & 19.8(21.4) & 8.4(8.4) & 5.2(5.3) & 15.5(15.0) & 9.7(9.6) & 6.3(6.3) \\
          &       &       & 2     & 9.9(9.9) & 5.9(5.9) & 4.4(4.4) & 20.8(22.8) & 8.5(8.5) & 5.3(5.4) & 16.1(15.5) & 9.5(9.4) & 6.2(6.1) \\
          &       &       & 3     & 11.3(11.3) & 6.1(6.1) & 4.5(4.5) & 22.0(24.2) & 8.8(8.9) & 5.5(5.6) & 16.9(16.1) & 9.5(9.3) & 6.1(6.0) \\
\cline{2-13}
          & 200   &$\widehat{\bSigma}_{u}^{\rm  thre} $       & 1     & 119(121) & 63.5(64.5) & 31.2(31.4) & 132(135) & 69.8(71.6) & 34.0(34.7) & 81.1(82.0) & 46.7(47.8) & 24.5(24.5) \\
          &       &       & 2     & 132(135) & 72.6(73.2) & 34.9(35.1) & 149(152) & 79.5(81.3) & 38.1(38.8) & 89.9(91.1) & 51.7(52.7) & 26.6(26.7) \\
          &       &       & 3     & 147(150) & 81.8(82.0) & 38.7(38.8) & 167(170) & 89.3(91.1) & 42.2(42.9) & 99.5(101) & 57.1(57.7) & 28.9(28.9) \\
          &       &$\widehat{\bSigma}_{u}^{\rm bc, thre} $       & 1     & 14.3(14.3) & 7.2(6.7) & 5.1(4.9) & 25.5(27.1) & 9.7(9.7) & 5.3(5.1) & 18.9(18.9) & 10.8(10.7) & 6.9(6.8) \\
          &       &       & 2     & 15.6(15.7) & 7.4(6.8) & 5.2(5.0) & 27.2(28.9) & 10.0(9.9) & 5.5(5.2) & 20.1(20.2) & 10.8(10.7) & 6.7(6.7) \\
          &       &       & 3     & 17.1(17.1) & 7.8(6.9) & 5.4(5.1) & 28.8(30.5) & 10.4(10.2) & 5.8(5.3) & 21.4(21.7) & 10.9(10.8) & 6.7(6.7) \\

\hline
  		\hline
  	\end{tabular}
  	\label{Error}
  \end{table}%

   \begin{landscape}
\begin{table}[htbp]
\tiny
	\centering
	\caption{The empirical true positive rates ($\times100$) and false positive rates ($\times100$) with synchronous data when jumps exist and do not exist (in parentheses) based on 1000 repetitions.}\label{TPRFPR_Syn}
		\begin{tabular}{ccccccccccccccccccccc}
	\hline
 \hline
	                &  &  &      &  \multicolumn{6}{c}{Model 1}  &  \multicolumn{6}{c}{Model 2}  \\
                   &    \multicolumn{3}{c}{}      &  \multicolumn{3}{c}{TPR}  &  \multicolumn{3}{c}{FPR}  &  \multicolumn{3}{c}{TPR }  &  \multicolumn{3}{c}{FPR} \\
          $252\sigma_e^2/\bar{\sigma}^2$      &         $p$ & Estimators & $\Delta_K$& $\Delta=3$     & $\Delta=2$     & $\Delta=1$     & $\Delta=3$     & $\Delta=2$     & $\Delta=1$     & $\Delta=3$     & $\Delta=2$     & $\Delta=1$ & $\Delta=3$     & $\Delta=2$     & $\Delta=1$\\
                         \hline
                 0.063 & 50    &$\widehat{\bSigma}_{u}^{\rm  thre} $       & 1     & 100(100) & 100(100) & 100(100) & 0.5(0.5) & 0.9(0.9) & 0.6(0.5) & 100(100) & 100(100) & 100(100) & 0.5(0.5) & 0.9(0.9) & 0.5(0.5) \\
          &       &       & 2     & 100(100) & 100(100) & 100(100) & 0.8(0.8) & 1.2(1.2) & 0.8(0.8) & 100(100) & 100(100) & 100(100) & 0.8(0.8) & 1.2(1.2) & 0.8(0.8) \\
          &       &       & 3     & 100(100) & 100(100) & 100(100) & 1.2(1.2) & 1.6(1.6) & 1.2(1.2) & 100(100) & 100(100) & 100(100) & 1.2(1.2) & 1.6(1.7) & 1.2(1.2) \\
          &       &$\widehat{\bSigma}_{u}^{\rm bc, thre} $       & 1     & 100(100) & 100(100) & 100(100) & 0.5(0.6) & 1.0(1.0)  & 0.6(0.6) & 100(100) & 100(100) & 100(100) & 0.5(0.6) & 1.0(1.0)  & 0.6(0.6) \\
          &       &       & 2     & 100(100) & 100(100) & 100(100) & 0.8(0.9) & 1.4(1.3) & 0.9(0.9) & 100(100) & 100(100) & 100(100) & 0.8(0.9) & 1.3(1.4) & 0.9(0.9) \\
          &       &       & 3     & 100(100) & 100(100) & 100(100) & 1.3(1.3) & 1.8(1.8) & 1.4(1.4) & 100(100) & 100(100) & 100(100) & 1.3(1.4) & 1.8(1.9) & 1.4(1.4) \\
\cline{2-16}
          & 100   &$\widehat{\bSigma}_{u}^{\rm thre} $       & 1     & 100(100) & 100(100) & 100(100) & 0.6(0.6) & 0.5(0.4) & 0.6(0.5) & 100(100) & 100(100) & 100(100) & 0.6(0.6) & 0.5(0.5) & 0.6(0.6) \\
          &       &       & 2     & 100(100) & 100(100) & 100(100) & 0.8(0.8) & 0.6(0.6) & 0.7(0.7) & 100(100) & 100(100) & 100(100) & 0.8(0.8) & 0.7(0.6) & 0.7(0.7) \\
          &       &       & 3     & 100(100) & 100(100) & 100(100) & 1.0(1.0)  & 0.9(0.9) & 1.0(1.0)  & 100(100) & 100(100) & 100(100) & 1.0(1.0)  & 0.9(0.9) & 1.0(1.0) \\
          &       &$\widehat{\bSigma}_{u}^{\rm bc, thre} $       & 1     & 100(100) & 100(100) & 100(100) & 0.7(0.7) & 0.5(0.5) & 0.7(0.6) & 100(100) & 100(100) & 100(100) & 0.7(0.7) & 0.5(0.5) & 0.7(0.7) \\
          &       &       & 2     & 100(100) & 100(100) & 100(100) & 0.9(0.9) & 0.7(0.7) & 0.9(0.9) & 100(100) & 100(100) & 100(100) & 0.9(0.9) & 0.7(0.7) & 0.9(0.9) \\
          &       &       & 3     & 100(100) & 100(100) & 100(100) & 1.2(1.2) & 1.1(1.1) & 1.2(1.2) & 100(100) & 100(100) & 100(100) & 1.2(1.2) & 1.1(1.1) & 1.2(1.2) \\
          \cline{2-16}
          & 200   &$\widehat{\bSigma}_{u}^{\rm  thre} $       & 1     & 100(100) & 100(100) & 100(100) & 0.5(0.5) & 0.5(0.5) & 0.5(0.5) & 100(100) & 100(100) & 100(100) & 0.5(0.5) & 0.5(0.5) & 0.5(0.5) \\
          &       &       & 2     & 100(100) & 100(100) & 100(100) & 0.6(0.6) & 0.6(0.6) & 0.6(0.6) & 100(100) & 100(100) & 100(100) & 0.6(0.6) & 0.6(0.6) & 0.6(0.6) \\
          &       &       & 3     & 100(100) & 100(100) & 100(100) & 0.7(0.7) & 0.7(0.7) & 0.7(0.7) & 100(100) & 100(100) & 100(100) & 0.7(0.7) & 0.7(0.7) & 0.7(0.7) \\
          &       &$\widehat{\bSigma}_{u}^{\rm bc, thre} $       & 1     & 100(100) & 100(100) & 100(100) & 0.6(0.6) & 0.5(0.5) & 0.6(0.6) & 100(100) & 100(100) & 100(100) & 0.6(0.6) & 0.6(0.5) & 0.6(0.6) \\
          &       &       & 2     & 100(100) & 100(100) & 100(100) & 0.7(0.8) & 0.7(0.7) & 0.8(0.8) & 100(100) & 100(100) & 100(100) & 0.8(0.8) & 0.7(0.7) & 0.8(0.8) \\
          &       &       & 3     & 100(100) & 100(100) & 100(100) & 0.9(0.9) & 0.8(0.8) & 0.9(0.9) & 100(100) & 100(100) & 100(100) & 0.9(0.9) & 0.9(0.8) & 0.9(0.9) \\
\hline
    0.00252 & 50    &$\widehat{\bSigma}_{u}^{\rm  thre} $       & 1     & 100(100) & 100(100) & 100(100) & 5.6(4.4) & 7.6(6.2) & 4.0(3.0)  & 100(100) & 100(100) & 100(100) & 11.0(9.8) & 13.4(12.1) & 8.3(7.0) \\
          &       &       & 2     & 100(100) & 100(100) & 100(100) & 8.2(6.8) & 10.0(8.6) & 5.8(4.5) & 100(100) & 100(100) & 100(100) & 14.2(13.1) & 16.1(15.0) & 10.9(9.5) \\
          &       &       & 3     & 100(100) & 100(100) & 100(100) & 10.9(9.7) & 12.5(11.1) & 8.0(6.6) & 100(100) & 100(100) & 100(100) & 17.2(16.4) & 18.7(17.8) & 13.6(12.4) \\
          &       &$\widehat{\bSigma}_{u}^{\rm bc, thre} $       & 1     & 100(100) & 100(100) & 100(100) & 0.5(0.5) & 0.8(0.8) & 0.5(0.5) & 100(100) & 100(100) & 100(100) & 0.5(0.5) & 0.8(0.8) & 0.5(0.5) \\
          &       &       & 2     & 100(100) & 100(100) & 100(100) & 0.7(0.7) & 1.1(1.1) & 0.8(0.8) & 100(100) & 100(100) & 100(100) & 0.7(0.7) & 1.0(1.1) & 0.7(0.8) \\
          &       &       & 3     & 100(100) & 100(100) & 100(100) & 1.1(1.1) & 1.5(1.5) & 1.2(1.2) & 100(100) & 100(100) & 100(100) & 1.1(1.1) & 1.5(1.5) & 1.1(1.2) \\
\cline{2-16}
          & 100   &$\widehat{\bSigma}_{u}^{\rm  thre} $       & 1     & 100(100) & 100(100) & 100(100) & 4.1(4.0) & 3.0(2.8) & 2.8(2.7) & 100(100) & 100(100) & 100(100) & 7.5(7.3) & 6.0(5.8) & 5.6(5.3) \\
          &       &       & 2     & 100(100) & 100(100) & 100(100) & 5.4(5.2) & 4.2(4.0) & 3.8(3.6) & 100(99.9) & 100(100) & 100(100) & 8.9(8.8) & 7.5(7.2) & 6.9(6.6) \\
          &       &       & 3     & 100(100) & 100(100) & 100(100) & 6.7(6.6) & 5.5(5.3) & 4.9(4.7) & 99.9(99.8) & 100(100) & 100(100) & 10.3(10.2) & 9.0(8.7) & 8.2(8.0) \\
          &       &$\widehat{\bSigma}_{u}^{\rm bc, thre} $       & 1     & 100(100) & 100(100) & 100(100) & 0.5(0.5) & 0.4(0.4) & 0.5(0.5) & 100(100) & 100(100) & 100(100) & 0.5(0.5) & 0.4(0.4) & 0.5(0.5) \\
          &       &       & 2     & 100(100) & 100(100) & 100(100) & 0.7(0.7) & 0.6(0.6) & 0.7(0.7) & 100(100) & 100(100) & 100(100) & 0.7(0.7) & 0.6(0.6) & 0.7(0.7) \\
          &       &       & 3     & 100(100) & 100(100) & 100(100) & 0.9(0.9) & 0.8(0.8) & 1.0(1.0)  & 100(100) & 100(100) & 100(100) & 0.9(0.9) & 0.8(0.8) & 1.0(1.0) \\
\cline{2-16}
          & 200   &$\widehat{\bSigma}_{u}^{\rm  thre} $       & 1     & 100(100) & 100(100) & 100(100) & 2.8(2.8) & 2.2(2.2) & 2.0(2.0)  & 99.8(99.8) & 100(100) & 100(100) & 4.6(4.7) & 3.9(4.0) & 3.6(3.6) \\
          &       &       & 2     & 100(100) & 100(100) & 100(100) & 3.3(3.4) & 2.7(2.8) & 2.4(2.5) & 99.8(99.7) & 100(100) & 100(100) & 5.2(5.3) & 4.6(4.6) & 4.1(4.2) \\
          &       &       & 3     & 100(100) & 100(100) & 100(100) & 3.9(3.9) & 3.3(3.4) & 2.9(3.0) & 99.6(99.6) & 100(99.9) & 100(100) & 5.8(5.8) & 5.2(5.2) & 4.7(4.7) \\
          &       &$\widehat{\bSigma}_{u}^{\rm bc, thre} $       & 1     & 100(100) & 100(100) & 100(100) & 0.5(0.5) & 0.4(0.4) & 0.5(0.5) & 100(100) & 100(100) & 100(100) & 0.5(0.5) & 0.4(0.4) & 0.5(0.5) \\
          &       &       & 2     & 100(100) & 100(100) & 100(100) & 0.6(0.6) & 0.5(0.5) & 0.6(0.6) & 100(100) & 100(100) & 100(100) & 0.6(0.6) & 0.5(0.5) & 0.6(0.6) \\
          &       &       & 3     & 100(100) & 100(100) & 100(100) & 0.7(0.7) & 0.7(0.7) & 0.8(0.8) & 100(100) & 100(100) & 100(100) & 0.7(0.7) & 0.7(0.7) & 0.8(0.8) \\

           \hline
 \hline

	\end{tabular}

\end{table}%

   \end{landscape}

   \begin{landscape}
\begin{table}[htbp]
\tiny
	\centering
	\caption{The empirical true positive rates ($\times100$) and false positive rates ($\times100$) in different settings with asynchronous data when jumps exist and do not exist (in parentheses) based on 1000 repetitions.}\label{TPRFPR_Asyn}
		\begin{tabular}{ccccccccccccccccccccc}
	\hline
 \hline
	                &  &  &      &  \multicolumn{6}{c}{Model 1}  &  \multicolumn{6}{c}{Model 2}  \\
                   &    \multicolumn{3}{c}{}      &  \multicolumn{3}{c}{TPR}  &  \multicolumn{3}{c}{FPR}  &  \multicolumn{3}{c}{TPR }  &  \multicolumn{3}{c}{FPR} \\
          $252\sigma_e^2/\bar{\sigma}^2$      &         $p$ & Estimators & $\Delta_K$& $\lambda=3$ &$\lambda=2$& $\lambda=1$ & $\lambda=3$ &$\lambda=2$& $\lambda=1$ & $\lambda=3$ &$\lambda=2$& $\lambda=1$ &$\lambda=3$ &$\lambda=2$& $\lambda=1$ \\
                         \hline
	     0.063 & 50    &$\widehat{\bSigma}_{u}^{\rm  thre} $       & 1     & 100(100) & 100(100) & 100(100) & 0.6(0.6) & 0.9(0.9) & 0.4(0.4) & 99.3(98.9) & 100(100) & 100(100) & 0.6(0.6) & 0.9(0.9) & 0.5(0.4) \\
          &       &       & 2     & 100(100) & 100(100) & 100(100) & 0.9(0.9) & 1.2(1.2) & 0.7(0.7) & 99.7(99.5) & 100(100) & 100(100) & 0.9(0.9) & 1.3(1.2) & 0.7(0.7) \\
          &       &       & 3     & 100(100) & 100(100) & 100(100) & 1.3(1.3) & 1.6(1.6) & 1.1(1.1) & 99.8(99.7) & 100(100) & 100(100) & 1.4(1.3) & 1.7(1.7) & 1.1(1.1) \\
          &       &$\widehat{\bSigma}_{u}^{\rm bc, thre} $       & 1     & 100(100) & 100(100) & 100(100) & 0.6(0.6) & 1.0(1.0)  & 0.5(0.5) & 99.3(98.9) & 100(100) & 100(100) & 0.6(0.6) & 1.0(1.0)  & 0.5(0.5) \\
          &       &       & 2     & 100(100) & 100(100) & 100(100) & 1.0(1.0)  & 1.4(1.4) & 0.8(0.8) & 99.7(99.4) & 100(100) & 100(100) & 1.0(1.0)  & 1.4(1.4) & 0.8(0.8) \\
          &       &       & 3     & 100(100) & 100(100) & 100(100) & 1.5(1.5) & 1.9(1.9) & 1.3(1.3) & 99.8(99.7) & 100(100) & 100(100) & 1.5(1.5) & 1.9(1.9) & 1.3(1.3) \\
\cline{2-16}
          & 100   &$\widehat{\bSigma}_{u}^{\rm  thre} $       & 1     & 100(100) & 100(100) & 100(100) & 0.6(0.6) & 0.4(0.4) & 0.5(0.5) & 99.4(99.1) & 100(100) & 100(100) & 0.6(0.6) & 0.4(0.4) & 0.5(0.5) \\
          &       &       & 2     & 100(100) & 100(100) & 100(100) & 0.8(0.8) & 0.6(0.6) & 0.7(0.7) & 99.7(99.5) & 100(100) & 100(100) & 0.9(0.9) & 0.7(0.7) & 0.7(0.7) \\
          &       &       & 3     & 100(100) & 100(100) & 100(100) & 1.1(1.1) & 0.9(0.9) & 1.0(1.0)  & 99.8(99.7) & 100(100) & 100(100) & 1.1(1.1) & 0.9(0.9) & 1.0(1.0) \\
          &       &$\widehat{\bSigma}_{u}^{\rm bc, thre} $       & 1     & 100(100) & 100(100) & 100(100) & 0.7(0.7) & 0.5(0.5) & 0.6(0.6) & 99.4(99.0) & 100(100) & 100(100) & 0.7(0.7) & 0.5(0.5) & 0.6(0.6) \\
          &       &       & 2     & 100(100) & 100(100) & 100(100) & 1.0(1.0)  & 0.8(0.8) & 0.9(0.9) & 99.7(99.4) & 100(100) & 100(100) & 1.0(1.0)  & 0.8(0.8) & 0.9(0.8) \\
          &       &       & 3     & 100(100) & 100(100) & 100(100) & 1.3(1.3) & 1.1(1.1) & 1.2(1.2) & 99.8(99.6) & 100(100) & 100(100) & 1.3(1.3) & 1.1(1.1) & 1.2(1.2) \\
\cline{2-16}
          & 200   &$\widehat{\bSigma}_{u}^{\rm  thre} $       & 1     & 100(100) & 100(100) & 100(100) & 0.5(0.5) & 0.5(0.5) & 0.5(0.5) & 99.4(98.8) & 100(100) & 100(100) & 0.5(0.5) & 0.5(0.5) & 0.5(0.5) \\
          &       &       & 2     & 100(100) & 100(100) & 100(100) & 0.6(0.6) & 0.6(0.6) & 0.6(0.6) & 99.6(99.2) & 100(100) & 100(100) & 0.6(0.6) & 0.6(0.6) & 0.6(0.6) \\
          &       &       & 3     & 100(100) & 100(100) & 100(100) & 0.7(0.7) & 0.7(0.7) & 0.7(0.7) & 99.7(99.4) & 100(100) & 100(100) & 0.7(0.7) & 0.7(0.7) & 0.7(0.7) \\
          &       &$\widehat{\bSigma}_{u}^{\rm bc, thre} $       & 1     & 100(100) & 100(100) & 100(100) & 0.6(0.6) & 0.6(0.6) & 0.6(0.6) & 99.3(98.6) & 100(100) & 100(100) & 0.6(0.6) & 0.6(0.6) & 0.6(0.6) \\
          &       &       & 2     & 100(100) & 100(100) & 100(100) & 0.8(0.8) & 0.7(0.7) & 0.8(0.8) & 99.6(99.1) & 100(100) & 100(100) & 0.8(0.8) & 0.7(0.7) & 0.8(0.8) \\
          &       &       & 3     & 100(100) & 100(100) & 100(100) & 1.0(1.0)  & 0.9(0.9) & 0.9(0.9) & 99.7(99.3) & 100(100) & 100(100) & 1.0(1.0)  & 0.9(0.9) & 0.9(0.9) \\
\hline
    0.00252 & 50    &$\widehat{\bSigma}_{u}^{\rm  thre} $       & 1     & 99.8(100) & 100(100) & 100(100) & 6.6(5.6) & 8.9(7.9) & 4.7(3.6) & 88.4(89.8) & 97.0(97.9) & 100(100) & 12.7(12.0) & 15.3(14.6) & 9.8(8.5) \\
          &       &       & 2     & 98.4(99.6) & 100(100) & 100(100) & 9.3(8.4) & 11.6(10.7) & 7.1(5.7) & 87.8(89.9) & 96.6(97.4) & 100(100) & 15.8(15.4) & 18.1(17.6) & 12.9(11.7) \\
          &       &       & 3     & 95.5(98.7) & 100(100) & 100(100) & 11.8(11.2) & 14.1(13.4) & 9.8(8.4) & 87.1(89.7) & 96.3(97.0) & 100(100) & 18.4(18.3) & 20.5(20.4) & 15.9(14.9) \\
          &       &$\widehat{\bSigma}_{u}^{\rm bc, thre} $       & 1     & 100(100) & 100(100) & 100(100) & 0.2(0.2) & 0.5(0.5) & 0.4(0.4) & 91.8(92.9) & 99.8(99.8) & 100(100) & 0.2(0.2) & 0.4(0.5) & 0.4(0.4) \\
          &       &       & 2     & 100(100) & 100(100) & 100(100) & 0.3(0.3) & 0.6(0.6) & 0.6(0.6) & 91.2(92.7) & 99.8(99.8) & 100(100) & 0.3(0.3) & 0.6(0.6) & 0.6(0.6) \\
          &       &       & 3     & 99.9(100) & 100(100) & 100(100) & 0.3(0.4) & 0.7(0.8) & 0.9(0.9) & 90.0(92.3) & 99.7(99.8) & 100(100) & 0.3(0.4) & 0.7(0.8) & 0.9(0.9) \\
\cline{2-16}
          & 100   &$\widehat{\bSigma}_{u}^{\rm  thre} $       & 1     & 96.2(97.4) & 100(100) & 100(100) & 4.5(4.4) & 3.6(3.5) & 3.6(3.4) & 86.1(80.9) & 98.1(97.4) & 100(100) & 8.2(8.1) & 7.0(6.9) & 6.8(6.6) \\
          &       &       & 2     & 93.4(94.8) & 100(100) & 100(100) & 5.6(5.6) & 5.0(4.9) & 4.8(4.6) & 85.4(80.2) & 97.9(97.0) & 100(100) & 9.3(9.3) & 8.6(8.4) & 8.2(8.1) \\
          &       &       & 3     & 89.9(91.8) & 100(99.8) & 100(100) & 6.6(6.6) & 6.3(6.2) & 6.1(6.0) & 84.3(79.3) & 97.7(96.6) & 100(100) & 10.3(10.3) & 10.0(9.9) & 9.7(9.5) \\
          &       &$\widehat{\bSigma}_{u}^{\rm bc, thre} $       & 1     & 99.9(99.9) & 100(100) & 100(100) & 0.2(0.2) & 0.2(0.2) & 0.4(0.4) & 87.5(84.5) & 99.8(99.6) & 100(100) & 0.2(0.2) & 0.2(0.2) & 0.4(0.4) \\
          &       &       & 2     & 99.8(99.8) & 100(100) & 100(100) & 0.2(0.2) & 0.3(0.3) & 0.6(0.6) & 86.5(83.4) & 99.8(99.7) & 100(100) & 0.2(0.2) & 0.3(0.3) & 0.6(0.6) \\
          &       &       & 3     & 99.5(99.5) & 100(100) & 100(100) & 0.2(0.2) & 0.4(0.4) & 0.8(0.8) & 85.1(82.0) & 99.8(99.7) & 100(100) & 0.2(0.2) & 0.4(0.4) & 0.8(0.8) \\
\cline{2-16}
          & 200   &$\widehat{\bSigma}_{u}^{\rm  thre} $       & 1     & 89.3(89.6) & 100(99.9) & 100(100) & 2.7(2.8) & 2.6(2.6) & 2.5(2.6) & 76.1(70.5) & 98.0(97.0) & 100(99.9) & 4.7(4.7) & 4.4(4.5) & 4.3(4.4) \\
          &       &       & 2     & 85.6(86.0) & 99.9(99.8) & 100(100) & 3.1(3.2) & 3.1(3.2) & 3.0(3.1) & 73.9(68.4) & 97.9(96.9) & 100(99.9) & 5.0(5.1) & 5.0(5.1) & 4.9(5.0) \\
          &       &       & 3     & 82.4(82.7) & 99.5(99.7) & 100(100) & 3.5(3.5) & 3.6(3.7) & 3.6(3.7) & 71.9(66.6) & 97.7(96.7) & 99.9(99.8) & 5.4(5.5) & 5.6(5.6) & 5.5(5.6) \\
          &       &$\widehat{\bSigma}_{u}^{\rm bc, thre} $       & 1     & 99.1(99.1) & 100(100) & 100(100) & 0.1(0.1) & 0.2(0.2) & 0.4(0.4) & 76.1(70.9) & 99.2(98.6) & 100(100) & 0.1(0.1) & 0.2(0.2) & 0.4(0.4) \\
          &       &       & 2     & 98.4(98.3) & 100(100) & 100(100) & 0.1(0.1) & 0.2(0.2) & 0.5(0.5) & 73.7(68.5) & 99.2(98.5) & 100(100) & 0.1(0.1) & 0.2(0.2) & 0.5(0.5) \\
          &       &       & 3     & 97.4(97.3) & 100(100) & 100(100) & 0.1(0.1) & 0.3(0.3) & 0.6(0.6) & 71.4(66.3) & 99.1(98.4) & 100(100) & 0.1(0.1) & 0.3(0.3) & 0.6(0.6) \\
     \hline
    \hline

	\end{tabular}
\end{table}
   \end{landscape}

\subsection{Real data analysis}
We analyzed a real high-dimensional data set, studying the statistical properties of microstructure noises that contaminate the trading prices (log-prices) of the constituent stocks of S\&P 500. Intra-day tick-by-tick trading data on two days, November  4 and 22, 2016,  were downloaded from the TAQ database.

Besides the prices themselves,  the  Global Industry Classification Standard (GICS) codes are available to classify the companies in S\&P 500\footnote{The code is 8-digits and each company has its unique code.  Digits 1-2 of the code describe  the company's sector; digits 3-4 describe the industry group; digits 5-6 describe the industry; digits 7-8 describe the sub-industry.}.
 Based on the GICS codes,  there are 36, 27, 71, 84, 36, 58, 64, 65, 5, 28, and 26 companies respectively belonging to the 11 different  sectors -- Energy (E), Materials (M), Industrials (I),   Consumer Discretionary (D), Consumer Staples (S), Health Care (H), Financial (F), Information Technology (T), Telecommunication Services (C), Utilities (U), and Real Estate (R). Since there are only 5 companies belonging to Telecommunication Services, we therefore combined the companies  belonging to the  Information Technology  and Telecommunication Services together and denoted them as `T'.
Our analysis does not assume any information from the GICS classifications; we use it  for validating and interpreting the outcomes from applying the proposed method.

\begin{figure}	
	\begin{center}
	\caption{Magnitudes of the elements in the estimated correlation matrix of the microstructure noises on November 4 and 22, 2016}
	\includegraphics[height=6.2in,
	width=7.2in] {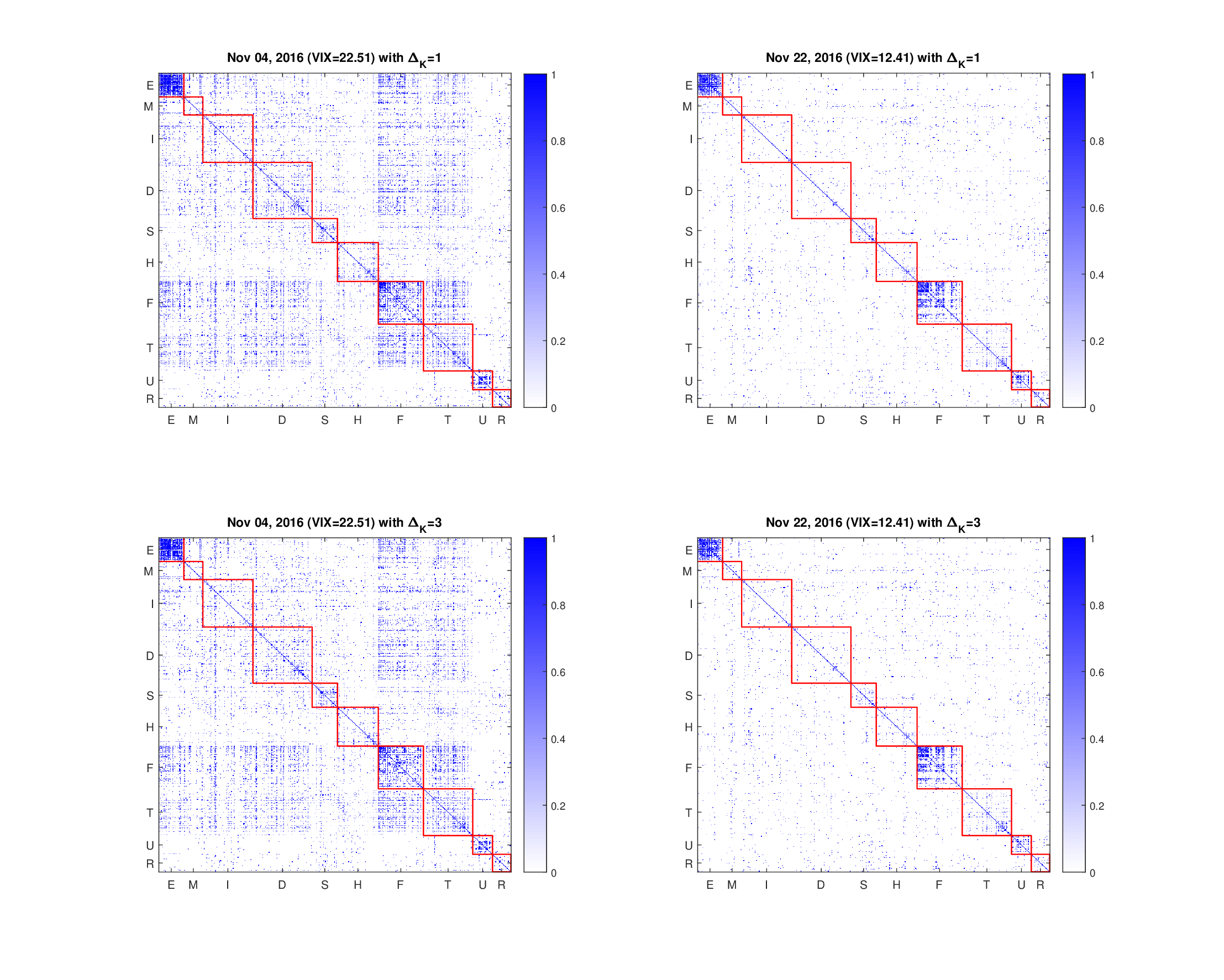}\label{Emprical Graph 1}
	\end{center}
	Note: The graph describes the magnitudes of the elements in the estimated correlation matrices of the microstructure noises. Different colors denote  different values of the pairwise correlations. The red squares along the diagonal denote the sectors---Energy (E), Materials (M), Industrials (I),   Consumer Discretionary (D), Consumer Staples (S), Health Care (H), Financial (F), Information Technology  and Telecommunication Services (T), Utilities (U), and Real Estate (R).
\end{figure}

Upon applying the proposed methods,  we report in Figure \ref{Emprical Graph 1} the magnitudes of the elements in the  correlation matrices of the microstructure noises  
estimated from $\widehat{\bSigma}_u^{{\rm bc, thre}}$  in (\ref{eq:praest}),  respectively for November  4 and 22, 2016 with tuning parameters suggested in Section \ref{sec:tuning}.
Here the companies are sorted by the categories defined by the GICS codes.
The red blocks along the diagonal in Figure \ref{Emprical Graph 1} represent the industrial classifications according to the digits 1-2 of GICS codes.
Hence, we can examine both the within- and between-block correlations as revealed by our estimator.

We remark with some interesting findings from Figure \ref{Emprical Graph 1}.
Overall, we can see that the estimated correlation matrices is sparse with many components  estimated as zero,
indicating that our approach  achieved  the goal of parsimonious  covariance estimations that can (i)  effectively identify  nonzero components, and (ii) support providing meaningful interpretations. 
More specifically, we see that the correlations differ substantially on these two days; and such difference is related to the level of the Chicago Board Options Exchange Volatility Index (VIX),  a popular measure of the stock market's expectation of volatility implied by S\&P 500 index options.
On November 4 when the VIX level was higher,  and the overall correlation level between different components of the microstructure noise is also found to be higher than that on November 22.
Upon examining the  within- and between- category correlations with categories defined by the GICS codes, we see that the correlations within each industrial sector are clear, especially for the Energy and Financial sectors. In contrast, the correlations between different industrial sector are much weaker.
 Meanwhile, we observe that the  correlation estimations have no substantial difference between the cases with $\Delta_K=1$ and 3, an indication that our method is not  sensitive for the choice of $\Delta_K$.

These findings suggest us that it is practically meaningful by studying the high-dimensional statistical properties of the noises. For example,  the between-day difference in correlations is helpful to understand the changes in the market sentiment  under different market conditions.   Furthermore, an interesting feature is the pattern found in the within- and between- industrial sector  correlations, which is seen interrelated to the market conditions.
Broad questions include how the correlations between noises vary associated with the prices and/or the volatility of different assets,  how the sparse  covariance matrix of the noises can help in solving practical problems, and so on.  
Supported by our new methods, we expect more future investigations along this line.



\section{Discussion}
\label{s5}

In this paper, we consider estimating the covariance matrix of the high-dimensional noise in high-frequency data. We propose an estimator with appropriate localization and thresholding to achieve the minimax optimal convergence rates under two kinds of loss. Although all theoretical properties of the proposed estimator are derived under the continuous-time model (\ref{eq:model}), the method developed in this paper could be applied to other types of process $\bX_t$, such as the smooth ones typically encountered in the functional data
literature. The key property that makes our method work is the continuity of the underlying
process $\bX_t$, but the convergence rate of the proposed estimator depends on more specific assumptions,
such as those implied by the model (\ref{eq:model}). On the other hand, Assumption \ref{as:space} can be replaced by a weaker assumption that $\min_{i,j\in[p]}n_{i,j}\rightarrow\infty$ and $\max_{i,j\in[p]}\max_{k\in[n_{i,j}-1]}\Delta t_{i,j,k}\rightarrow0$ as $n\rightarrow\infty$. If we write $n_{\min}=\min_{i,j\in[p]}n_{i,j}$ and assume $\max_{i,j\in[p]}\max_{k\in[n_{i,j}-1]}\Delta t_{i,j,k}\asymp n_{\min}^{-\epsilon}$ for some $\epsilon\in(0,1]$, with suitable selection of $K$, Theorems \ref{pn:1} and \ref{tm:3} in Section \ref{sec:theoryfU} still hold with replacing $n_*$ by $n_{\min}^{g(\epsilon)}$, where $g(\epsilon)\in(0,1]$ is a function of $\epsilon$. More specifically, if $\epsilon=1$, then $g(1)=1$. However, whether such rates are minimax optimal under the associated losses or not is unclear.

In the analysis of this study, we focus on noises with  homoscedastic covariance matrix, so that our target is ${\rm Cov}(\bU_{t_k})\equiv\bSigma_u$ for each $k\in[n]$. As an interesting problem, our framework can be extensively developed to handle time-varying heteroscedastic noise where the target covariance is time-dependent ${\rm Cov}(\bU_{t_{s}})$ at a given time point $t_{s}$ instead. Denote by $\circ$ the Hadamard product. Assume $\bU_{t_k}=\bgamma_{t_k}\circ\bchi_{t_k}$ with $\bgamma_{t_k}=(\gamma_{1,t_k},\ldots,\gamma_{p,t_k})^\top$ and $\bchi_{t_k}=(\chi_{1,t_k},\ldots,\chi_{p,t_k})^\top$, where $\gamma_{1,t},\ldots,\gamma
_{p,t}$ are $p$ nonnegative continuous-time processes, and each sequence $\{\chi_{j,t_k}\}_{k\geq1}$ is $\alpha$-mixing. Without loss of generality, we assume $\mathbb{E}(\chi_{j,t_k})\equiv0$ and $\mathbb{E}(\chi_{j,t_k}^2)=1$ for each $j\in[p]$ and $k\in[n]$. Recall $\mathcal{G}_i\cap\mathcal{G}_j=\{t_{i,j,1},\ldots,t_{i,j,n_{i,j}}\}$ with $t_{i,j,1}<\cdots<t_{i,j,n_{i,j}}$ for any $i,j\in[p]$. 
We outline a framework  as follows.

Step 1.  Given some $\tilde{\xi}=o(1)$, define $\mathcal{G}_i^{(s)}=\{t_{i,\ell}\in\mathcal{G}_i:|t_{i,\ell}-t_{s}|\leq \tilde{\xi}\}$ and $n_{i,j}^{(s)}=|\mathcal{G}_i^{(s)}\cap\mathcal{G}_j^{(s)}|$. For each $t_{i,j,k}\in \mathcal{G}_i^{(s)}\cap\mathcal{G}_j^{(s)}$, write $S_{i,j,k}^{(s)}=\{t_{i,j,\ell}\in\mathcal{G}_i^{(s)}\cap\mathcal{G}_j^{(s)}:K\leq |\ell-k|\leq K+\Delta_K\}$ and $N_{i,j,k}^{(s)}=|S_{i,j,k}^{(s)}|$. We then define $\hat{\sigma}_{u,i,j}^{(s)}$ in the same manner as $\hat{\sigma}_{u,i,j}$ in \eqref{eq:estasyn} with replacing $(n_{i,j},N_{i,j,k},S_{i,j,k})$ by $(n_{i,j}^{(s)},N_{i,j,k}^{(s)},S_{i,j,k}^{(s)})$. Under the independence between the process $\bgamma_t$ and the sequence $\{\bchi_{t_k}\}_{k\geq1}$, following our current technical arguments, we have $\hat{\sigma}_{u,i,j}^{(s)}\rightarrow \gamma_{i,t_{s}}\gamma_{j,t_{s}}\mathbb{E}(\chi_{i,t_{s}}\chi_{j,t_{s}})$ in probability with suitable selection of $\tilde{\xi}$.

Step 2.  Given some integers $\tilde{K}\geq1$ and $h\geq1$, if $\mathbb{E}(\gamma_{i,t_k}\gamma_{j,t_k})$ and $\mathbb{E}(\chi_{i,t_k}\chi_{j,t_k})$ are slowly varying with $k$, we know $\{\tilde{N}_{i,j}^{(s)}\}^{-1}\sum_{b=-\tilde{K}}^{\tilde{K}}\hat{\sigma}_{u,i,j}^{(s+bh)}$ with $\tilde{N}_{i,j}^{(s)}=|\{s+bh\in[n]:|b|\leq\tilde{K}\}|$ will provide a consistent estimator for ${\rm Cov}(U_{i,t_{s}},U_{j,t_{s}})$ under some regularity conditions with suitable selections of $\tilde{K}$ and $h$.

Clearly, this problem differs substantially from our current investigation.
The technical analysis of such estimator would require an extensive framework including additional assumptions on $\{\bgamma_{t_k}\}_{k\geq1}$ and $\{\bchi_{t_k}\}_{k\geq1}$ which are beyond the scope of this study.  We plan to carefully investigate this  problem in a future project.


\section{Proofs}\label{proof}

In the sequel, we use $C$ to denote a generic positive finite universal constant that may be different in different uses.


\subsection{Proof of Theorem \ref{pn:1}}\label{se:pftm1}

For any $k\in[n_{i,j}]$, let $S_{i,j,k}=\{t_{i, j, \ell}:K\leq|\ell-k|\leq K+\Delta_K \}$. For any $i,j\in[p]$, we have that
\begin{align*}
\hat{\sigma}_{u,i,j}-\sigma_{u,i,j}=&~\underbrace{\frac{1}{2n_{i,j}}\sum_{k=1}^{n_{i,j}}\frac{1}{N_{i,j,k}}\sum_{t_{i,j,\ell}\in S_{i,j,k}}(U_{i,t_{i,j,\ell}}-U_{i,t_{i,j,k}})(U_{j,t_{i,j,\ell}}-U_{j,t_{i,j,k}})-\sigma_{u,i,j}}_{\textrm{I}(i,j)}\\
&+\underbrace{\frac{1}{2n_{i,j}}\sum_{k=1}^{n_{i,j}}\frac{1}{N_{i,j,k}}\sum_{t_{i,j,\ell}\in S_{i,j,k}}(X_{i,t_{i,j,\ell}}-X_{i,t_{i,j,k}})(X_{j,t_{i,j,\ell}}-X_{j,t_{i,j,k}})}_{\textrm{II}(i,j)}\\
&+\underbrace{\frac{1}{2n_{i,j}}\sum_{k=1}^{n_{i,j}}\frac{1}{N_{i,j,k}}\sum_{t_{i,j,\ell}\in S_{i,j,k}}(X_{i,t_{i,j,\ell}}-X_{i,t_{i,j,k}})(U_{j,t_{i,j,\ell}}-U_{j,t_{i,j,k}})}_{\textrm{III}(i,j)}\\
&+\underbrace{\frac{1}{2n_{i,j}}\sum_{k=1}^{n_{i,j}}\frac{1}{N_{i,j,k}}\sum_{t_{i,j,\ell}\in S_{i,j,k}}(U_{i,t_{i,j,\ell}}-U_{i,t_{i,j,k}})(X_{j,t_{i,j,\ell}}-X_{j,t_{i,j,k}})}_{\textrm{IV}(i,j)}\,.
\end{align*}
Define $\xi=\max_{i,j\in[p]}\max_{ k\in[n_{i,j}]}\max_{t_{i,j,\ell}\in S_{i,j,k}}|t_{i,j,\ell}-t_{i,j,k}|$. To prove Theorem \ref{pn:1}, we need the following three propositions whose proofs are given in Sections \ref{sec:pfpn1}--\ref{sec:pfpn3}, respectively.
\begin{proposition}\label{laa:1}
Under Assumptions {\rm\ref{as:alphmix}}--{\rm\ref{as:moment}}, we have that
\begin{align*}
\max_{i,j\in[p]}\mathbb{P}\{|{\rm I}(i,j)|>v\}\lesssim&~\{1+n_*(K+L_n)^{-1}\rho^{-1}v^2\}^{-\rho/2}+v^{-1}\exp\{-C(n_{*}L_n^{-1}\rho^{-1}v)^{\varphi/(\varphi+1)}\}\\
&+v^{-1}\exp(-Cn_{*}K^{-1}\rho^{-1}v)
\end{align*}
for any $v\gg \exp(-CL_n^{-\varphi}K^{\varphi})$ and $\rho\geq 1$.
Furthermore, it holds that  $\max_{ i,j\in[p]}\mathbb{E}\{|{\rm I}(i,j)|^{m}\}\lesssim 1$ for any fixed positive integer $m$.
\end{proposition}

\begin{rek}\label{rek:1}
As shown in Section {\rm\ref{sec:pfpn1}}, the upper bound stated in Proposition {\rm \ref{laa:1}} holds for any $v\gg \delta:=\max_{i,j\in[p]}|n_{i,j}^{-1}\sum_{k=1}^{n_{i,j}}\mathbb{E}(N_{i,j,k}^{-1}U_{i,t_{i,j,k}}\sum_{t_{i,j,\ell}\in S_{i,j,k}}U_{j,t_{i,j,\ell}})|$. Since $\delta\lesssim \exp(-CL_n^{-\varphi}K^{\varphi})$ by Davydov's inequality, we need the restriction $v\gg \exp(-CL_n^{-\varphi}K^{\varphi})$ in general settings.

{\rm (i)} If $\{\bU_{t_k}\}$ is an independent sequence, we can select $L_n=1/2$ and $\varphi=\infty$ in Assumption {\rm\ref{as:alphmix}}. Then we have $\delta=0$ in this setting, and the upper bound stated in Proposition {\rm\ref{laa:1}} can be refined as $\{1+n_*(K+L_n)^{-1}\rho^{-1}v^2\}^{-\rho/2}+v^{-1}\exp(-Cn_{*}K^{-1}\rho^{-1}v)$ which holds for any $v>0$ and $\rho\geq 1$.

{\rm (ii)} If $\{\bU_{t_k}\}$ is an $L_n$-dependent sequence, we can select $\varphi=\infty$ in Assumption {\rm\ref{as:alphmix}}. If $K>L_n$, we have $\delta=0$ in this setting, and the upper bound stated in Proposition {\rm\ref{laa:1}} can be refined as $\{1+n_*(K+L_n)^{-1}\rho^{-1}v^2\}^{-\rho/2}+v^{-1}\exp(-Cn_{*}K^{-1}\rho^{-1}v)$ which holds for any $v>0$ and $\rho\geq 1$.

\end{rek}

\begin{proposition}\label{la:2}
 Under Assumptions {\rm \ref{as:space}}, {\rm\ref{as:drifdiff}} and {\rm\ref{as:diff}}, we have that
\begin{align*}
\max_{i,j\in[p]}\mathbb{P}\{|{\rm II}(i,j)|>v\}\lesssim\exp(-Cn_*K^{-1}v^2)+\exp\{-C(n_*K^{-1})^{\gamma/(\gamma+4)}\}
\end{align*}
for any $v=o\{(n_*^{-1}K)^{2/(\gamma+4)}\}$, where $\gamma$ is specified in Assumption {\rm\ref{as:diff}}. Furthermore, it holds that $\max_{ i,j\in[p]}\mathbb{E}\{|{\rm II}(i,j)|^{m}\}\lesssim n_*^{-m}K^{m}$ for any fixed positive integer $m$.
\end{proposition}

\begin{proposition}\label{la:3}
 Under Assumptions {\rm\ref{as:alphmix}}--{\rm\ref{as:diff}}, if $K(\log n_*)^{1+2/\gamma}=o(n_*)$ and $K\gtrsim L_n$, we have that
\begin{align*}
\max_{i,j\in[p]}\mathbb{P}\{|{\rm{III}}(i,j)|\geq v\}\lesssim&~\exp(-Cn_*L_n^{-1}v^2)+\exp\{-C(n_*L_n^{-1}v)^{\varphi/(2\varphi+1)}\}\\
&+\exp\{-C(n_*K^{-1})^{\gamma/(\gamma+2)}\}
\end{align*}
for any $v>0$, where $\gamma$ is specified in Assumption {\rm\ref{as:diff}}.
Furthermore,  it holds that $\max_{ i,j\in[p]}\mathbb{E}\{|{\rm III}(i,j)|^{m}\}\lesssim n_*^{-m/2}K^{m/2}$ for any fixed positive integer $m$.
The same results also hold for ${\rm IV}(i,j)$.
\end{proposition}

\begin{rek}
If $\varphi=\infty$, the upper bound stated in Proposition {\rm\ref{la:3}} can be refined as $\exp(-Cn_*L_n^{-1}v^2)+\exp\{-C(n_*L_n^{-1}v)^{1/2}\}+\exp\{-C(n_*K^{-1})^{\gamma/(\gamma+2)}\}
$.
\end{rek}

Write $\aleph=(n_*^{-1}K\log p)^{1/2}$. We first consider the case with $\varphi<\infty$. Notice that $(1+cx^{-1})^{-x}\geq e^{-c}$ for any $x>0$ and $c>0$. By Propositions \ref{laa:1}--\ref{la:3}, if $K(\log n_*)^{1+2/\gamma}=o(n_*)$ and $K\gtrsim L_n$, we have
\begin{align}\label{eq:tailb1}
\notag\max_{i,j\in[p]}\mathbb{P}\big(|\hat{\sigma}_{u,i,j}-\sigma_{u,i,j}|>v\big)
\lesssim& ~\{1+n_*(K+L_n)^{-1}\rho^{-1}v^2\}^{-\rho/2}+v^{-1}\exp\{-C(n_{*}L_n^{-1}\rho^{-1}v)^{\varphi/(\varphi+1)}\}\\
&+v^{-1}\exp(-Cn_{*}K^{-1}\rho^{-1}v)+\exp\{-C(n_*K^{-1})^{\gamma/(\gamma+4)}\}\\
&+\exp\{-C(n_*L_n^{-1}v)^{\varphi/(2\varphi+1)}\}\notag
\end{align}
for any $v\gg \exp(-CL_n^{-\varphi}K^{\varphi})$, $v=o\{(n_*^{-1}K)^{2/(\gamma+4)}\}$ and $\rho\geq 1$. Since $K^{-\varphi}L_n^{\varphi}\log\{n_*(K\log p)^{-1}\}=o(1)$ and $\log p=o\{(n_*K^{-1})^{\gamma/(\gamma+4)}\}$, then  $\exp(-CL_n^{-\varphi}K^{\varphi})=o(\aleph)$ and $\aleph=o\{(n_*^{-1}K)^{2/(\gamma+4)}\}$. Given a sufficiently large constant $\alpha>0$, it holds that
\begin{align*}
\mathbb{E}\big(|\widehat{\bSigma}_{u}-\bSigma_{u}|_\infty\big)
\leq&~\mathbb{E}\bigg\{\max_{i,j\in [p]}|\hat{\sigma}_{u,i,j}-\sigma_{u,i,j}|{I}\big(|\hat{\sigma}_{u,i,j}-\sigma_{u,i,j}|\leq \alpha\aleph\big)\bigg\}
\\&+\mathbb{E}\bigg\{\max_{i,j\in[p]}|\hat{\sigma}_{u,i,j}-\sigma_{u,i,j}|{I}\big(|\hat{\sigma}_{u,i,j}-\sigma_{u,i,j}|> \alpha\aleph\big)\bigg\}
\\=&:A_1+A_2\,.
\end{align*}
It is easy to see that $A_1\leq \alpha\aleph$. By Cauchy-Schwarz inequality, we have
\begin{align*}
A_2\leq&~\sum_{i,j=1}^p\mathbb{E}\big\{|\hat{\sigma}_{u,i,j}-\sigma_{u,i,j}|{I}\big(|\hat{\sigma}_{u,i,j}-\sigma_{u,i,j}|> \alpha\aleph\big)\big\}\\
\leq&~p^2\max_{i,j\in[p]}\big\{\mathbb{E}\big(|\hat{\sigma}_{u,i,j}-\sigma_{u,i,j}|^2\big)\big\}^{1/2}\cdot\max_{i,j\in[p]}\big\{\mathbb{P}\big(|\hat{\sigma}_{u,i,j}-\sigma_{u,i,j}|> \alpha\aleph\big)\big\}^{1/2}\,.
\end{align*}
Let $\rho\asymp\log p\geq 1$. Since $K\log p=o(n_*)$, then it follows from (\ref{eq:tailb1}) that
\begin{align*}
\max_{i,j\in[p]}\mathbb{P}\big(|\hat{\sigma}_{u,i,j}-\sigma_{u,i,j}|> \alpha\aleph\big)
\lesssim&~p^{-2w}+\exp\{-C(n_*KL_n^{-2}\log p)^{\varphi/(4\varphi+2)}\}\\
&+\exp[-C\{n_*K^{-1}(\log p)^{-1}\}^{1/2}]+\exp\{-C(n_*K^{-1})^{\gamma/(\gamma+4)}\}
\\&+
\exp[-C\{n_*KL_n^{-2}(\log p)^{-1}\}^{\varphi/(2\varphi+2)}]
\end{align*}
with some sufficiently large $w>0$, where $w\rightarrow\infty$ as $\alpha\rightarrow\infty$. Due to
$
\max_{i,j\in[p]}\mathbb{E}(|\hat{\sigma}_{u,i,j}-\sigma_{u,i,j}|^2)\lesssim 1$, if $\log p=o[\min\{(n_*L_n^{-2}K)^{\varphi/(3\varphi+2)},(n_*K^{-1})^{\chi}\}]$ with $\chi=\min\{\gamma/(\gamma+4),1/3\}$, then
\begin{align*}
A_2\lesssim&~p^{2-w}+p^2\exp\{-C(n_*KL_n^{-2}\log p)^{\varphi/(4\varphi+2)}\}+p^2\exp[-C\{n_*K^{-1}(\log p)^{-1}\}^{1/2}]\\
&+p^2\exp\{-C(n_*K^{-1})^{\gamma/(\gamma+4)}\}+p^2\exp[-C\{n_*KL_n^{-2}(\log p)^{-1}\}^{\varphi/(2\varphi+2)}]\\
\lesssim&~p^{2-w}+\exp\{-C(n_*KL_n^{-2}\log p)^{\varphi/(4\varphi+2)}\}+\exp[-C\{n_*K^{-1}(\log p)^{-1}\}^{1/2}]\\
&+\exp\{-C(n_*K^{-1})^{\gamma/(\gamma+4)}\}+\exp[-C\{n_*KL_n^{-2}(\log p)^{-1}\}^{\varphi/(2\varphi+2)}]\\
=&~o\{(n_*^{-1}K\log p)^4\}\,.
\end{align*}
Hence,
$
\sup_{\mathcal{P}_1}\mathbb{E}(|\widehat{\bSigma}_{u}-\bSigma_{u}|_\infty)\lesssim(n_*^{-1}K\log p)^{1/2}$ provided that $K(\log n_*)^{1+2/\gamma}=o(n_*)$,  $K\gtrsim L_n$, $K^{-\varphi}L_n^{\varphi}\log\{n_*(K\log p)^{-1}\}=o(1)$ and $\log p=o[\min\{(n_*L_n^{-2}K)^{\varphi/(3\varphi+2)},(n_*K^{-1})^{\chi}\}]$.

Now we consider the case with $\varphi=\infty$. As we discussed in Remark \ref{rek:1}(i), if $\{\bU_{t_k}\}$ is an independent sequence, we can select $L_n=1/2$. Due to $K\geq 1$, we have $K>L_n$ in this case. Without loss of generality, we can always assume $K>L_n$ when $\varphi=\infty$. Based on Remark \ref{rek:1}, it holds that $\{1+n_*(K+L_n)^{-1}\rho^{-1}v^2\}^{-\rho/2}+v^{-1}\exp(-Cn_{*}K^{-1}\rho^{-1}v)$ for any $v>0$ and $\rho\geq 1$ under either of the scenarios: (i) $\{\bU_{t_k}\}$ is an independent sequence, and (ii) $\{\bU_{t_k}\}$ is an $L_n$-dependent sequence. Repeating the arguments for $\varphi<\infty$, we have $
\sup_{\mathcal{P}_1}\mathbb{E}(|\widehat{\bSigma}_{u}-\bSigma_{u}|_\infty)\lesssim(n_*^{-1}K\log p)^{1/2}$ provided that $K(\log n_*)^{1+2/\gamma}=o(n_*)$,  $K>L_n$ and $\log p=o\{(n_*K^{-1})^{\chi}\}$. We complete the proof of Theorem \ref{pn:1}. $\hfill\Box$

\subsection{Proof of Proposition \ref{laa:1}}\label{sec:pfpn1}

To prove Proposition \ref{laa:1}, we need the following lemma whose proof is given in Section \ref{sec:pfberninq}.

\begin{lemma}\label{bern.inq}
	Let $\{z_t\}_{t=1}^{\tilde{n}} $ be an $\alpha$-mixing sequence of real-valued and centered random variables with $\alpha$-mixing coefficients $\{\alpha(k)\}_{k\geq1}$. Assume there exist some universal constants $a_1>1$, $a_2, b_1,b_2>0$, $r\in (0,2]$ and $\varphi>0$ such that {\rm(i)} $
		\max_{t\in[\tilde{n}]}\mathbb{P}(|z_t|>u)\leq b_1\exp(-b_2u^{r})$
		for any $u>0$, {\rm(ii)} $
		\alpha(k)\leq a_1\exp(-a_2\tilde{L}_{\tilde{n}}^{-\varphi}|k-m|_+^{\varphi})$
		for any integer $k\geq 1$, where $\tilde{L}_{\tilde{n}}>0$ and $m=o(\tilde{n})\geq 0$ may diverge with $\tilde{n}$.	Let $s_{\tilde{n}}^2=\sum_{t_1,t_2=1}^{\tilde{n}}|\mathrm{Cov}(z_{t_1},z_{t_2})|$ and $r_*=r\varphi/(r+\varphi)$. It holds that
\begin{align*}
\mathbb{P}\bigg(\sup_{k\in[\tilde{n}]}\bigg|\sum_{t=1}^{k}z_t\bigg|\geq\tilde{n}x\bigg)\lesssim&~(1+\tilde{n}^2x^2\rho^{-1}s_{\tilde{n}}^{-2})^{-\rho/2}+x^{-1}\exp(-C\tilde{n}^{r_*}x^{r_*}\rho^{-r_*}\tilde{L}_{\tilde{n}}^{-r_*})\\
&+x^{-1}\exp(-Cm^{-r}\tilde{n}^{r}x^{r}\rho^{-r})
\end{align*}
for any $x>0$ and $\rho\geq 1$, where we adopt the convention $\exp(-C0^{-r}x^r)=0$ for any $x>0$.
\end{lemma}
	
\begin{rek}
If $\varphi=\infty$, we have $r_*=r$. Then the upper bound in Lemma {\rm\ref{bern.inq}} can be refined as $
	(1+\tilde{n}^2x^2\rho^{-1}s_{\tilde{n}}^{-2})^{-\rho/2}+x^{-1}\exp\{-C\tilde{n}^{r}x^{r}\rho^{-r}(m+\tilde{L}_{\tilde{n}})^{-r}\}$.
\end{rek}

Now we begin to prove Proposition \ref{laa:1}. Recall that
\begin{align}\label{eq:I(i,j)}
{\rm I}(i,j)
=&~\underbrace{\frac{1}{n_{i,j}}\sum_{\ell=1}^{n_{i,j}}\bigg(\frac{1}{2}+\frac{1}{2}\sum_{t_{i,j,k}\in S_{i,j,\ell}}\frac{1}{N_{i,j,k}}\bigg)(U_{i,t_{i,j,\ell}}U_{j,t_{i,j,\ell}}-\sigma_{u,i,j})}_{\textrm{I}_1(i,j)}\notag\\
&-\underbrace{\frac{1}{2n_{i,j}}\sum_{k=1}^{n_{i,j}}\frac{U_{i,t_{i,j,k}}}{N_{i,j,k}}\sum_{t_{i,j,\ell}\in S_{i,j,k}}U_{j,t_{i,j,\ell}}}_{\textrm{I}_2(i,j)}-\underbrace{\frac{1}{2n_{i,j}}\sum_{k=1}^{n_{i,j}}\frac{U_{j,t_{i,j,k}}}{N_{i,j,k}}\sum_{t_{i,j,\ell}\in S_{i,j,k}}U_{i,t_{i,j,\ell}}}_{\textrm{I}_3(i,j)}\,.
\end{align}
In the sequel, we will bound the tail probabilities of $\textrm{I}_1(i,j)$, $\textrm{I}_2(i,j)$ and $\textrm{I}_3(i,j)$, respectively.

For each $\ell\in[n_{i,j}]$, let
$
\zeta_{i,j,\ell}=2^{-1}(1+\sum_{t_{i,j,k}\in S_{i,j,\ell}}N_{i,j,k}^{-1})(U_{i,t_{i,j,\ell}}U_{j,t_{i,j,\ell}}-\sigma_{u,i,j})$.
Then we have
$
\text{I}_{1}(i,j)=n_{i,j}^{-1}\sum_{\ell=1}^{n_{i,j}}\zeta_{i,j,\ell}$.
Recall that $N_{i,j,k}=|S_{i,j,k}|$ with $S_{i,j,k}=\{t_{i,j,\ell}:K\leq|\ell-k|\leq K+\Delta_K\}$. Since $K<K+\Delta_K=o(n_*)$ and $n_{i,j}\asymp n_*\to \infty$ as $n\to \infty$, we then have $2(K+\Delta_K)<n_{i,j}$ for sufficiently large $n$. Thus for sufficiently large $n$, it holds that
\begin{equation}\label{eq:Nijk}
N_{i,j,k}=\left\{
\begin{aligned}
\Delta_K+1\,, ~~~~~~~~~~~~&\textrm{if}~1\leq k \leq K\,,  \\
\Delta_K-K+k+1\,, ~~~~~~&\textrm{if}~K+1\leq k\leq K+\Delta_K\,,\\
2\Delta_K+2\,, ~~~~~~~~~~~& \textrm{if}~K+\Delta_K+1\leq k\leq n_{i,j}-K-\Delta_K\,,\\
n_{i,j}+\Delta_K-K-k+2\,, ~~~& \textrm{if}~n_{i,j}-K-\Delta_K+1\leq k\leq n_{i,j}-K\,,
\\
\Delta_K+1\,, ~~~~~~~~~~~~& \textrm{if}~n_{i,j}-K+1\leq k\leq n_{i,j}\,,
\end{aligned}
\right.
\end{equation}
which implies that $\min_{k \in[n_{i,j}]}N_{i,j,k}=\Delta_K+1$ and $\max_{k\in[n_{i,j}]}N_{i,j,k}=2\Delta_K+2$. Therefore, we have that $C^{-1}<\sum_{t_{i,j,k}\in S_{i,j,\ell}}N_{i,j,k}^{-1}<C$ holds uniformly over $\ell\in[n_{i,j}]$ and $i,j\in[p]$.
By Lemma 2 of \cite{Changetal_2013_AOS}, Assumption \ref{as:moment} yields
$
\max_{i,j\in[p],\ell\in[n_{i,j}]}\mathbb{P}(|\zeta_{i,j,\ell}|>v)\leq C\exp(-Cv)$
for any $v>0$, which implies $\max_{i,j\in[p],\ell\in[n_{i,j}]}\mathbb{E}(|\zeta_{i,j,\ell}|^4)\leq C$.
It follows from Davydov's inequality that
$
\sum_{\ell_1,\ell_2=1}^{n_{i,j}}|\mathrm{Cov}(\zeta_{i,j,\ell_1},\zeta_{i,j,\ell_2})|\lesssim\sum_{\ell_2=1}^{n_{i,j}}1+\sum_{\ell_2=1}^{n_{i,j}}\sum_{\ell_1=\ell_2+1}^{n_{i,j}}\exp(-CL_n^{-\varphi}|\ell_1-\ell_2|^{\varphi})\lesssim n_*L_n$.
By Lemma \ref{bern.inq} with $m=0$ and $\tilde{L}_{\tilde{n}}=L_n$, we have
\begin{align}\label{bound I1}
\max_{i,j\in[p]}\mathbb{P}\{|\mathrm{I}_1(i,j)|\geq v\}\lesssim (1+n_*L_n^{-1}\rho^{-1}v^2)^{-\rho/2}+v^{-1}\exp\{-C(n_{*}L_n^{-1}\rho^{-1}v)^{\varphi/(\varphi+1)}\}
\end{align}
for any $v>0$ and $\rho\geq 1$.

Define
$
\eta_{i,j,k}=(2N_{i,j,k})^{-1}U_{i,t_{i,j,k}}\sum_{t_{i,j,\ell}\in S_{i,j,k}}U_{j,t_{i,j,\ell}}$.
 Then
$
 \text{I}_2(i,j)= n_{i,j}^{-1}\sum_{k=1}^{n_{i,j}}\eta_{i,j,k}$. Since
 \begin{align*}
 \max_{i,j\in[p],k\in[n_{i,j}]}N_{i,j,k}\leq 2\Delta_K+2
 \end{align*} for sufficiently large $n$, we know $\max_{i,j\in[p],k\in[n_{i,j}]}N_{i,j,k}$ is uniformly bounded away from infinity due to the fact $\Delta_K$ is a fixed constant.
It follows from Assumption \ref{as:moment} that
\begin{align*}
\max_{i,j\in[p],k\in[n_{i,j}]}\mathbb{P}(|\eta_{i,j,k}|\geq x)\leq&~\max_{i,j\in[p],k\in[n_{i,j}]}\mathbb{P}\bigg(\bigg|\frac{1}{N_{i,j,k}}\sum_{t_{i,j,\ell}\in S_{i,j,k}}U_{j,t_{i,j,\ell}}\bigg|\geq \sqrt{x}\bigg) \\
&+\max_{i,j\in[p],k\in[n_{i,j}]}\mathbb{P}\big(|U_{i,t_{i,j,k}}|\geq 2\sqrt{x}\big)\\
\leq&~\max_{i,j\in[p],k\in[n_{i,j}]}\sum_{t_{i,j,\ell}\in S_{i,j,k}}\mathbb{P}\big(|U_{j,t_{i,j,\ell}}|\geq \sqrt{x}\big) \\
&+\max_{i,j\in[p],k\in[n_{i,j}]}\mathbb{P}\big(|U_{i,t_{i,j,k}}|\geq 2\sqrt{x}\big)\\
\lesssim&~
\exp(-Cx)
\end{align*}
for any $x>0$. By Davydov's inequality, we have
$
\sum_{k_1,k_2=1}^{n_{i,j}}|\mathrm{Cov}(\eta_{i,j,k_1},\eta_{i,j,k_2})|\lesssim
n_{*}(K+L_n)$.
Write $\mathring{\eta}_{i,j,k}=\eta_{i,j,k}-\mu_{i,j,k}$ with $\mu_{i,j,k}=\mathbb{E}(\eta_{i,j,k})$. By Lemma \ref{bern.inq} with $m=2(K+\Delta_K)$ and $\tilde{L}_{\tilde{n}}=L_n$,
$
\max_{i,j\in[p]}\mathbb{P}(|n_{i,j}^{-1}\sum_{k=1}^{n_{i,j}}\mathring{\eta}_{i,j,k}|\geq x)\lesssim \{1+n_*(K+L_n)^{-1}\rho^{-1}x^2\}^{-\rho/2}+x^{-1}\exp\{-C(n_{*}L_n^{-1}\rho^{-1}x)^{\varphi/(\varphi+1)}\}+x^{-1}\exp(-Cn_{*}K^{-1}\rho^{-1}x)$
for any $x>0$ and $\rho\geq 1$. Let $S_{i,j,k,1}=\{t_{i,j,\ell}\in S_{i,j,k}:\ell<k\}$ and $S_{i,j,k,2}=\{t_{i,j,\ell}\in S_{i,j,k}:\ell>k\}$. Then
$
\mu_{i,j,k}=
(2N_{i,j,k})^{-1}\{\mathbb{E}(U_{i,t_{i,j,k}}\sum_{t_{i,j,\ell}\in S_{i,j,k,1}}U_{j,t_{i,j,\ell}})+\mathbb{E}(U_{i,t_{i,j,k}}\sum_{t_{i,j,\ell}\in S_{i,j,k,2}}U_{j,t_{i,j,\ell}})\}$.
Applying Davydov's inequality and Jensen's inequality, it holds that
$
|\mathbb{E}(U_{i,t_{i,j,k}}\sum_{t_{i,j,\ell}\in S_{i,j,k,1}}U_{j,t_{i,j,\ell_1}})|\lesssim\{\mathbb{E}(|U_{i,t_{i,j,k}}|^{4})\}^{1/4}\{\mathbb{E}(|\sum_{t_{i,j,\ell}\in S_{i,j,k,1}}U_{j,t_{i,j,\ell}}|^{4})\}^{1/4}\alpha_n^{1/2}(K)\lesssim\exp(-CL_n^{-\varphi}K^{\varphi})$.
Analogously, we also have $|\mathbb{E}\big(U_{i,t_{i,j,k}}\sum_{t_{i,j,\ell}\in S_{i,j,k,2}}U_{j,t_{i,j,\ell}}\big)|\lesssim \exp(-CL_n^{-\varphi}K^{\varphi})$. Thus, it holds that $\max_{i,j\in[p],k\in[n_{i,j}]}|\mu_{i,j,k}|\lesssim \exp(-CL_n^{-\varphi}K^{\varphi})$, which implies $
	\max_{i,j\in[p]}|n_{i,j}^{-1}\sum_{k=1}^{n_{i,j}}\mu_{i,j,k}|\lesssim\exp(-CL_n^{-\varphi}K^{\varphi})$. For any $v\gg \exp(-CL_n^{-\varphi}K^{\varphi})$ and $\rho\geq 1$, we have
\begin{align*}
\max_{i,j\in[p]}\mathbb{P}\{|\mathrm{I}_2(i,j)|\geq v\}\leq&~\max_{i,j\in[p]}\mathbb{P}\bigg(\bigg|\frac{1}{n_{i,j}}\sum_{k=1}^{n_{i,j}}\mathring{\eta}_{i,j,k}\bigg|\geq v-\bigg|\frac{1}{n_{i,j}}\sum_{k=1}^{n_{i,j}}\mu_{i,j,k}\bigg|\bigg)\\
\leq&~\max_{i,j\in[p]}\mathbb{P}\biggl(\bigg|\frac{1}{n_{i,j}}\sum_{k=1}^{n_{i,j}}\mathring{\eta}_{i,j,k}\bigg|\geq v/2\bigg)\\
\lesssim &~\{1+n_*(K+L_n)^{-1}\rho^{-1}v^2\}^{-\rho/2}+v^{-1}\exp\{-C(n_{*}L_n^{-1}\rho^{-1}v)^{\varphi/(\varphi+1)}\}\\
&+v^{-1}\exp(-Cn_{*}K^{-1}\rho^{-1}v)\,.
\end{align*}
Analogously,
$
\max_{i,j\in[p]}\mathbb{P}\{|\mathrm{I}_3(i,j)|\geq v\}\lesssim \{1+n_*(K+L_n)^{-1}\rho^{-1}v^2\}^{-\rho/2}+v^{-1}\exp\{-C(n_{*}L_n^{-1}\rho^{-1}v)^{\varphi/(\varphi+1)}\}\\
+v^{-1}\exp(-Cn_{*}K^{-1}\rho^{-1}v)$
for any $v\gg \exp(-CL_n^{-\varphi}K^{\varphi})$ and $\rho\geq 1$. Note that ${\rm I}(i,j)={\rm I}_1(i,j)+{\rm I}_2(i,j)+{\rm I}_3(i,j)$. Together with (\ref{bound I1}),  it holds that
$
\max_{i,j\in[p]}\mathbb{P}\{|\mathrm{I}(i,j)|\geq v\}
\lesssim \{1+n_*(K+L_n)^{-1}\rho^{-1}v^2\}^{-\rho/2}+v^{-1}\exp\{-C(n_{*}L_n^{-1}\rho^{-1}v)^{\varphi/(\varphi+1)}\}
+v^{-1}\exp(-Cn_{*}K^{-1}\rho^{-1}v)$
for any $v\gg \exp(-CL_n^{-\varphi}K^{\varphi})$ and $\rho\geq 1$.

By \eqref{eq:I(i,j)},
$|{\rm I}(i,j)|^{m}\lesssim |{\rm I}_1(i,j)|^{m}+|{\rm I}_2(i,j)|^{m}+|{\rm I}_3(i,j)|^{m}$ for any fixed positive integer $m$. By Assumption \ref{as:moment},  $\mathbb{E}\{|{\rm I}_1(i,j)|^{m}\}\lesssim
n_{i,j}^{-1}\sum_{\ell=1}^{n_{i,j}}\mathbb{E}(|U_{i,t_{i,j,\ell}}U_{j,t_{i,j,\ell}}-\sigma_{u,i,j}|^{m})\lesssim  1$ and  $\mathbb{E}\{|{\rm I}_2(i,j)|^{m}\}\lesssim
n_{i,j}^{-1}\sum_{k=1}^{n_{i,j}}N_{i,j,k}^{-1}\sum_{t_{i,j,\ell}\in S_{i,j,k}}\mathbb{E}(|U_{i,t_{i,j,k}}U_{j,t_{i,j,\ell}}|^{m})\lesssim 1
$. Analogously, we also have $\mathbb{E}\{|{\rm I}_3(i,j)|^{m}\}\lesssim 1$. Thus,  $\max_{ i,j\in[p]}\mathbb{E}\{|{\rm I}(i,j)|^{m}\}\lesssim 1$.  We complete the proof of  Proposition \ref{laa:1}.
$\hfill\Box$

\subsection{Proof of Proposition \ref{la:2}}\label{sec:pfpn2}

Notice that ${\rm d}X_{i,t}=\mu_{i,t}\,{\rm d}t+\sigma_{i,t}\,{\rm d}B_{i,t}$. Then
\begin{align}\label{eq:IIij}
{\rm II}(i,j)=&~\underbrace{\frac{1}{2n_{i,j}}\sum_{k=1}^{n_{i,j}}\frac{1}{N_{i,j,k}}\sum_{t_{i,j,\ell}\in S_{i,j,k}}\bigg(\int_{t_{i,j,k}}^{t_{i,j,\ell}}\mu_{i,s}\,{\rm d}s\bigg)\bigg(\int_{t_{i,j,k}}^{t_{i,j,\ell}}\mu_{j,s}\,{\rm d}s\bigg)}_{{\rm II}_1(i,j)}\notag\\
&+\underbrace{\frac{1}{2n_{i,j}}\sum_{k=1}^{n_{i,j}}\frac{1}{N_{i,j,k}}\sum_{t_{i,j,\ell} \in S_{i,j,k}}\bigg(\int_{t_{i,j,k}}^{t_{i,j,\ell}}\sigma_{i,s}\,{\rm d}B_{i,s}\bigg)\bigg(\int_{t_{i,j,k}}^{t_{i,j,\ell}}\sigma_{j,s}\,{\rm d}B_{j,s}\bigg)}_{{\rm II}_2(i,j)}\notag\\
&+\underbrace{\frac{1}{2n_{i,j}}\sum_{k=1}^{n_{i,j}}\frac{1}{N_{i,j,k}}\sum_{t_{i,j,\ell} \in S_{i,j,k}}\bigg(\int_{t_{i,j,k}}^{t_{i,j,\ell}}\mu_{i,s}\,{\rm d}s\bigg)\bigg(\int_{t_{i,j,k}}^{t_{i,j,\ell}}\sigma_{j,s}\,{\rm d}B_{j,s}\bigg)}_{{\rm II}_3(i,j)}\\
&+\underbrace{\frac{1}{2n_{i,j}}\sum_{k=1}^{n_{i,j}}\frac{1}{N_{i,j,k}}\sum_{t_{i,j,\ell} \in S_{i,j,k}}\bigg(\int_{t_{i,j,k}}^{t_{i,j,\ell}}\sigma_{i,s}\,{\rm d}B_{i,s}\bigg)\bigg(\int_{t_{i,j,k}}^{t_{i,j,\ell}}\mu_{j,s}\,{\rm d}s\bigg)}_{{\rm II}_4(i,j)}\,.\notag
\end{align}
Recall $\xi=\max_{i,j\in[p]}\max_{ k\in[n_{i,j}]}\max_{t_{i,j,\ell}\in S_{i,j,k}}|t_{i,j,\ell}-t_{i,j,k}|$. In the sequel, we will bound the tail probabilities of
$\max_{i,j\in[p]}|\textrm{II}_1(i,j)|$,
 $\max_{i,j\in[p]}|\textrm{II}_2(i,j)|$, $\max_{i,j\in[p]}|\textrm{II}_3(i,j)|$ and $\max_{i,j\in[p]}|\textrm{II}_4(i,j)|$, respectively.

For any $k=[n_{i,j}]$, define $
\zeta_{i,j,k}^*=N_{i,j,k}^{-1}\sum_{t_{i,j,\ell}\in S_{i,j,k}}(\int_{t_{i,j,k}}^{t_{i,j,\ell}}\mu_{i,s}\,{\rm d}s)(\int_{t_{i,j,k}}^{t_{i,j,\ell}}\mu_{j,s}\,{\rm d}s)$. Then we have $
{\rm II}_1(i,j)=(2n_{i,j})^{-1}\sum_{k=1}^{n_{i,j}}\zeta_{i,j,k}^*$. We will first bound $\mathbb{E}\{\exp(\theta\zeta_{i,j,k}^*)\}$ for any $|\theta|\in(0,C_5\xi^{-2}]$, where $C_5$ is specified in Assumption \ref{as:drifdiff}. By Jensen's inequality and Cauchy-Schwarz inequality,
\begin{align}\label{eq:b1}
\mathbb{E}\{\exp(\theta\zeta_{i,j,k}^*)\}\leq&~\frac{1}{N_{i,j,k}}\sum_{t_{i,j,\ell}\in S_{i,j,k}}\mathbb{E}\bigg[\exp\bigg\{\theta\bigg(\int_{t_{i,j,k}}^{t_{i,j,\ell}}\mu_{i,s}\,{\rm d}s\bigg)\bigg(\int_{t_{i,j,k}}^{t_{i,j,\ell}}\mu_{j,s}\,{\rm d}s\bigg)\bigg\}\bigg]\notag\\
\leq&~\frac{1}{N_{i,j,k}}\sum_{t_{i,j,\ell}\in S_{i,j,k}}\frac{1}{|t_{i,j,\ell}-t_{i,j,k}|^2}\notag\\
&~~~~~~~~~~~\times\int_{t_{i,j,k}\wedge t_{i,j,\ell}}^{t_{i,j,k}\vee t_{i,j,\ell}}\int_{t_{i,j,k}\wedge t_{i,j,\ell}}^{t_{i,j,k}\vee t_{i,j,\ell}}\mathbb{E}\{\exp(\theta|t_{i,j,\ell}-t_{i,j,k}|^2\mu_{i,s_1}\mu_{j,s_2})\}\,{\rm d}s_1{\rm d}s_2\\
\leq&~\frac{1}{N_{i,j,k}}\sum_{t_{i,j,\ell} \in S_{i,j,k}}\frac{1}{|t_{i,j,\ell}-t_{i,j,k}|^2}\notag\\
&~~~~~~~~~~~\times\int_{t_{i,j,k}\wedge t_{i,j,\ell}}^{t_{i,j,k}\vee t_{i,j,\ell}}\big[\mathbb{E}\{\exp(|\theta||t_{i,j,\ell}-t_{i,j,k}|^2\mu_{i,s_1}^2)\}\big]^{1/2}\,{\rm d}s_1\notag\\
&~~~~~~~~~~~\times\int_{t_{i,j,k}\wedge t_{i,j,\ell}}^{t_{i,j,k}\vee t_{i,j,\ell}}\big[\mathbb{E}\{\exp(|\theta||t_{i,j,\ell}-t_{i,j,k}|^2\mu_{j,s_2}^2)\}\big]^{1/2}\,{\rm d}s_2\notag\\
\leq&~\max_{i\in[p]}\sup_{0\leq s\leq T}\mathbb{E}\{\exp(|\theta|\xi^2\mu_{i,s}^2)\}\notag\,.
\end{align}
Recall that $\xi\asymp Kn_*^{-1}=o(1)$. By Assumption \ref{as:drifdiff},
$
\max_{i\in[p]}\sup_{0\leq s\leq T}\mathbb{E}\{\exp(|\theta|\xi^2\mu_{i,s}^2)\}\leq \exp(C_5C_7)\\\cdot\exp(C_6\xi^4\theta^2)\lesssim\exp(C\xi\theta^2)
$
for any $|\theta|\in(0,C_5\xi^{-2}]$. Therefore, by (\ref{eq:b1}), $
\max_{i,j\in[p],k\in[n_{i,j}]}\mathbb{E}\{\exp(\theta\zeta_{i,j,k}^*)\}\lesssim\exp(C\xi\theta^2)$ for any $|\theta|\in(0,C_5\xi^{-2}]$. By Lemma 2 of \cite{Fan2012}, it holds that
\begin{equation}\label{eq:c1}
\max_{i,j\in[p]}\mathbb{P}\{|{\rm II}_1(i,j)|>v\}\lesssim\exp(-Cv^2\xi^{-1})
\end{equation}
for any $v=o(\xi^{-1})$.

For any $k\in[n_{i,j}]$, define $
\zeta_{i,j,k}^{**}=N_{i,j,k}^{-1}\sum_{t_{i,j,\ell} \in S_{i,j,k}}(\int_{t_{i,j,k}}^{t_{i,j,\ell}}\sigma_{i,s}\,{\rm d}B_{i,s})(\int_{t_{i,j,k}}^{t_{i,j,\ell}}\sigma_{j,s}\,{\rm d}B_{j,s})$. Then we have
${\rm II}_2(i,j)=(2n_{i,j})^{-1}\sum_{k=1}^{n_{i,j}}\zeta_{i,j,k}^{**}$.
 For any constant $d\in(0,\xi^{-1/2}]$, define a stopping time $
 \Gamma_{i,d}=T\wedge\inf\{t>0:\sup_{0\leq s\leq t}\sigma_{i,s}>d\}$.
  For any $|\theta|\in(0,d^{-2}\xi^{-1}/4]$, by Jensen's inequality and Cauchy-Schwarz inequality, it holds that
\begin{align}\label{eq:b2}
&\mathbb{E}\{\exp(\theta\zeta_{i,j,k}^{**}){I}(\Gamma_{i,d}=\Gamma_{j,d}=T)\}\notag\\
&~~~~~~~\leq\frac{1}{N_{i,j,k}}\sum_{t_{i,j,\ell}\in S_{i,j,k}}\mathbb{E}\bigg[\exp\bigg\{\theta\bigg(\int_{t_{i,j,k}}^{t_{i,j,\ell}}\sigma_{i,s}\,{\rm d}B_{i,s}\bigg)\notag\\
&~~~~~~~~~~~~~~~~~~~~~~~~~~~~~~~~~~~~~~~~~~~\times\bigg(\int_{t_{i,j,k}}^{t_{i,j,\ell}}\sigma_{j,s}\,{\rm d}B_{j,s}\bigg)\bigg\}{I}(\Gamma_{i,d}=\Gamma_{j,d}=T)\bigg]\\
&~~~~~~~\leq\frac{1}{N_{i,j,k}}\sum_{t_{i,j,\ell} \in S_{i,j,k}}\bigg(\mathbb{E}\bigg[\exp\bigg\{|\theta|\bigg(\int_{t_{i,j,k}}^{t_{i,j,\ell}}\sigma_{i,s}\,{\rm d}B_{i,s}\bigg)^2\bigg\}{I}(\Gamma_{i,d}=T)\bigg]\bigg)^{1/2}\notag\\
&~~~~~~~~~~~~~~~~~~~~~~~~~~~~~~~~~~~~~~\times\bigg(\mathbb{E}\bigg[\exp\bigg\{|\theta|\bigg(\int_{t_{i,j,k}}^{t_{i,j,\ell}}\sigma_{j,s}\,{\rm d}B_{j,s}\bigg)^2\bigg\}{I}(\Gamma_{j,d}=T)\bigg]\bigg)^{1/2}\,.\notag
\end{align}
Restricted on the event $\{\Gamma_{i,d}=T\}$, we have $\sup_{0\leq s\leq T}\sigma_{i,s}\leq d$. For any $|\theta|\in(0,d^{-2}\xi^{-1}/4]$, it holds that
\begin{align}\label{eq:re1}
&\exp\bigg\{|\theta|\bigg(\int_{t_{i,j,k}}^{t_{i,j,\ell}}\sigma_{i,s}\,{\rm d}B_{i,s}\bigg)^2\bigg\}{I}(\Gamma_{i,d}=T)\notag\\
&~~~~~~~~~~~\leq\exp\bigg[|\theta|\bigg\{\bigg(\int_{t_{i,j,\ell} \wedge t_{i,j,k}}^{t_{i,j,\ell} \vee t_{i,j,k}}\sigma_{i,s}\,{\rm d}B_{i,s}\bigg)^2-\int_{t_{i,j,\ell} \wedge t_{i,j,k}}^{t_{i,j,\ell} \vee t_{i,j,k}}\sigma_{i,s}^2\,{\rm d}s\bigg\}\bigg]{I}(\Gamma_{i,d}=T)\\
&~~~~~~~~~~~~~~~~~~\times\exp\bigg(|\theta|\int_{t_{i,j,\ell} \wedge t_{i,j,k}}^{t_{i,j,\ell} \vee t_{i,j,k}}\sigma_{i,s}^2\,{\rm d}s\bigg){I}(\Gamma_{i,d}=T)\notag\\
&~~~~~~~~~~~\leq C\exp\bigg[|\theta|\bigg\{\bigg(\int_{t_{i,j,\ell} \wedge t_{i,j,k}}^{t_{i,j,\ell} \vee t_{i,j,k}}\sigma_{i,s}\,{\rm d}B_{i,s}\bigg)^2-\int_{t_{i,j,\ell} \wedge t_{i,j,k}}^{t_{i,j,\ell} \vee t_{i,j,k}}\sigma_{i,s}^2\,{\rm d}s\bigg\}\bigg]{I}(\Gamma_{i,d}=T)\,.\notag
\end{align}
Recall $d\leq\xi^{-1/2}$. Following the arguments of Equation (A.5) in \cite{Fan2012}, we have that
\begin{align}\label{eq:re2}
&\mathbb{E}\bigg\{\exp\bigg[|\theta|\bigg\{\bigg(\int_{t_{i,j,\ell} \wedge t_{i,j,k}}^{t_{i,j,\ell} \vee t_{i,j,k}}\sigma_{i,s}\,{\rm d}B_{i,s}\bigg)^2\notag-\int_{t_{i,j,\ell} \wedge t_{i,j,k}}^{t_{i,j,\ell} \vee t_{i,j,k}}\sigma_{i,s}^2\,{\rm d}s\bigg\}\bigg]{I}(\Gamma_{i,d}=T)\,\bigg|\,\mathcal{F}^*_{i,t_{i,j,\ell}\wedge t_{i,j,k}}\bigg\}\notag\\
&~~~~~~~~\leq \mathbb{E}[\exp\{|\theta|(B^2_{i,d^2|t_{i,j,\ell}-t_{i,j,k}|}-d^2|t_{i,j,\ell}-t_{i,j,k}|)\}]=\mathbb{E}[\exp\{|\theta| d^2|t_{i,j,\ell}-t_{i,j,k}|(Z^2-1)\}]\notag\\
&~~~~~~~~\leq \exp(Cd^4\xi^2\theta^2)\leq \exp\{C(d\xi^{1/2})^{4-\tau}\theta^2\}
\end{align}
for any $\tau\in[0,4)$, where $Z\sim N(0,1)$ and $\mathcal{F}_{i,t}^*$ is the $\sigma$-field generated by $(\sigma_{i,s},B_{i,s})_{0\leq s\leq t}$. Thus, by (\ref{eq:b2}) and (\ref{eq:re1}), we have
$
\max_{i,j\in[p],k\in[n_{i,j}]}\mathbb{E}\{\exp(\theta\zeta_{i,j,k}^{**}){I}(\Gamma_{i,d}=\Gamma_{j,d}=T)\}\lesssim\exp\{C(d\xi^{1/2})^{4-\tau}\theta^2\}$
for any $|\theta|\in(0,d^{-2}\xi^{-1}/4]$. By Lemma 2 of \cite{Fan2012}, it holds that
$
\max_{i,j\in[p]}\mathbb{P}\{|{\rm II}_2(i,j)|>v,\Gamma_{i,d}=\Gamma_{j,d}=T\}\lesssim\exp\{-Cv^2(d\xi^{1/2})^{\tau-4}\}$
for any $v=o\{(d\xi^{1/2})^{2-\tau}\}$. Note that $\mathbb{P}\{|{\rm II}_2(i,j)|>v\}\leq \mathbb{P}\{|{\rm II}_2(i,j)|>v, \Gamma_{i,d}=\Gamma_{j,d}=T\}+\mathbb{P}(\Gamma_{i,d}\neq T)+\mathbb{P}(\Gamma_{j,d}\neq T)$. Since $\Gamma_{i,d}=T\wedge\inf\{t>0:\sup_{0\leq s\leq t}\sigma_{i,s}>d\}$, by Assumption \ref{as:diff}, we have
$
\max_{i\in[p]}\mathbb{P}(\Gamma_{i,d}\neq T)\leq \max_{i\in[p]}\mathbb{P}(\sup_{0\leq s\leq T}\sigma_{i,s}>d)\lesssim\exp(-Cd^{\gamma})$.
Then
\begin{equation}\label{eq:c2}
\max_{i,j\in[p]}\mathbb{P}\{|{\rm II}_2(i,j)|>v\}\lesssim\exp\{-Cv^2(d\xi^{1/2})^{\tau-4}\}+\exp(-Cd^\gamma)
\end{equation}
for any $v=o\{(d\xi^{1/2})^{2-\tau}\}$ with $d\leq \xi^{-1/2}$.

Due to $xy\leq 2^{-1}(x^2+y^2)$ for any $x,y>0$, we have
\begin{align*}
2|{\rm II}_3(i,j)|\leq&~\frac{1}{2n_{i,j}}\sum_{k=1}^{n_{i,j}}\frac{1}{N_{i,j,k}}\sum_{t_{i,j,\ell} \in S_{i,j,k}}\bigg(\int_{t_{i,j,k}}^{t_{i,j,\ell}}\mu_{i,s}\,{\rm d}s\bigg)^2\\
&+\frac{1}{2n_{i,j}}\sum_{k=1}^{n_{i,j}}\frac{1}{N_{i,j,k}}\sum_{t_{i,j,\ell} \in S_{i,j,k}}\bigg(\int_{t_{i,j,k}}^{t_{i,j,\ell}}\sigma_{j,s}\,{\rm d}B_{j,s}\bigg)^2\,.
\end{align*}
Identical to deriving (\ref{eq:c1}) and (\ref{eq:c2}), it holds that
$
\max_{i,j\in[p]}\mathbb{P}\{|{\rm II}_3(i,j)|>v\}\lesssim \exp(-Cv^2\xi^{-1})+\exp\{-Cv^2(d\xi^{1/2})^{\tau-4}\}+\exp(-Cd^\gamma)
$
for any $v=o\{(d\xi^{1/2})^{2-\tau}\}$ with $d\leq \xi^{-1/2}$. Such upper bound also holds for $\max_{i,j\in[p]}\mathbb{P}\{|{\rm II}_4(i,j)|>v\}$. Together with (\ref{eq:c1}) and (\ref{eq:c2}), we have
\begin{align}\label{eq:lm2f}
\max_{i,j\in[p]}\mathbb{P}\{|{\rm II}(i,j)|>v\}\lesssim\exp(-Cv^2\xi^{-1})+\exp\{-Cv^2(d\xi^{1/2})^{\tau-4}\}+\exp(-Cd^\gamma)
\end{align}
for any $v=o\{(d\xi^{1/2})^{2-\tau}\wedge\xi^{-1}\}$ with $d\leq \xi^{-1/2}$. To make $\max_{i,j\in[p]}|{\rm II}(i,j)|=O_\p(n_*^{-1/2}K^{1/2}\log^{1/2}p)$, it suffices to require $d^{2\tau-4}\xi^{\tau-1}\log p=o(1)$, $\xi^3\log p=o(1)$, $d^{8-2\tau}\xi^{2-\tau}=O(1)$, $d^2\xi\leq 1$ and $\log p=o(d^{\gamma})$. Due to $\tau\in[0,4)$, $\xi=o(1)$ and $d\rightarrow\infty$, $d^{8-2\tau}\xi^{2-\tau}\rightarrow\infty$ when $\tau\in[2,4)$. Thus, we need to restrict $\tau\in[0,2)$, which leads to $\xi^{(\tau-1)/(2-\tau)}(\log p)^{1/(2-\tau)}\ll d^2\ll \xi^{(\tau-2)/(4-\tau)}$ and $\log p=o(d^\gamma \wedge \xi^{-3})$. It follows from $\xi^{(\tau-1)/(2-\tau)}(\log p)^{1/(2-\tau)}\ll \xi^{(\tau-2)/(4-\tau)}$ that $\xi^{\tau/(4-\tau)}\log p=o(1)$.  Selecting $d$ sufficiently close to $\xi^{(\tau-2)/(8-2\tau)}$, we have $\log p =o[\min\{\xi^{-\tau/(4-\tau)},\xi^{\gamma(\tau-2)/(8-2\tau)}\}]$.  To make $p$ diverge as fast as possible, we can choose $\tau=2\gamma/(\gamma+2)$ and $d=\xi^{-1/(4+\gamma)}$. Note that $\xi^{-1}\asymp n_*K^{-1}$. It follows from (\ref{eq:lm2f}) that
$
\max_{i,j\in[p]}\mathbb{P}\{|{\rm II}(i,j)|>v\}\lesssim\exp(-Cn_*K^{-1}v^2)+\exp\{-C(n_*K^{-1})^{\gamma/(\gamma+4)}\}
$
for any $v=o\{(n_*^{-1}K)^{2/(\gamma+4)}\}$.

By \eqref{eq:IIij}, ${\rm II}(i,j)={\rm II}_1(i,j)+{\rm II}_2(i,j)+{\rm II}_3(i,j)+{\rm II}_4(i,j)$.  Notice that for sufficiently large $n$, $\min_{k \in[n_{i,j}]}N_{i,j,k}=\Delta_K+1$ and $\max_{k\in[n_{i,j}]}N_{i,j,k}=2\Delta_K+2$. Since  $\Delta_K$ is a fixed constant, by Assumption \ref{as:drifdiff}, for any fixed positive integer $m$, Jensen's inequality implies that
\begin{align}\label{eq:ub}
\max_{ i,j\in[p]}\mathbb{E}\{|{\rm II}_1(i,j)|^{m}\}\lesssim&~ \max_{ i,j\in[p]}\max_{ k\in[n_{i,j}]}\max_{t_{i,j,\ell}\in S_{i,j,k}}\mathbb{E}\bigg\{\bigg(\int_{t_{i,j,k}}^{t_{i,j,\ell}}\mu_{i,s}\,{\rm d}s\bigg)^{2m}+\bigg(\int_{t_{i,j,k}}^{t_{i,j,\ell}}\mu_{j,s}\,{\rm d}s\bigg)^{2m}\bigg\}
\notag\\\lesssim &~
\xi^{2m-1}\max_{ i,j\in[p]}\max_{ k\in[n_{i,j}]}\max_{t_{i,j,\ell}\in S_{i,j,k}} \int_{t_{i,j,k}\wedge t_{i,j,\ell}}^{t_{i,j,k}\vee t_{i,j,\ell}}\sup_{0\leq s\leq  T}\mathbb{E}(|\mu_{i,s}|^{2m})\,{\rm d}s\lesssim \xi^{2m}\,.
\end{align}
Meanwhile, by Assumption \ref{as:drifdiff} and Burkholder-Davis-Gundy inequality, it holds that
\begin{align}\label{eq:sigmab}
\max_{ i,j\in[p]}\mathbb{E}\{|{\rm II}_2(i,j)|^{m}\}\lesssim&~ \max_{ i,j\in[p]}\max_{ k\in[n_{i,j}]}\max_{t_{i,j,\ell}\in S_{i,j,k}}\mathbb{E}\bigg\{\bigg(\int_{t_{i,j,k}}^{t_{i,j,\ell}}\sigma_{i,s}\,{\rm d}B_{i,s}\bigg)^{2m}\bigg\}\lesssim \xi^{m}\,.
\end{align}
Due to
$
\max_{ i,j\in[p]}\mathbb{E}\{|{\rm II}_3(i,j)|^{m}+|{\rm II}_4(i,j)|^{m}\}\lesssim\max_{ i,j\in[p]}\max_{ k\in[n_{i,j}]}\max_{t_{i,j,\ell}\in S_{i,j,k}}\mathbb{E}\{(\int_{t_{i,j,k}}^{t_{i,j,\ell}}\mu_{i,s}\,{\rm d}s)^{2m}+(\int_{t_{i,j,k}}^{t_{i,j,\ell}}\sigma_{j,s}\,{\rm d}B_{j,s})^{2m}\}$,
together with \eqref{eq:ub} and \eqref{eq:sigmab},
$\max_{ i,j\in[p]}\mathbb{E}\{|{\rm II}_3(i,j)|^{m}+|{\rm II}_4(i,j)|^{m}\}\lesssim \xi^{m}$.  Hence, $\max_{ i,j\in[p]}\mathbb{E}\{|{\rm II}(i,j)|^{m}\}\lesssim \xi^{m}$.
 We complete the proof of  Proposition \ref{la:2}. $\hfill\Box$

\subsection{Proof of Proposition \ref{la:3}}\label{sec:pfpn3}

To prove Proposition \ref{la:3}, we need the following lemma whose proof is given in Section \ref{se:pflemma2}.
\begin{lemma}\label{eq:mixingbound}
	Let $\{z_t\}_{t=1}^{\tilde{n}} $ be an $\alpha$-mixing sequence of real-valued and centered random variables with $\alpha$-mixing coefficients $\{\alpha(k)\}_{k\geq1}$.  Assume there exist some universal  constants $a_1>1$, $a_2>0$, $\tilde{r}\geq0$ and $\varphi>0$ such that {\rm(i)} $\max_{t\in[\tilde{n}]}\mathbb{E}(|z_t|^k)\leq (k!)^{1+\tilde{r}}H_{\tilde{n}}^k$ for any integer $k\geq2$, where $H_{\tilde{n}}>0$ may diverge with $\tilde{n}$, {\rm(ii)} $\alpha(k)\leq a_1\exp\{-a_2(\tilde{L}_{\tilde{n}}^{-1}k)^{\varphi}\} $ for any integer $k\geq 1$, where $ \tilde{L}_{\tilde{n}}>0$ may diverge with $\tilde{n}$. Let $S_{\tilde{n}}=\sum_{t=1}^{\tilde{n}}z_t$.  It holds that 	\begin{align*}
	\mathbb{P}(|S_{\tilde{n}}|\geq \tilde{n}x)\lesssim \exp(-C\tilde{n}\tilde{L}_{\tilde{n}}^{-1}H_{\tilde{n}}^{-2}x^2)+\exp\{-C(\tilde{n}\tilde{L}_{\tilde{n}}^{-1}H_{\tilde{n}}^{-1}x)^{1/(1+\check{r})}\}
	\end{align*}
	for any $x>0$, where $\check{r}=1+\tilde{r}+\varphi^{-1}$.
\end{lemma}

\begin{rek}
If $\varphi=\infty$, the upper bound in Lemma {\rm\ref{eq:mixingbound}} holds with $\check{r}=1+\tilde{r}$.
\end{rek}

By the definition of ${\rm III}(i,j)$, we can reformulate it as
\[
\begin{split}
{\rm III}(i,j)=\frac{1}{2n_{i,j}}\sum_{k=1}^{n_{i,j}}\sum_{t_{i,j,\ell}\in S_{i,j,k}}\bigg(\frac{1}{N_{i,j,\ell}}+\frac{1}{N_{i,j,k}}\bigg)(X_{i,t_{i,j,k}}-X_{i,t_{i,j,\ell}})U_{j,t_{i,j,k}}\,.
\end{split}
\]
For each $i,j\in[p]$ and $k\in[n_{i,j}]$, define
$
G_{i,j,k}=\sum_{t_{i,j,\ell}\in S_{i,j,k}}(N_{i,j,\ell}^{-1}+N_{i,j,k}^{-1})(X_{i,t_{i,j,k}}-X_{i,t_{i,j,\ell}})$
and $D_{i,j}=\max_{ k\in[n_{i,j}]}|G_{i,j,k}|$. Recall $\xi=\max_{i,j\in[p]}\max_{ k\in[n_{i,j}]}\max_{t_{i,j,\ell}\in S_{i,j,k}}|t_{i,j,\ell}-t_{i,j,k}|$.

We will first consider the tail probability $\mathbb{P}(D_{i,j}>v)$. By Bonferroni inequality, we have
\begin{align}\label{eq:di}
\mathbb{P}(D_{i,j}>v)\leq&~\sum_{k=1}^{n_{i,j}}\mathbb{P}\bigg\{\bigg|\sum_{t_{i,j,\ell}\in S_{i,j,k}}\bigg(\frac{1}{N_{i,j,\ell}}+\frac{1}{N_{i,j,k}}\bigg)(X_{i,t_{i,j,k}}-X_{i,t_{i,j,\ell}})\bigg|>v\bigg\}\notag\\
\leq&~\sum_{k=1}^{n_{i,j}}\mathbb{P}\bigg\{\bigg|\sum_{t_{i,j,\ell}\in S_{i,j,k}}\bigg(\frac{1}{N_{i,j,\ell}}+\frac{1}{N_{i,j,k}}\bigg)\int_{t_{i,j,\ell}}^{t_{i,j,k}}\mu_{i,s}\,{\rm d}s\bigg|>\frac{v}{2}\bigg\}\\
&+\sum_{k=1}^{n_{i,j}}\mathbb{P}\bigg\{\bigg|\sum_{t_{i,j,\ell}\in S_{i,j,k}}\bigg(\frac{1}{N_{i,j,\ell}}+\frac{1}{N_{i,j,k}}\bigg)\int_{t_{i,j,\ell}}^{t_{i,j,k}}\sigma_{i,s}\,{\rm d}B_{i,s}\bigg|>\frac{v}{2}\bigg\}\,.\notag
\end{align}
For any $\theta>0$, by Triangle inequality and Jensen's inequality, similar to \eqref{eq:b1}, it holds that
$\max_{ k\in[n_{i,j}]}\mathbb{E}[\exp\{\theta|\sum_{t_{i,j,\ell}\in S_{i,j,k}}(N_{i,j,\ell}^{-1}+N_{i,j,k}^{-1})\int_{t_{i,j,\ell}}^{t_{i,j,k}}\mu_{i,s}\,{\rm d}s|\}]\leq \sup_{0\leq t\leq T}\mathbb{E}\{\exp(C\theta\xi|\mu_{i,t}|)\}$.
It follows from Assumption \ref{as:drifdiff} that
$
\sup_{0\leq t\leq T}\mathbb{E}\{\exp(C\theta\xi|\mu_{i,t}|)\}\leq C\exp(C\xi^2\theta^2)$.
Selecting $\theta\asymp \xi^{-1}$ and applying Markov's inequality, we have
\begin{align}
&\max_{ k\in[n_{i,j}]}\mathbb{P}\bigg\{\bigg|\sum_{t_{i,j,\ell}\in S_{i,j,k}}\bigg(\frac{1}{N_{i,j,\ell}}+\frac{1}{N_{i,j,k}}\bigg)\int_{t_{i,j,\ell}}^{t_{i,j,k}}\mu_{i,s}\,{\rm d}s\bigg|>\frac{v}{2}\bigg\}\notag\\
&~~~~~~~~~\leq \exp(-\theta v)\max_{k\in[n_{i,j}]}\mathbb{E}\bigg[\exp\bigg\{2\theta\bigg|\sum_{t_{i,j,\ell}\in S_{i,j,k}}\bigg(\frac{1}{N_{i,j,\ell}}+\frac{1}{N_{i,j,k}}\bigg)\int_{t_{i,j,\ell}}^{t_{i,j,k}}\mu_{i,s}\,{\rm d}s\bigg|\bigg\}\bigg]\label{eq:di1}\\
&~~~~~~~~~\lesssim\exp(-C\xi^{-1}v)\notag
\end{align}
for any $v>0$. For any constant $d\in(0,\xi^{-1/2}]$, define a stopping time $
\Gamma_{i,d}=T\wedge\inf\{t>0:\sup_{0\leq s\leq t}\sigma_{i,s}>d\}$. By Cauchy-Schwarz inequality,
$
|\sum_{t_{i,j,\ell}\in S_{i,j,k}}(N_{i,j,\ell}^{-1}+N_{i,j,k}^{-1})\int_{t_{i,j,\ell}}^{t_{i,j,k}}\sigma_{i,s}\,{\rm d}B_{i,s}|^2\lesssim N_{i,j,k}^{-1}\sum_{t_{i,j,\ell}\in S_{i,j,k}}(\int_{t_{i,j,\ell}}^{t_{i,j,k}}\sigma_{i,s}\,{\rm d}B_{i,s})^2$,
which implies that
\begin{equation}\label{eq:eb1}
\begin{split}
&\mathbb{P}\bigg\{\bigg|\sum_{t_{i,j,\ell}\in S_{i,j,k}}\bigg(\frac{1}{N_{i,j,\ell}}+\frac{1}{N_{i,j,k}}\bigg)\int_{t_{i,j,\ell}}^{t_{i,j,k}}\sigma_{i,s}\,{\rm d}B_{i,s}\bigg|>\frac{v}{2},\Gamma_{i,d}=T\bigg\}\\
&~~~~~~~~\lesssim \exp(-C\theta v^2)\mathbb{E}\bigg[\exp\bigg\{\frac{\theta}{N_{i,j,k}}\sum_{t_{i,j,\ell}\in S_{i,j,k}}\bigg(\int_{t_{i,j,\ell}}^{t_{i,j,k}}\sigma_{i,s}\,{\rm d}B_{i,s}\bigg)^2\bigg\}{I}(\Gamma_{i,d}=T)\bigg]
\end{split}
\end{equation}
for any $\theta>0$. By Jensen's inequality,
$\mathbb{E}[\exp\{\theta N_{i,j,k}^{-1}\sum_{t_{i,j,\ell}\in S_{i,j,k}}(\int_{t_{i,j,\ell}}^{t_{i,j,k}}\sigma_{i,s}\,{\rm d}B_{i,s})^2\}{I}(\Gamma_{i,d}=T)]\leq N_{i,j,k}^{-1}\sum_{t_{i,j,\ell}\in S_{i,j,k}}\mathbb{E}[\exp\{\theta(\int_{t_{i,j,\ell}}^{t_{i,j,k}}\sigma_{i,s}\,{\rm d}B_{i,s})^2\}{I}(\Gamma_{i,d}=T)]$.
Same as (\ref{eq:re1}) and (\ref{eq:re2}), for any $\theta\in(0,d^{-2}\xi^{-1}/4]$, it holds that
$
\mathbb{E}[\exp\{\theta(\int_{t_{i,j,\ell}}^{t_{i,j,k}}\sigma_{i,s}\,{\rm d}B_{i,s})^2\}{I}(\Gamma_{i,d}=T)]\leq \exp(Cd^4\xi^2\theta^2)\lesssim 1$.
Selecting $\theta=d^{-2}\xi^{-1}/4$, together with (\ref{eq:eb1}), we have
$
\max_{ k\in[n_{i,j}]}\mathbb{P}\{|\sum_{t_{i,j,\ell}\in S_{i,j,k}}(N_{i,j,\ell}^{-1}+N_{i,j,k}^{-1})\int_{t_{i,j,\ell}}^{t_{i,j,k}}\sigma_{i,s}\,{\rm d}B_{i,s}|>v/2,\Gamma_{i,d}=T\}\lesssim \exp(-Cd^{-2}\xi^{-1}v^2)$
for any $v>0$. It follows from Assumption \ref{as:diff} that
\begin{align*}
&\max_{ k\in[n_{i,j}]}\mathbb{P}\bigg\{\bigg|\sum_{t_{i,j,\ell}\in
S_{i,j,k}}\bigg(\frac{1}{N_{i,j,\ell}}+\frac{1}{N_{i,j,k}}\bigg)\int_{t_{i,j,\ell}}^{t_{i,j,k}}\sigma_{i,s}\,{\rm d}B_{i,s}\bigg|>\frac{v}{2}\bigg\}\\
&~~~~~\leq\max_{k\in[n_{i,j}]}\mathbb{P}\bigg\{\bigg|\sum_{t_{i,j,\ell}\in S_{i,j,k}}\bigg(\frac{1}{N_{i,j,\ell}}+\frac{1}{N_{i,j,k}}\bigg)\int_{t_{i,j,\ell}}^{t_{i,j,k}}\sigma_{i,s}\,{\rm d}B_{i,s}\bigg|>\frac{v}{2},\Gamma_{i,d}=T\bigg\}+\mathbb{P}(\Gamma_{i,d}\neq T)\\
&~~~~~\lesssim\exp(-Cd^{\gamma})+\exp(-Cd^{-2}\xi^{-1}v^2)
\end{align*}
for any $v>0$. Letting $d\rightarrow\infty$, together with (\ref{eq:di1}), (\ref{eq:di}) implies that
\begin{equation}\label{eq:D}
\max_{i,j\in[p]}\mathbb{P}(D_{i,j}>v)\lesssim n_*\exp(-Cd^{\gamma})+n_*\exp(-Cd^{-2}\xi^{-1}v^2)
\end{equation}
for any $0<v\leq C$.
 Let $\bar{\mathbb{E}}(\cdot)=\mathbb{E}\{\cdot\,|\,(\bX_t)_{t\in[0,T]}\}$. For any integer $s\geq 2$, we have
$
\bar{\mathbb{E}}(|G_{i,j,k}U_{j,t_{i,j,k}}|^s)\leq  D_{i,j}^s\cdot\mathbb{E}(|U_{j,t_{i,j,k}}|^s)\leq D_{i,j}^s\cdot C^ss^{(s+1)/2}\leq s!\cdot(CD_{i,j})^s$.
By Lemma \ref{eq:mixingbound} with $\tilde{L}_{\tilde{n}}=L_n$,
$
\mathbb{P}\{|{\rm{III}}(i,j)|\geq x\,|\,(\bX_t)_{t\in[0,T]}\}\lesssim \exp(-Cn_*L_n^{-1}D_{i,j}^{-2}x^2)+\exp\{-C(n_*L_n^{-1}D_{i,j}^{-1}x)^{\varphi/(2\varphi+1)}\}$ for any $x>0$,
which implies that
\begin{align*}
\mathbb{P}\{|{\rm{III}}(i,j)|\geq x\}\lesssim \mathbb{E}\bigg\{\exp\bigg(-\frac{Cn_*L_n^{-1}x^2}{D_{i,j}^2}\bigg)\bigg\}+\mathbb{E}\bigg[\exp\bigg\{-C\bigg(\frac{n_*L_n^{-1}x}{D_{i,j}}\bigg)^{\varphi/(2\varphi+1)}\bigg\}\bigg]
\end{align*}
for any $x>0$. Therefore, from (\ref{eq:D}) with $v\asymp 1$, it holds that
\begin{align*}
\mathbb{E}\biggl\{\exp\bigg(-\frac{Cn_*L_n^{-1}x^2}{D_{i,j}^2}\bigg)\biggr\}\lesssim&~ \exp(-Cv^{-2}n_*L_n^{-1}x^2)+\mathbb{P}(D_{i,j}>v)\\\lesssim&~
\exp(-Cn_*L_n^{-1}x^2)+n_*\exp(-Cd^{\gamma})+n_*\exp(-Cd^{-2}\xi^{-1})\,.
\end{align*}
Analogously, we also have
\begin{align*}
\mathbb{E}\bigg[\exp\bigg\{-C\bigg(\frac{n_*L_n^{-1}x}{D_{i,j}}\bigg)^{\varphi/(2\varphi+1)}\bigg\}\bigg]\lesssim&~
\exp\{-C(n_*L_n^{-1}x)^{\varphi/(2\varphi+1)}\}\\
&+n_*\exp(-Cd^{\gamma})+n_*\exp(-Cd^{-2}\xi^{-1})\,.
\end{align*}
Thus, it holds that $
\max_{ i,j\in[p]}\mathbb{P}\{|{\rm{III}}(i,j)|\geq x\}\lesssim\exp(-Cn_*L_n^{-1}x^2)+\exp\{-C(n_*L_n^{-1}x)^{\varphi/(2\varphi+1)}\}+n_*\exp(-Cd^{\gamma})+n_*\exp(-Cd^{-2}\xi^{-1})$
for any $x>0$.  Recall that $\xi\asymp Kn_*^{-1}$.
To make $\max_{i,j\in[p]}|{\rm{III}}(i,j)|=O_\p(n_*^{-1/2}K^{1/2}\log ^{1/2}p)$, it suffices to require $K\gtrsim L_n$, $\log p=o[\min\{d^{\gamma},d^{-2}n_*K^{-1},(n_*L_n^{-2}K)^{\varphi/(3\varphi+2)}\}]$ and $\log n_*=o[\min\{d^\gamma, d^{-2}n_*K^{-1}\}]$.  In order to make $p$ diverge as fast as possible, we can select $d=(n_*K^{-1})^{1/(2+\gamma)}$. If $K(\log n_*)^{1+2/\gamma}=o(n_*)$, then
$
\max_{i,j\in[p]}\mathbb{P}\{|{\rm{III}}(i,j)|\geq v\}\lesssim\exp(-Cn_*L_n^{-1}v^2)+\exp\{-C(n_*L_n^{-1}v)^{\varphi/(2\varphi+1)}\}
+\exp\{-C(n_*K^{-1})^{\gamma/(2+\gamma)}\}$
for any $v>0$.

By Jensen's inequality and Cauchy-Schwarz inequality, it holds that
\begin{align*}
\mathbb{E}\{|{\rm III}(i,j)|^{m}\}\lesssim&~
\max_{ k\in[n_{i,j}]}\max_{t_{i,j,\ell}\in S_{i,j,k}}\mathbb{E}(|X_{i,t_{i,j,k}}-X_{i,t_{i,j,\ell}}|^{m}|U_{j,t_{i,j,k}}|^{m})
\\= &~
\max_{ k\in[n_{i,j}]}\max_{t_{i,j,\ell}\in S_{i,j,k}}\mathbb{E}\bigg(\bigg|\int_{t_{i,j,\ell}}^{t_{i,j,k}}\mu_{i,s}\,{\rm d}s+\int_{t_{i,j,\ell}}^{t_{i,j,k}}\sigma_{i,s}\,{\rm d}B_{i,s}\bigg|^{m}|U_{j,t_{i,j,k}}|^{m}\bigg)
\\\lesssim &~
\max_{ k\in[n_{i,j}]}\max_{t_{i,j,\ell}\in S_{i,j,k}}\bigg\{\mathbb{E}\bigg(\bigg|\int_{t_{i,j,\ell}}^{t_{i,j,k}}\mu_{i,s}\,{\rm d}s\bigg|^{2m}+\bigg|\int_{t_{i,j,\ell}}^{t_{i,j,k}}\sigma_{i,s}\,{\rm d}B_{i,s}\bigg|^{2m}\bigg)\bigg\}^{1/2}\{\mathbb{E}(U_{j,t_{i,j,k}}^{2m})\}^{1/2}\,.
\end{align*}
By Assumption \ref{as:moment}, $\max_{ k\in[n_{i,j}]}\mathbb{E}(U_{j,t_{i,j,k}}^{2m})\lesssim1$. Together with \eqref{eq:ub} and \eqref{eq:sigmab},  $\max_{ i,j\in[p]}\mathbb{E}\{|{\rm III}(i,j)|^{m}\}\lesssim \xi^{m/2}$.
 Analogously, we can show the same results hold for ${\rm{IV}}(i,j)$. We complete the proof of Proposition \ref{la:3}. $\hfill\Box$

\subsection{Proof of Theorem \ref{tm:2}}

To prove Theorem \ref{tm:2}, we need the Le Cam's lemma as stated in Lemma \ref{la:lecam} below. Its proof can be found in \cite{LeCam_1973} and \cite{DonohoLiu_1991}. Let $\mathcal{Z}$ be an observation from a distribution $\mathbb{P}_\theta$ where $\theta$ belongs to a parameter space $\Theta$. For two distributions $\mathbb{Q}_0$ and $\mathbb{Q}_1$ with densities $q_0$ and $q_1$ with respect to any common dominating measure $\mu$, the total variation affinity is given by $\|\mathbb{Q}_0\wedge\mathbb{Q}_1\|=\int q_0\wedge q_1\,{\rm d}\mu$. Let $\Theta=\{\theta_0,\theta_1,\ldots,\theta_D\}$ and denote by $L$ the loss function. Define $l_{\min}=\min_{d\in[D]}\inf_t\{L(t,\theta_0)+L(t,\theta_d)\}$ and denote $\bar{\mathbb{P}}=D^{-1}\sum_{d=1}^D\mathbb{P}_{\theta_d}$.

\begin{lemma}\label{la:lecam}
{\rm(Le Cam's lemma)} Let $T$ be any estimator of $\theta$ based on an observation $\mathcal{Z}$ from a distribution ${P}_\theta$ with $\theta\in\Theta=\{\theta_0,\theta_1,\ldots,\theta_D\}$, then
$
\sup_{\theta\in\Theta}\mathbb{E}_{\mathcal{Z}|\theta}\{L(T,\theta)\}\geq{2}^{-1}l_{\min}\|\mathbb{P}_{\theta_0}\wedge\bar{\mathbb{P}}\|$.
\end{lemma}

For each $k\in[n]$, define $\mathcal{C}_k=\{i\in[p]:t_k\in\mathcal{G}_i\}$ where $\mathcal{G}_i$ is the grid of time points where
we observe the noisy data of the $i$th component process. For any $s$-dimensional vector $\ba$ and an index set $\mathcal{C}\subset[s]$, denote by $\ba_{\mathcal{S}}$ the subvector of $\ba$ with components indexed by $\mathcal{C}$. The data we have is $\mathcal{Z}=\{\bY_{t_1,\mathcal{C}_1},\ldots,\bY_{t_n,\mathcal{C}_n}\}$. Select the loss function $
L(T,\theta)=\max_{i,j\in[p]}|\omega_{i,j}-\theta_{i,j}|$
 for any $T=(\omega_{i,j})_{p\times p}$ and $\theta=(\theta_{i,j})_{p\times p}\in\Theta$.
Select $D=p$, $\theta_0=\bSigma_{u,0}=\bI_p$ and \[
\theta_d=\bSigma_{u,d}=\bI_p+(v^{1/2}n_*^{-1/2}\log^{1/2}D)\textrm{diag}(\underbrace{0,\ldots,0}_{d-1},1,\underbrace{0,\ldots,0}_{p-d})\,,
\]
for any $d\in[D]$, where $v>0$ is a sufficiently small constant. For each $d=0,1,\ldots,D$, we write $\theta_d=\bSigma_{u,d}=(\sigma_{u,i,j,d})_{p\times p}$. Then
\begin{equation}\label{eq:mini1}
\begin{split}
l_{\min}=&~\min_{d\in[D]}\inf_t\{L(t,\theta_0)+L(t,\theta_d)\}\\
\geq&~\min_{d\in[D]}\max_{i,j\in[p]}|\sigma_{u,i,j,0}-\sigma_{u,i,j,d}|\geq\sqrt{\frac{v\log p}{n_*}}\asymp \sqrt{\frac{\log p}{n_*}}\,.
\end{split}
\end{equation}

To prove the lower bound stated in Theorem \ref{tm:2}, it suffices to construct a specific model which makes the stated lower bound be achievable. To do this, we select $\mu_{i,t}=0$ and $\sigma_{i,t}=0$ for any $t\in[0,T]$. Then the associated $\bX_t={\bf0}$ for any $t\in[0,T]$. In this special case, $\bY_{t_k}=\bU_{t_k}$. Given $(n,n_*)$ with $n\geq n_*$, and $0\leq t_1<\cdots<t_n=T$, we define $\mathcal{G}_*=\{\tilde{t}_1,\ldots,\tilde{t}_{n_*}\}$ with each $\tilde{t}_j\in\{t_1,\ldots,t_n\}$ and $\tilde{t}_j<\tilde{t}_{j+1}$. For each $t_j\in\mathcal{G}_*$, we assume all $p$ component processes are observed. For any $t_j\notin\mathcal{G}_*$, we assume only one component process are observed. Without loss of generality, we assume $\mathcal{G}_*=\{t_1,\ldots,t_{n_*}\}$. Let $n-n_*=ap+q_*$ where $a\geq 0$ and $0\leq q_*<p$ are two integers. We assume the $i$th component process is observed at $t_{n_*+jp+i}$'s with $j=0,\ldots,a$ and $i\in[p]$.
Then
$
\mathcal{G}_i=\mathcal{G}_*\cup\{t_{n_*+i},\ldots,t_{n_*+ap+i}\}$.

Let $\bU_{t_k}\sim_{{\rm i.i.d.}} N(\bzero,\bSigma_{u,d})$, and denote the joint density of $\bU_{t_1,\mathcal{C}_1},\ldots,\bU_{t_n,\mathcal{C}_n}$ by $f_d$. Denote by $\phi_\sigma$ the density of $N(0,\sigma^2)$. Write $\sigma_*^2=1+v^{1/2}n_*^{-1/2}\log^{1/2}D$. Then
$
f_0=\prod_{k=1}^n\prod_{j\in\mathcal{C}_k}\phi_1(u_{k,j})
$
and
$
f_d=\prod_{k=1}^n\prod_{j\in\mathcal{C}_k\backslash\{d\}}\phi_1(u_{k,j})\cdot\prod_{k=1}^n\prod_{d\in\mathcal{C}_k}\phi_{\sigma_*}(u_{k,d})
$
for each $d\in[D]$. Here we adopt the convention $\prod_{d\in\mathcal{C}_k}\phi_{\sigma_*}(u_{k,d})\equiv1$ if $d\notin\mathcal{C}_k$.
We will show $\|\mathbb{P}_{\theta_0}\wedge\bar{\mathbb{P}}\|\geq c$ for some uniform constant $c>0$.

For any two densities $q_0$ and $q_1$, by Cauchy-Schwarz inequality, we have
$
(\int |q_0-q_1|\,{\rm d}\mu)^2\leq\int{(q_0-q_1)^2}/{q_1}\,{\rm d}\mu=\int{q_0^2}/{q_1}\,{\rm d}\mu-1$,
which implies that
$
\int q_0\wedge q_1\,{\rm d}\mu=1-{2}^{-1}\int|q_0-q_1|\,{\rm d}\mu\geq 1-{2}^{-1}(\int {q_0^2}/{q_1}\,{\rm d}\mu-1)^{1/2}$.
 In order to show $\|\mathbb{P}_{\theta_0}\wedge\bar{\mathbb{P}}\|\geq c$ for some uniform constant $c>0$, it suffices to show that $\int(D^{-1}\sum_{d=1}^Df_d)^2f_0^{-1}\,{\rm d}\mu-1\rightarrow0$, that is,
\begin{equation}\label{eq:minitoshow}
\frac{1}{D^2}\sum_{d=1}^D\bigg(\int\frac{f_d^2}{f_0}\,{\rm d}\mu-1\bigg)+\frac{1}{D^2}\sum_{d_1\neq d_2}\bigg(\int\frac{f_{d_1}f_{d_2}}{f_0}\,{\rm d}\mu-1\bigg)\rightarrow0\,.
\end{equation}
Notice that
$
f_{d_1}f_{d_2}/f_0=\prod_{k=1}^n[\prod_{d_1\in\mathcal{C}_k}\phi_{\sigma_*}(u_{k,d_1})\cdot\prod_{d_2\in\mathcal{C}_k}\phi_{\sigma_*}(u_{k,d_2})\cdot\prod_{j\in\mathcal{C}_k\backslash\{ d_1,d_2\}}\phi_{1}(u_{k,j})]
$ for any $d_1\neq d_2$,
then $
\int f_{d_1}f_{d_2}/f_0\,{\rm d}\mu=1$,
which implies
$
D^{-2}\sum_{d_1\neq d_2}(\int{f_{d_1}f_{d_2}}/{f_0}\,{\rm d}\mu-1)=0$.
For any $d=1,\ldots,D$, we have
\[
\frac{f_d^2}{f_0}=\prod_{k=1}^n\prod_{j\in\mathcal{C}_k\backslash\{d\}}\phi_1(u_{k,j})\cdot\prod_{k=1}^n\bigg[\frac{1}{\sqrt{2\pi}\sigma_*^2}\exp\bigg\{-\frac{(2-\sigma_*^2)u_{k,d}^2}{2\sigma_*^2}\bigg\}\bigg]^{{I}(d\in\mathcal{C}_k)}\,,
\]
which implies
\begin{align*}
\int\frac{f_d^2}{f_0}\,{\rm d}\mu=&~\bigg(\frac{1}{\sigma_*\sqrt{2-\sigma_*^2}}\bigg)^{\sum_{k=1}^n{I}(d\in\mathcal{C}_k)}\prod_{k=1}^n\prod_{j\in\mathcal{C}_k\backslash\{d\}}\bigg\{\int\phi_1(u_{k,j})\,{\rm d}u_{k,j}\bigg\}\\
&~~~~~~~~~~~~~~~\times\prod_{k=1}^n\bigg[\int\frac{\sqrt{2-\sigma_*^2}}{\sqrt{2\pi}\sigma_*}\exp\bigg\{-\frac{(2-\sigma_*^2)u_{k,d}^2}{2\sigma_*^2}\bigg\}\,{\rm d}u_{k,d}\bigg]^{{I}(d\in\mathcal{C}_k)}\\
=&~\bigg(\frac{1}{\sigma_*\sqrt{2-\sigma_*^2}}\bigg)^{\sum_{k=1}^n{I}(d\in\mathcal{C}_k)}=\bigg(1-\frac{v\log D}{n_*}\bigg)^{-\sum_{k=1}^n{I}(d\in\mathcal{C}_k)/2}\,.
\end{align*}
Notice that $\sum_{k=1}^n{I}(d\in\mathcal{C}_k)\leq n_*+a+1$ for each $d\in[D]$. Therefore,
$
\int{f_d^2}/{f_0}\,{\rm d}\mu\leq (1-vn_*^{-1}\log D)^{-(n_*+a+1)/2}
$
for each $d\in[D]$. Due to $n/n_*\lesssim p$, we know $a\lesssim n_*$. Applying the inequality $\log(1-x)\geq -2x$ for any $0<x<1/2$, we have
\[
\begin{split}
0\leq&~\frac{1}{D^2}\sum_{d=1}^D\bigg(\int\frac{f_d^2}{f_0}\,{\rm d}\mu-1\bigg)\leq\exp\bigg[-\bigg\{1-v\bigg(1+\frac{a+1}{n_*}\bigg)\bigg\}\log D\bigg]\rightarrow0
\end{split}
\]
for sufficiently small $v>0$. Then (\ref{eq:minitoshow}) holds. Hence $\|\mathbb{P}_{\theta_0}\wedge\bar{\mathbb{P}}\|\geq c$ for some uniform constant $c>0$. Together with (\ref{eq:mini1}), we can obtain Theorem \ref{tm:2} by Lemma \ref{la:lecam}. $\hfill\Box$

\subsection{Proof of Theorem \ref{tm:3}}\label{sec:pfthm3}

 We first consider the case with $\varphi<\infty$. Write $\aleph=(n_*^{-1}K\log p)^{1/2}$. For each $i,j\in[p]$, we  define the event
$
A_{i,j}=\{|\hat{\sigma}_{u,i,j}^{\rm thre}-\sigma_{u,i,j}|\leq 4\min(|\sigma_{u,i,j}|,\alpha\aleph)\}
$
with some $\alpha>0$, and
$
d_{i,j}=(\hat{\sigma}_{u,i,j}^{\rm thre}-\sigma_{u,i,j}){I}(A_{i,j}^c)$.
Write $\bD=(d_{i,j})_{p\times p}$. Due to $\|\bW\|_2\leq \|\bW\|_\infty$ for any $p\times p$ symmetric matrix $\bW$, it holds that
\begin{equation}\label{eq:speupp}
\begin{split}
\|\widehat{\bSigma}_{u}^{\rm thre}-\bSigma_{u}\|_2^2\leq&~\bigg\{\max_{i\in[p]}\sum_{j=1}^p|\hat{\sigma}_{u,i,j}^{\rm thre}-\sigma_{u,i,j}|\bigg\}^2\\
\leq&~2\bigg(\max_{i\in[p]}\sum_{j=1}^p|d_{i,j}|\bigg)^2+2\bigg\{\max_{i\in[p]}\sum_{j=1}^p|\hat{\sigma}_{u,i,j}^{\rm thre}-\sigma_{u,i,j}|{I}(A_{i,j})\bigg\}^2\,.
\end{split}
\end{equation}
For the second term on the right-hand side of (\ref{eq:speupp}), we have that
$
\sum_{j=1}^p|\hat{\sigma}_{u,i,j}^{\rm thre}-\sigma_{u,i,j}|{I}(A_{i,j})\leq4\sum_{j=1}^{p}\alpha\aleph{I}(|\sigma_{u,i,j}|>\alpha\aleph)+4\sum_{j=1}^p|\sigma_{u,i,j}|{I}(|\sigma_{u,i,j}|\leq\alpha\aleph)$. 
Due to $\sum_{j=1}^p|\sigma_{u,i,j}|^q\leq c_p$, we then have
$
\sum_{j=1}^p|\sigma_{u,i,j}|{I}(|\sigma_{u,i,j}|\leq\alpha\aleph)\leq\sum_{j=1}^p|\sigma_{u,i,j}|^q(\alpha\aleph)^{1-q}\leq\alpha^{1-q}c_p\aleph^{1-q}$
and
$
\sum_{j=1}^{p}\alpha\aleph{I}(|\sigma_{u,i,j}|>\alpha\aleph)\leq\sum_{j=1}^p|\sigma_{u,i,j}|^q(\alpha\aleph)^{1-q}\leq\alpha^{1-q}c_p\aleph^{1-q}$.
 Therefore, we have that
$
\sum_{j=1}^p|\hat{\sigma}_{u,i,j}^{\rm thre}-\sigma_{u,i,j}|{I}(A_{i,j})\leq 8\alpha^{1-q}c_p\aleph^{1-q}$ holds uniformly over $i\in[p]$.
It follows from (\ref{eq:speupp}) that
\begin{equation}\label{eq:speupp1}
\mathbb{E}\big(\|\widehat{\bSigma}_{u}^{\rm thre}-\bSigma_{u}\|_2^2\big)\lesssim \mathbb{E}\bigg\{\bigg(\max_{i\in[p]}\sum_{j=1}^p|d_{i,j}|\bigg)^2\bigg\}+c_p^2\aleph^{2(1-q)}\,.
\end{equation}
Recall $\hat{\sigma}_{u,i,j}$ is defined as \eqref{eq:estasyn}. It holds that
\begin{align}\label{eq:Edij}
\mathbb{E}\bigg\{\bigg(\max_{i\in[p]}\sum_{j=1}^p|d_{i,j}|\bigg)^2\bigg\}\leq&~p\sum_{i,j=1}^p\mathbb{E}\big\{|\hat{\sigma}_{u,i,j}^{\rm thre}-\sigma_{u,i,j}|^2{I}(A_{i,j}^c)\big\}\notag\\
=&~\underbrace{p\sum_{i,j=1}^p\mathbb{E}\big(|\hat{\sigma}_{u,i,j}^{\rm thre}-\sigma_{u,i,j}|^2{I}[A_{i,j}^c\cap\{\hat{\sigma}_{u,i,j}^{\rm thre}=0\}]\big)}_{\textrm{I}}\\
&+\underbrace{p\sum_{i,j=1}^p\mathbb{E}\big(|\hat{\sigma}_{u,i,j}^{\rm thre}-\sigma_{u,i,j}|^2{I}[A_{i,j}^c\cap\{\hat{\sigma}_{u,i,j}^{\rm thre}=\hat{\sigma}_{u,i,j}\}]\big)}_{\textrm{II}}\,.\notag
\end{align}
Recall $\hat{\sigma}_{u,i,j}^{\rm thre}=\hat{\sigma}_{u,i,j}{I}(|\hat{\sigma}_{u,i,j}|\geq \beta\aleph)$ for any $i,j\in[p]$. Then
\begin{align}\label{eq:th3I}
{\rm I}=&~p\sum_{i,j=1}^p\sigma_{u,i,j}^2\mathbb{P}\big[\{|\sigma_{u,i,j}|\geq4\alpha\aleph\}\cap\{|\hat{\sigma}_{u,i,j}|< \beta\aleph\}\big]\notag\\
\leq&~p\sum_{i,j=1}^p\sigma_{u,i,j}^2\mathbb{P}\big[\{|\sigma_{u,i,j}|\geq4\alpha\aleph\}\cap\{|\sigma_{u,i,j}|-|\hat{\sigma}_{u,i,j}-\sigma_{u,i,j}|< \beta\aleph\}\big]\\
\leq&~p\sum_{i,j=1}^p\sigma_{u,i,j}^2\mathbb{P}\big\{|\hat{\sigma}_{u,i,j}-\sigma_{u,i,j}|\geq(4\alpha-\beta)\aleph\big\}\,.\notag
\end{align}
Selecting $\alpha=\beta/2$ and $\beta$ being sufficiently large,  identical to the arguments used in Section \ref{se:pftm1} for bounding the convergence rate of $A_2$, we have
$
{\rm I}\lesssim  o(\aleph^4)$
provided that  $K(\log n_*)^{1+2/\gamma}=o(n_*)$, $K\gtrsim L_n$, $K^{-\varphi}L_n^{\varphi}\log\{n_*(K\log p)^{-1}\}=o(1)$, $\log p=o[\min\{(n_*L_n^{-2}K)^{\varphi/(3\varphi+2)},(n_*K^{-1})^{\chi}\}]$ with $\chi=\min\{\gamma/(\gamma+4),1/3\}$. Also, by Cauchy-Schwarz inequality, it holds that
\begin{align}\label{eq:IIthm3}
{\rm II}
=&~p\sum_{i,j=1}^p\mathbb{E}\big\{|\hat{\sigma}_{u,i,j}-\sigma_{u,i,j}|^2{I}(|\hat{\sigma}_{u,i,j}-\sigma_{u,i,j}|>4\alpha\aleph){I}(|\sigma_{u,i,j}|>\alpha\aleph){I}(|\hat{\sigma}_{u,i,j}|\geq \beta\aleph)\}\notag\\
&+p\sum_{i,j=1}^p\mathbb{E}\big\{|\hat{\sigma}_{u,i,j}-\sigma_{u,i,j}|^2{I}(|\hat{\sigma}_{u,i,j}-\sigma_{u,i,j}|>4|\sigma_{u,i,j}|){I}(|\sigma_{u,i,j}|\leq\alpha\aleph){I}(|\hat{\sigma}_{u,i,j}|\geq \beta\aleph)\big\}\notag\\
\leq&~p\sum_{i,j=1}^p\big\{\mathbb{E}(|\hat{\sigma}_{u,i,j}-\sigma_{u,i,j}|^4)\big\}^{1/2}\big\{\mathbb{P}(|\hat{\sigma}_{u,i,j}-\sigma_{u,i,j}|>4\alpha\aleph)\big\}^{1/2}\\
&+p\sum_{i,j=1}^p\big\{\mathbb{E}(|\hat{\sigma}_{u,i,j}-\sigma_{u,i,j}|^4)\big\}^{1/2}\big[\mathbb{P}\{|\hat{\sigma}_{u,i,j}-\sigma_{u,i,j}|>(\beta-\alpha)\aleph\}\big]^{1/2}\,.\notag
\end{align}
Notice that
$
\max_{i,j\in[p]}\mathbb{E}(|\hat{\sigma}_{u,i,j}-\sigma_{u,i,j}|^4)\lesssim1$.
Since $\alpha=\beta/2$, repeating the arguments used in Section \ref{se:pftm1} for bounding the convergence rate of $A_2$ again, we have
\begin{align*}
p^3\bigg\{\max_{i,j\in[p]}\mathbb{P}(|\hat{\sigma}_{u,i,j}-\sigma_{u,i,j}|>4\alpha\aleph)\bigg\}^{1/2}\leq p^3\bigg[\max_{i,j\in[p]}\mathbb{P}\{|\hat{\sigma}_{u,i,j}-\sigma_{u,i,j}|>(\beta-\alpha)\aleph\}
\bigg]^{1/2}\lesssim o(\aleph^2)
\end{align*}
 provided that  $K(\log n_*)^{1+2/\gamma}=o(n_*)$, $K\gtrsim L_n$, $K^{-\varphi}L_n^{\varphi}\log\{n_*(K\log p)^{-1}\}=o(1)$ and  $\log p=o[\min\{(n_*L_n^{-2}K)^{\varphi/(3\varphi+2)},(n_*K^{-1})^{\chi}\}]$.
With sufficiently large $\beta$, we have
$
{\rm II}\lesssim o(\aleph^2)$. Together with ${\rm I}\lesssim o(\aleph^4)$, we have
$
\mathbb{E}\{(\max_{i\in[p]}\sum_{j=1}^p|d_{i,j}|)^2\}\leq{\rm I}+{\rm II}\lesssim o(\aleph
^2)$.
It follows from (\ref{eq:speupp1}) that
$
\sup_{\mathcal{P}_2}\mathbb{E}(\|\widehat{\bf\Sigma}_{u}^{\rm thre}-{\bf\Sigma}_{u}\|_2^2)\lesssim c_p^2(n_*^{-1}K\log p)^{1-q}
$. Analogously, in the case with $\varphi=\infty$, we have $\sup_{\mathcal{P}_2}\mathbb{E}(\|\widehat{\bf\Sigma}_{u}^{\rm thre}-{\bf\Sigma}_{u}\|_2^2)\lesssim c_p^2(n_*^{-1}K\log p)^{1-q}$ provided that $K(\log n_*)^{1+2/\gamma}=o(n_*)$, $K>L_n$ and  $\log p=o\{(n_*K^{-1})^{\chi}\}$.
 We complete the proof of Theorem \ref{tm:3}. $\hfill\Box$

\subsection{Proof of Theorem 4}

Same as the proof of Theorem 2, we also select $\mu_{i,t}=0$ and $\sigma_{i,t}=0$ for any $t\in[0,T]$. Then the associated $\bX_t={\bf 0}$ for any $t\in[0,T]$. In this special case, $\bY_{t_k}=\bU_{t_k}$. Given $(n,n_*)$ with $n\geq n_*$, and $0\leq t_1<\cdots<t_n=T$, we define $\mathcal{G}_*=\{\tilde{t}_1,\ldots,\tilde{t}_{n_*}\}$ with each $\tilde{t}_j\in\{t_1,\ldots,t_n\}$ and $\tilde{t}_j<\tilde{t}_{j+1}$. For each $t_j\in\mathcal{G}_*$, we assume all $p$ component processes are observed. For any $t_j\notin\mathcal{G}_*$, we assume only one component process are observed. Without loss of generality, we assume $\mathcal{G}_*=\{t_1,\ldots,t_{n_*}\}$. Let $n-n_*=ap+q_*$ where $a\geq 0$ and $0\leq q_*<p$ are two integers. We assume the $i$th component process is observed at $t_{n_*+jp+i}$'s with $j=0,\ldots,a$ and $i\in[p]$.
Then
$
\mathcal{G}_i=\mathcal{G}_*\cup\{t_{n_*+i},\ldots,t_{n_*+ap+i}\}
$.
The data we have is $\mathcal{Z}=\{\bY_{t_1,\mathcal{C}_1},\ldots,\bY_{t_n,\mathcal{C}_n}\}$ where $\mathcal{C}_k=\{i\in [p]:t_k\in\mathcal{G}_i\}$.

Let $r=\lfloor p/2\rfloor$, where $\lfloor x\rfloor$ denotes the largest integer less than or equal to $x$. Let $\mathcal{B}$ be the collection of all $p$-dimensional row vectors $v=(v_1,\ldots,v_p)$ such that $v_j=0$ for $1\leq j\leq p-r$ and $v_j=0$ or $1$ for $p-r+1\leq j\leq p$ under the restriction $\sum_{j=1}^p|v_j|=K_*$. We will specify $K_*$ later. If each $\lambda_j\in\mathcal{B}$, we say $\lambda=(\lambda_1,\ldots,\lambda_r)\in\mathcal{B}^r$. Set $\Gamma=\{0,1\}^r$ and $\Lambda\subset\mathcal{B}^r$. For each $\lambda=(\lambda_1,\ldots,\lambda_r)\in\Lambda$, we define $p\times p$ symmetric matrices $\bA_1(\lambda_1),\ldots,\bA_r(\lambda_r)$ where $\bA_m(\lambda_m)$ is a matrix with the $m$th row and $m$th column being $\lambda_m$ and $\lambda_m^{\top}$, respectively, and the rest of the entries being $0$. Define $\Theta=\Gamma\otimes\Lambda$. For each $\theta\in\Theta$, we write $\theta=\{\gamma(\theta),\lambda(\theta)\}$ with $\gamma(\theta)=\{\gamma_1(\theta),\ldots,\gamma_r(\theta)\}\in\Gamma$ and $\lambda(\theta)=\{\lambda_1(\theta),\ldots,\lambda_r(\theta)\}\in\Lambda$. We select $K_*=\lfloor c_p(n_*/\log p)^{q/2}\rfloor$ and define a collection $\mathcal{M}(\alpha,\nu)$ of covariance matrices as
\[
\mathcal{M}(\alpha,\nu)=\bigg\{\bSigma(\theta):\bSigma(\theta)=\alpha \bI_p+\sqrt{\frac{\nu\log p}{n_*}}\sum_{m=1}^r\gamma_m(\theta)\bA_m\{\lambda_m(\theta)\},~\theta\in\Theta\bigg\}\,,
\]
where $\alpha>0$ and $\nu>0$ are two constants. Notice that each $\bSigma\in\mathcal{M}(\alpha,\nu)$ has value $\alpha$ along the main diagonal, and contains an $r\times r$ submatrix, say $A$, at the upper right corner, $A^{\top}$ at the lower left corner and zero elsewhere. Write $\bSigma(\theta)=\{\sigma_{i,j}(\theta)\}_{p\times p}$. It holds that
$
\max_{\theta\in\Theta}\max_{i\in[p]}\sigma_{i,i}(\theta)=\alpha$ and $\max_{\theta\in\Theta}\max_{i\in[p]}\sum_{j=1}^p|\sigma_{i,j}(\theta)|^q\leq \alpha^q+c_p\nu^{q/2}$.
For sufficiently small $\alpha$ and $\nu$, we have $\mathcal{M}(\alpha,\nu)\subset\mathcal{H}(q,c_p,M)$ for $\mathcal{H}(q,c_p,M)$ defined as \eqref{eq:sparse-set}. Without loss of generality, we assume $\alpha=1$ in the sequel and write $\mathcal{M}(1,\nu)$ as $\mathcal{M}$ for simplification.

Let $\bU_{t_k}\sim N\{\bzero,\bSigma(\theta)\}$ with $\bSigma(\theta)\in\mathcal{M}$. When $\bU_{t_k}\sim N\{\bzero,\bSigma(\theta)\}$, we write the distribution of $\mathcal{Z}$ as $\mathbb{P}_\theta$. More specifically, the joint density of $\mathcal{Z}$ is
\begin{align*}
f_\theta=&~\prod_{k=1}^{n_*}\frac{1}{(2\pi)^{p/2}|\bSigma(\theta)|^{1/2}}\exp\bigg\{-\frac{1}{2}u_k^{\top }\bSigma^{-1}(\theta)u_k\bigg\}\times\prod_{k=n_*+1}^n\frac{1}{\sqrt{2\pi\sigma_{k,k}(\theta)}}\exp\bigg\{-\frac{u_{k,k}^2}{2\sigma_{k,k}(\theta)}\bigg\}\\
=&~\prod_{k=1}^{n_*}\frac{1}{(2\pi)^{p/2}|\bSigma(\theta)|^{1/2}}\exp\bigg\{-\frac{1}{2}u_k^{\top}\bSigma^{-1}(\theta)u_k\bigg\}\times\prod_{k=n_*+1}^n\frac{1}{\sqrt{2\pi}}\exp\bigg(-\frac{u_{k,k}^2}{2}\bigg)
\end{align*}
where $u_k=(u_{k,1},\ldots,u_{k,p})^{\top}$. It follows from Lemma 3 of \cite{Cai2012} with $s=2$ and $d$ being the matrix spectral norm $\|\cdot\|_2$ that
\begin{align*}
&\inf_{\widehat{\bSigma}}\max_{\theta\in\Theta}\mathbb{E}_{\mathcal{Z}|\theta}\big\{\|\widehat{\bSigma}-\bSigma(\theta)\|_2^2\big\}\geq \min_{(\theta,\theta'):H\{\gamma(\theta),\gamma(\theta')\}\geq 1}\frac{\|\bSigma(\theta)-\bSigma(\theta')\|_2^2}{H\{\gamma(\theta),\gamma(\theta')\}}\cdot\frac{r}{8}\cdot\min_{i\in[r]}\|\bar{\mathbb{P}}_{i,0}\wedge\bar{\mathbb{P}}_{i,1}\|\,,
\end{align*}
where $H(\cdot,\cdot)$ is the Hamming distance, and
$
\bar{\mathbb{P}}_{i,a}={2^{-(r-1)}|\Lambda|^{-1}}\sum_{\theta\in\{\theta\in\Theta:\gamma_i(\theta)=a\}}\mathbb{P}_\theta
$
for each $a\in\{0,1\}$. In the sequel, we will show the following two results:
\begin{equation}\label{eq:toprove1}
 \min_{(\theta,\theta'):H\{\gamma(\theta),\gamma(\theta')\}\geq 1}\frac{\|\bSigma(\theta)-\bSigma(\theta')\|_2^2}{H\{\gamma(\theta),\gamma(\theta')\}}\gtrsim \frac{c_p^2}{p}\bigg(\frac{\log p}{n_*}\bigg)^{1-q}
\end{equation}
and
 \begin{equation}\label{eq:toprove2}
\min_{i\in[r]}\|\bar{\mathbb{P}}_{i,0} \wedge \bar{\mathbb{P}}_{i,1}\|\gtrsim1\,.
 \end{equation}
Recall $r=\lfloor p/2\rfloor$. Then we will have Theorem \ref{tm:4}. The proofs of (\ref{eq:toprove1}) and (\ref{eq:toprove2}) are identical to that for Lemmas 5 and 6 in \cite{Cai2012}, respectively. Hence, we omit here. $\hfill\Box$

\subsection{Proof of Theorem \ref{tm:5}}

As shown in (\ref{eq:IIij}), ${\rm II}(i,j)={\rm II}_1(i,j)+{\rm II}_2(i,j)+{\rm II}_3(i,j)+{\rm II}_4(i,j)$. Notice that $\Delta_K$ is a fixed integer. By Jensen's inequality,
$
\mathbb{E}\{(\int_{t_{i,j,k}}^{t_{i,j,\ell}}\mu_{i,s}\,{\rm d}s)^2\}\leq|t_{i,j,\ell}-t_{i,j,k}|\int_{t_{i,j,\ell}\wedge t_{i,j,k}}^{t_{i,j,\ell}\vee t_{i,j,k}}\mathbb{E}(\mu_{i,s}^2)\,{\rm d}s\leq|t_{i,j,\ell}-t_{i,j,k}|^2\max_{0\leq s\leq T}\mathbb{E}(\mu_{i,s}^2)$,
which implies
$
\max_{i,j\in[p]}\max_{k\in[n_{i,j}]}\max_{\ell:\,K\leq|\ell-k|\leq K+\Delta_K}\mathbb{E}\{(\int_{t_{i,j,k}}^{t_{i,j,\ell}}\mu_{i,s}\,{\rm d}s)^2\}\lesssim K^2n_*^{-2}$.
Due to
$
\mathbb{E}\{(\int_{t_{i,j,k}}^{t_{i,j,\ell}}\sigma_{j,s}\,{\rm d}B_{j,s})^2\}=\mathbb{E}(\int_{t_{i,j,\ell}\wedge t_{i,j,k}}^{t_{i,j,\ell}\vee t_{i,j,k}}\sigma_{j,s}^2\,{\rm d}s)\leq|t_{i,j,\ell}-t_{i,j,k}|\max_{0\leq s\leq T}\mathbb{E}(\sigma_{j,s}^2)$,
we have
$
\max_{i,j\in[p]}\max_{k\in[n_{i,j}]}\max_{\ell:\,K\leq|\ell-k|\leq K+\Delta_K}\mathbb{E}\{(\int_{t_{i,j,k}}^{t_{i,j,\ell}}\sigma_{j,s}\,{\rm d}B_{j,s})^2\}\lesssim Kn_*^{-1}$.
By Cauchy-Schwarz inequality, we have that $
\max_{i,j\in[p]}|\mathbb{E}\{{\rm II}_1(i,j)\}|\lesssim K^2n_*^{-2}$, $\max_{i,j\in[p]}|\mathbb{E}\{{\rm II}_3(i,j)\}|\lesssim K^{3/2}n_*^{-3/2}$ and $\max_{i,j\in[p]}|\mathbb{E}\{{\rm II}_4(i,j)\}|\lesssim K^{3/2}n_*^{-3/2}$.
Notice that
\begin{align*}
\mathbb{E}\{{\rm II}_2(i,j)\}=&~\frac{1}{2n_{i,j}}\sum_{k=1}^{n_{i,j}}\frac{1}{N_{i,j,k}}\sum_{K\leq |\ell-k|\leq K+\Delta_K}\mathbb{E}\bigg(\int_{t_{i,j,k}\wedge t_{i,j,\ell}}^{t_{i,j,k}\vee t_{i,j,\ell}}\sigma_{i,s}\sigma_{j,s}\rho_{i,j,s}\,{\rm d}s\bigg)\\
=&~\frac{1}{2n_{i,j}}\sum_{\ell=1}^{n_{i,j}-1}\bigg[\sum_{k=\ell-K-\Delta_K+1}^{\ell-K}\frac{\min\{k+K+\Delta_K,n_{i,j}\}-\ell}{N_{i,j,k}}\\
&~~~~~~~~~~~~~~~~~+\sum_{k=\ell-K+1}^{\ell}\frac{|\min\{k+K+\Delta_K,n_{i,j}\}-k-K+1|_+}{N_{i,j,k}}\\
&~~~~~~~~~~~~~~~~~+\sum_{k=\ell+1}^{\ell+K}\frac{|k-K+1-\max\{k-K-\Delta_K,1\}|_+}{N_{i,j,k}}\\
&~~~~~~~~~~~~~~~~~+\sum_{k=\ell+K+1}^{\ell+K+\Delta_K}\frac{\ell+1-\max\{k-K-\Delta_K,1\}}{N_{i,j,k}}\bigg]\mathbb{E}\bigg(\int_{t_{i,j,\ell}}^{t_{i,j,\ell+1}}\sigma_{i,s}\sigma_{j,s}\rho_{i,j,s}\,{\rm d}s\bigg)\\
:=&~\frac{1}{2n_{i,j}}\sum_{\ell=1}^{n_{i,j}-1} Q_{i,j,\ell}\cdot\mathbb{E}\bigg(\int_{t_{i,j,\ell}}^{t_{i,j,\ell+1}}\sigma_{i,s}\sigma_{j,s}\rho_{i,j,s}\,{\rm d}s\bigg)\,,
\end{align*}
where we adopt the convention $N_{i,j,k}=\infty$ if $k>n_{i,j}$ or $k<0$.
For sufficiently large $n$, $N_{i,j,k}$ follows the formula \eqref{eq:Nijk}. Since $K=o(n_*)$ and $\Delta_K$ is a fixed integer, for sufficiently large $n$, we have
\begin{align*}
Q_{i,j,\ell}\,\left\{ \begin{aligned}
~\lesssim K\,,~~~~~~~&\textrm{if}~1\leq \ell\leq 2(K+\Delta_K)-1 \,, \\
=\frac{2K+\Delta_K}{2}\,,~~~&\textrm{if}~2(K+\Delta_K)\leq \ell\leq n_{i,j}-2(K+\Delta_K) \,, \\
\lesssim K\,,~~~~~~~&\textrm{if}~n_{i,j}-2(K+\Delta_K)+1\leq \ell\leq n_{i,j}-1\,,
\end{aligned} \right.
\end{align*}
which implies that
\[
\begin{split}
\mathbb{E}\{{\rm II}_2(i,j)\}=&~\frac{2K+\Delta_K}{4n_{i,j}}\sum_{\ell=2(K+\Delta_K)}^{n_{i,j}-2(K+\Delta_K)}\mathbb{E}\bigg(\int_{t_{i,j,\ell}}^{t_{i,j,\ell+1}}\sigma_{i,s}\sigma_{j,s}\rho_{i,j,s}\,{\rm d}s\bigg)+R_1(i,j)\\
=&~\frac{2K+\Delta_K}{4n_{i,j}}\mathbb{E}\bigg(\int_{t_{i,j,1}}^{t_{i,j,n_{i,j}}}\sigma_{i,s}\sigma_{j,s}\rho_{i,j,s}\,{\rm d}s\bigg)+R_2(i,j)\,,
\end{split}
\]
where $\max_{i,j\in[p]}|R_1(i,j)|=O(K^2n_*^{-2})$ and $\max_{i,j\in[p]}|R_2(i,j)|=O(K^2n_*^{-2})$. Therefore, it follows from \eqref{eq:IIij} that
\[
\begin{split}
\max_{i,j\in[p]}\bigg|\mathbb{E}\bigg\{{\rm II}(i,j)-\frac{2K+\Delta_K}{4n_{i,j}}\int_{t_{i,j,1}}^{t_{i,j,n_{i,j}}}\sigma_{i,s}\sigma_{j,s}\rho_{i,j,s}\,{\rm d}s\bigg\}\bigg|
\lesssim&~\frac{K^{3/2}}{n_*^{3/2}}\,.
\end{split}
\]
We complete the proof of Theorem \ref{tm:5}. $\hfill\Box$

\subsection{Proof of Theorem \ref{corr1}}

Write $\aleph=(n_*^{-1}{K\log p})^{1/2}$. For Part (i),  Theorem \ref{pn:1} implies $\max_{i,j\in[p]}|\hat{\sigma}_{u,i,j}-\sigma_{u,i,j}|=O_\p(\aleph)$. Due to  $\max_{i,j\in[p]}|\hat{\psi}_{i,j}|=O_\p(\log p)$ and $\aleph=o(1)$, by Triangle inequality,
$
\max_{i,j\in[p]}|\hat{\sigma}_{u,i,j}^{{\rm bc}}-\sigma_{u,i,j}|\leq  \max_{i,j\in[p]}|\hat{\sigma}_{u,i,j}-\sigma_{u,i,j}|+\max_{i,j\in[p]}(4n_{i,j})^{-1}(2K+\Delta_K)|\hat{\psi}_{i,j}|=O_\p(\aleph)$. Then Part (i) holds. To prove Part (ii), we define an event $\mathcal{E}=\{\max_{i,j\in[p]}(4n_{i,j})^{-1}(2K+\Delta_K)|\hat{\psi}_{i,j}|\leq \epsilon\aleph\}$ for some constant  $\epsilon>0$. Since  $\max_{i,j\in[p]}|\hat{\psi}_{i,j}|=O_\p(\log p)$ and $\aleph=o(1)$, then $\max_{i,j\in[p]}(4n_{i,j})^{-1}(2K+\Delta_K)|\hat{\psi}_{i,j}|=O_{\p}(\aleph^2)$, which implies $\mathbb{P}(\mathcal{E}^c)=o(1)$.
For any $C>0$, by Markov's inequality, we have
$
\mathbb{P}(\|\widehat{\bSigma}_u^{{\rm bc,thre}}-\bSigma_u\|_2\geq C c_p\aleph^{1-q})\leq \mathbb{P}(\|\widehat{\bSigma}_u^{{\rm bc,thre}}-\bSigma_u\|_2\geq C c_p\aleph^{1-q}\,,\mathcal{E})+\mathbb{P}(\mathcal{E}^c)\leq
C^{-2}c_p^{-2}\aleph^{2(q-1)}\mathbb{E}\{\|\widehat{\bSigma}_u^{{\rm bc,thre}}-\bSigma_u\|_2^2I(\mathcal{E})\}+\mathbb{P}(\mathcal{E}^c)$.
In the sequel, we will show $\mathbb{E}\{\|\widehat{\bSigma}_u^{{\rm bc,thre}}-\bSigma_u\|_2^2I(\mathcal{E})\}\lesssim c_p^2\aleph^{2(1-q)}$. Based on this result, we know Part (ii) holds.

 For each $i,j\in[p]$, we  define the event
$
A_{i,j}=\{|\hat{\sigma}_{u,i,j}^{\rm bc, thre}-\sigma_{u,i,j}|\leq 4\min(|\sigma_{u,i,j}|,\alpha\aleph)\}
$
with some $\alpha>0$, and
$
d_{i,j}=(\hat{\sigma}_{u,i,j}^{\rm bc, thre}-\sigma_{u,i,j}){I}(A_{i,j}^c)$.
Write $\bD=(d_{i,j})_{p\times p}$. Identical to \eqref{eq:speupp} and the arguments below it, $
\|\widehat{\bSigma}_{u}^{\rm bc, thre}-\bSigma_{u}\|_2^2\lesssim(\max_{i\in[p]}\sum_{j=1}^p|d_{i,j}|)^2+\{\max_{i\in[p]}\sum_{j=1}^p|\hat{\sigma}_{u,i,j}^{\rm bc, thre}-\sigma_{u,i,j}|{I}(A_{i,j})\}^2\lesssim(\max_{i\in[p]}\sum_{j=1}^p|d_{i,j}|)^2+c_p^2\aleph^{2(1-q)}$,
which implies
$
\mathbb{E}\{\|\widehat{\bSigma}_{u}^{\rm bc, thre}-\bSigma_{u}\|_2^2I(\mathcal{E})\}\lesssim \mathbb{E}\{(\max_{i\in[p]}\sum_{j=1}^p|d_{i,j}|)^2I(\mathcal{E})\}+c_p^2\aleph^{2(1-q)}$.
Identical to \eqref{eq:Edij}, we have
\begin{align*}
\mathbb{E}\bigg\{\bigg(\max_{i\in[p]}\sum_{j=1}^p|d_{i,j}|\bigg)^2I(\mathcal{E})\bigg\}\leq&~\underbrace{p\sum_{i,j=1}^p\mathbb{E}\big(|\hat{\sigma}_{u,i,j}^{\rm bc, thre}-\sigma_{u,i,j}|^2I(\mathcal{E}){I}[A_{i,j}^c\cap\{\hat{\sigma}_{u,i,j}^{\rm bc, thre}=0\}]\big)}_{\textrm{I}}\\
&+\underbrace{p\sum_{i,j=1}^p\mathbb{E}\big(|\hat{\sigma}_{u,i,j}^{\rm bc, thre}-\sigma_{u,i,j}|^2I(\mathcal{E}){I}[A_{i,j}^c\cap\{\hat{\sigma}_{u,i,j}^{\rm bc, thre}=\hat{\sigma}_{u,i,j}^{\rm bc}\}]\big)}_{\textrm{II}}\,.
\end{align*}
Recall $\hat{\sigma}_{u,i,j}^{\rm bc, thre}=\hat{\sigma}_{u,i,j}^{\rm bc}{I}(|\hat{\sigma}_{u,i,j}^{\rm bc}|\geq \beta\aleph)$ and $\max_{i,j\in[p]}|\hat{\sigma}_{u,i,j}^{\rm bc}-\hat{\sigma}_{u,i,j}|\leq\epsilon\aleph$ restricted on $\mathcal{E}$. Then
$
{\rm I}=p\sum_{i,j=1}^p\sigma_{u,i,j}^2\mathbb{P}(|\sigma_{u,i,j}|\geq4\alpha\aleph\,, |\hat{\sigma}_{u,i,j}^{\rm bc}|< {\beta}\aleph\,, \mathcal{E})
\leq p\sum_{i,j=1}^p\sigma_{u,i,j}^2\mathbb{P}\{|\hat{\sigma}_{u,i,j}-\sigma_{u,i,j}|\geq(4\alpha-{\beta}-\epsilon)\aleph\}$.
Notice that
$
\max_{i,j\in[p]}\mathbb{E}\{|\hat{\sigma}_{u,i,j}^{\rm bc}-\sigma_{u,i,j}|^4I(\mathcal{E})\}\lesssim \max_{i,j\in[p]}\mathbb{E}(|\hat{\sigma}_{u,i,j}-\sigma_{u,i,j}|^4)+\aleph^4
\lesssim1$. Identical to \eqref{eq:IIthm3}, we have
$
{\rm II}\lesssim p\sum_{i,j=1}^p\{\mathbb{P}(|\hat{\sigma}_{u,i,j}^{\rm bc}-\sigma_{u,i,j}|>4\alpha\aleph\,,\mathcal{E})\}^{1/2}+p\sum_{i,j=1}^p[\mathbb{P}\{|\hat{\sigma}_{u,i,j}^{\rm bc}-\sigma_{u,i,j}|>({\beta}-\alpha)\aleph\,, \mathcal{E}\}]^{1/2}\leq  p\sum_{i,j=1}^p[\mathbb{P}\{|\hat{\sigma}_{u,i,j}-\sigma_{u,i,j}|>(4\alpha-\epsilon)\aleph\,,\mathcal{E})\}]^{1/2}+p\sum_{i,j=1}^p[\mathbb{P}\{|\hat{\sigma}_{u,i,j}-\sigma_{u,i,j}|>({\beta}-\alpha-\epsilon)\aleph\,, \mathcal{E}\}]^{1/2}$. Selecting $\alpha=2\epsilon={\beta}/2$ for some sufficiently large $\beta>0$, applying the same arguments for bounding ${\rm I}$ and ${\rm II}$ there in Section \ref{sec:pfthm3}, we have
$
{\rm I}+{\rm II}\lesssim o(\aleph
^2)$, which implies
$
\mathbb{E}\{\|\widehat{\bf\Sigma}_{u}^{\rm bc, thre}-{\bf\Sigma}_{u}\|_2^2I(\mathcal{E})\}\lesssim c_p^2\aleph^{2(1-q)}
$.
We complete the proof of Theorem \ref{corr1}. $\hfill\Box$

\subsection{Proof of Part (i) in Theorem \ref{eq:the6}}\label{sec:prothe6}
For any $k\in[n_{i,j}]$, let $S_{i,j,k}=\{t_{i, j, \ell}:K\leq|\ell-k|\leq K+\Delta_K \}$. For any $i,j\in[p]$, we have that
\begin{align*}
\hat{\sigma}_{u,i,j}-\sigma_{u,i,j}=&~\underbrace{\frac{1}{2n_{i,j}}\sum_{k=1}^{n_{i,j}}\frac{1}{N_{i,j,k}}\sum_{t_{i,j,\ell}\in S_{i,j,k}}(U_{i,t_{i,j,\ell}}-U_{i,t_{i,j,k}})(U_{j,t_{i,j,\ell}}-U_{j,t_{i,j,k}})-\sigma_{u,i,j}}_{\textrm{I}(i,j)}\\
&+\underbrace{\frac{1}{2n_{i,j}}\sum_{k=1}^{n_{i,j}}\frac{1}{N_{i,j,k}}\sum_{t_{i,j,\ell}\in S_{i,j,k}}(X_{i,t_{i,j,\ell}}-X_{i,t_{i,j,k}})(X_{j,t_{i,j,\ell}}-X_{j,t_{i,j,k}})}_{\textrm{II}'(i,j)}\\
&+\underbrace{\frac{1}{2n_{i,j}}\sum_{k=1}^{n_{i,j}}\frac{1}{N_{i,j,k}}\sum_{t_{i,j,\ell}\in S_{i,j,k}}(X_{i,t_{i,j,\ell}}-X_{i,t_{i,j,k}})(U_{j,t_{i,j,\ell}}-U_{j,t_{i,j,k}})}_{\textrm{III}'(i,j)}\\
&+\underbrace{\frac{1}{2n_{i,j}}\sum_{k=1}^{n_{i,j}}\frac{1}{N_{i,j,k}}\sum_{t_{i,j,\ell}\in S_{i,j,k}}(U_{i,t_{i,j,\ell}}-U_{i,t_{i,j,k}})(X_{j,t_{i,j,\ell}}-X_{j,t_{i,j,k}})}_{\textrm{IV}'(i,j)}\,.
\end{align*}
Based on Proposition \ref{laa:1}, we can obtain the convergence rate of $\max_{ i,j\in[p]}|{\rm I}({i,j})|$.
Notice that ${\rm d}X_{i,t}=\mu_{i,t}\,{\rm d}t+\sigma_{i,t}\,{\rm d}B_{i,t} + J_{i,t}\,{\rm d}M_{i,t}$. In comparison to ${\rm II}({i,j})$ specified in \eqref{eq:IIij}, we have
\begin{align}\label{eq:IIij2}
&{\rm II}'(i,j)-\frac{1}{2n_{i,j}}\sum_{k=1}^{n_{i,j}}\frac{1}{N_{i,j,k}}\sum_{t_{i,j,\ell} \in S_{i,j,k}}\bigg(\int_{t_{i,j,k}}^{t_{i,j,\ell}}J_{i,s}\,{\rm d}M_{i,s}\bigg)\bigg(\int_{t_{i,j,k}}^{t_{i,j,\ell}}J_{j,s}\,{\rm d}M_{j,s}\bigg)\notag\\
&~~~~~~~={\rm II}(i,j)+
\underbrace{\frac{1}{2n_{i,j}}\sum_{k=1}^{n_{i,j}}\frac{1}{N_{i,j,k}}\sum_{t_{i,j,\ell}\in S_{i,j,k}}\bigg(\int_{t_{i,j,k}}^{t_{i,j,\ell}}\mu_{i,s}\,{\rm d}s\bigg)\bigg(\int_{t_{i,j,k}}^{t_{i,j,\ell}}J_{j,s}\,{\rm d}M_{j,s}\bigg)}_{{\rm II}_5(i,j)}\notag\\
&~~~~~~~~~~+\underbrace{\frac{1}{2n_{i,j}}\sum_{k=1}^{n_{i,j}}\frac{1}{N_{i,j,k}}\sum_{t_{i,j,\ell} \in S_{i,j,k}}\bigg(\int_{t_{i,j,k}}^{t_{i,j,\ell}}\sigma_{i,s}\,{\rm d}B_{i,s}\bigg)\bigg(\int_{t_{i,j,k}}^{t_{i,j,\ell}}J_{j,s}\,{\rm d}M_{j,s}\bigg)}_{{\rm II}_6(i,j)}\\
&~~~~~~~~~~+\underbrace{\frac{1}{2n_{i,j}}\sum_{k=1}^{n_{i,j}}\frac{1}{N_{i,j,k}}\sum_{t_{i,j,\ell} \in S_{i,j,k}}\bigg(\int_{t_{i,j,k}}^{t_{i,j,\ell}}J_{i,s}\,{\rm d}M_{i,s}\bigg)\bigg(\int_{t_{i,j,k}}^{t_{i,j,\ell}}\sigma_{j,s}\,{\rm d}B_{j,s}\bigg)}_{{\rm II}_7(i,j)}\notag\\
&~~~~~~~~~~+\underbrace{\frac{1}{2n_{i,j}}\sum_{k=1}^{n_{i,j}}\frac{1}{N_{i,j,k}}\sum_{t_{i,j,\ell} \in S_{i,j,k}}\bigg(\int_{t_{i,j,k}}^{t_{i,j,\ell}}J_{i,s}\,{\rm d}M_{i,s}\bigg)\bigg(\int_{t_{i,j,k}}^{t_{i,j,\ell}}\mu_{j,s}\,{\rm d}s\bigg)}_{{\rm II}_8(i,j)}\notag\,.
\end{align}
For ${\rm III}(i,j)$ specified in Section \ref{se:pftm1}, we have ${\rm III}'(i,j)={\rm III}(i,j)
+{(2n_{i,j})^{-1}}\sum_{k=1}^{n_{i,j}}\sum_{t_{i,j,\ell}\in S_{i,j,k}}({N_{i,j,\ell}^{-1}}+N_{i,j,k}^{-1})(\int_{t_{i,j,\ell}}^{t_{i,j,k}}J_{i,s}\,{\rm d}M_{i,s})U_{j,t_{i,j,k}}$.
Propositions  \ref{la:2} and \ref{la:3} give the convergence rates of $\max_{ i,j\in[p]}|{\rm II}(i,j)|$ and $\max_{ i,j\in[p]}|{\rm III}(i,j)|$, respectively.
Write ${\rm II}{''}(i,j)={\rm II}'(i,j)-{\rm II}(i,j)-\varpi_{i,j}$ and
$
{\rm III}''(i,j)={\rm III}'(i,j)-{\rm III}(i,j)
$.
 To prove Theorem \ref{eq:the6}, we need the following two propositions whose proofs are given in Sections \ref{sec:pro4} and \ref{sec:pro5}, respectively.
\begin{proposition}\label{eq:proj4}
Under Assumptions $\ref{as:space}$ and $\ref{as:drifdiff}$--$\ref{as:jumpsize}$, if  $K (\log n_*)^{(2\iota+2\gamma+\iota\gamma)/(\iota\gamma)}=o(n_*)$, we have that
\begin{align*}
\max_{i,j\in[p]}\mathbb{P}\{|{\rm II}''(i,j)|>v\}\lesssim \exp(-Cn_*K^{-1}v) +
\exp\{-C(n_*K^{-1})^{\iota\gamma/(2\iota+2\gamma+\iota\gamma)}\}
\end{align*}
for any $v\gg n_*^{-1}K\lambda_*$, where $\gamma$ and $\iota$ are specified in Assumptions $\ref{as:diff}$ and $\ref{as:jumpsize}$, respectively.
Furthermore, it holds that $\max_{ i,j\in[p]}\mathbb{E}\{|{\rm II}''(i,j)|^{m}\}\lesssim n_*^{(1-m)/2}K^{m/2}$
for any fixed positive integer $m$ provided that $n_*^{-1}K\lambda_*=o(1)$.
\end{proposition}

\begin{proposition}\label{eq:pfpro5}
Under Assumptions $\ref{as:space}$--$\ref{as:jumpsize}$, if $K(\log n_*)^{2+2/\iota}=o(n_*)$, we have that
\begin{align*}
\max_{i,j\in[p]}\mathbb{P}\{|{\rm{III}}''(i,j)|\geq v\}\lesssim\exp\{-C(n_*K^{-1})^{(3\iota+2)/(4\iota+4)}v\}+\exp\{-C(n_*K^{-1})^{\iota/(2\iota+2)}\}
\end{align*}
for any $v\gg (n_*^{-1}K)^{(3\iota+2)/(4\iota+4)}\lambda_*$, where  $\iota$ is specified in Assumption  $\ref{as:jumpsize}$.
Furthermore, it holds that $\max_{ i,j\in[p]}\mathbb{E}\{|{\rm III}''(i,j)|^{m}\}\lesssim n_*^{1/2}$
for any fixed positive integer $m$ provided that  $n_*^{-1}K\lambda_*=o(1)$.
\end{proposition}

Recall that $\hat{\sigma}_{u,i,j}^{{\rm jump}}=\hat{\sigma}_{u,i,j}-\varpi_{i,j}=\sigma_{u,i,j}+{\rm I}(i,j)+{\rm II}(i,j)+{\rm II}''(i,j)+{\rm III}'(i,j)+{\rm IV}'(i,j)$. Write $\aleph=(n_*^{-1}K\log p)^{1/2}$. We first consider the case with $\varphi<\infty$. Notice that $(1+cx^{-1})^{-x}\geq e^{-c}$ for any $x>0$ and $c>0$, and ${\rm III}'(i,j)={\rm III}(i,j)+{\rm III}''(i,j)$. Since the tail probability of $\max_{ i,j\in[p]}|{\rm{IV}}'(i,j)|$ is the same as that of $\max_{ i,j\in[p]}|{\rm{III}}'(i,j)|$,
by Propositions \ref{laa:1}--\ref{eq:pfpro5}, if $K\gtrsim L_n$ and $K (\log n_*)^{\chi_1}=o(n_*)$ with $\chi_1=\max\{(2\iota+2\gamma+\iota\gamma)/(\iota\gamma),2+2/\iota,1+2/\gamma\}$, we have
\begin{align}\label{eq:tailb12}
\notag\max_{i,j\in[p]}\mathbb{P}\big(|\hat{\sigma}_{u,i,j}^{{\rm jump}}-\sigma_{u,i,j}|>v\big)
\lesssim& ~\{1+n_*(K+L_n)^{-1}\rho^{-1}v^2\}^{-\rho/2}+v^{-1}\exp\{-C(n_{*}L_n^{-1}\rho^{-1}v)^{\varphi/(\varphi+1)}\}\\
&+v^{-1}\exp(-Cn_{*}K^{-1}\rho^{-1}v)+\exp\{-C(n_*L_n^{-1}v)^{\varphi/(2\varphi+1)}\}\notag\\&+\exp\{-C(n_*K^{-1})^{(3\iota+2)/(4\iota+4)}v\}+\exp\{-C(n_*K^{-1})^{\eta}\}
\end{align}
for any $\rho\geq 1$  and  $(n_*^{-1}K)^{2/(\gamma+4)} \gg v\gg \max\{\exp(-CL_n^{-\varphi}K^{\varphi}),(n_*^{-1}K)^{(3\iota+2)/(4\iota+4)}\lambda_*\}$, where $\eta=\min\{\gamma/(\gamma+4),\iota\gamma/(2\iota+2\gamma+\iota\gamma),\iota/(2\iota+2)\}$. Since $K^{-\varphi}L_n^{\varphi}\log\{n_*(K\log p)^{-1}\}=o(1)$, $\lambda_*^2(n_*^{-1}K)^{\iota/(2\iota+2)} (\log p)^{-1}=o(1)$ and $\log p=o\{(n_*K^{-1})^{\gamma/(\gamma+4)}\}$, then  $(n_*^{-1}K)^{2/(\gamma+4)} \gg \aleph\gg \max\{\exp(-CL_n^{-\varphi}K^{\varphi}),(n_*^{-1}K)^{(3\iota+2)/(4\iota+4)}\lambda_*\}$. Given a sufficiently large constant $\alpha>0$, we have
$
\mathbb{E}(|\widehat{\bSigma}_{u}^{{\rm jump}}-\bSigma_{u}|_\infty)
\leq\mathbb{E}\{\max_{i,j\in [p]}|\hat{\sigma}_{u,i,j}^{{\rm jump}}-\sigma_{u,i,j}|{I}(|\hat{\sigma}_{u,i,j}^{{\rm jump}}-\sigma_{u,i,j}|\leq \alpha\aleph)\}
+\mathbb{E}\{\max_{i,j\in[p]}|\hat{\sigma}_{u,i,j}^{{\rm jump}}-\sigma_{u,i,j}|{I}(|\hat{\sigma}_{u,i,j}^{{\rm jump}}-\sigma_{u,i,j}|> \alpha\aleph)\}
=:A_1^*+A_2^*$.
It is easy to see that $A_1^*\leq \alpha\aleph$. By Cauchy-Schwarz inequality, we have
$
A_2^*\leq\sum_{i,j=1}^p\mathbb{E}\{|\hat{\sigma}_{u,i,j}^{{\rm jump}}-\sigma_{u,i,j}|{I}(|\hat{\sigma}_{u,i,j}^{{\rm jump}}-\sigma_{u,i,j}|> \alpha\aleph)\}
\leq p^2\max_{i,j\in[p]}\{\mathbb{E}(|\hat{\sigma}_{u,i,j}^{{\rm jump}}-\sigma_{u,i,j}|^2)\}^{1/2}\cdot\max_{i,j\in[p]}\{\mathbb{P}(|\hat{\sigma}_{u,i,j}^{{\rm jump}}-\sigma_{u,i,j}|> \alpha\aleph)\}^{1/2}$.
Let $\rho\asymp\log p\geq 1$. Since $K\log p=o(n_*)$, (\ref{eq:tailb12}) implies
$
\max_{i,j\in[p]}\mathbb{P}(|\hat{\sigma}_{u,i,j}^{{\rm jump}}-\sigma_{u,i,j}|> \alpha\aleph)
\lesssim p^{-2w}+\exp[-C\{n_*KL_n^{-2}(\log p)^{-1}\}^{\varphi/(2\varphi+2)}]+\exp[-C\{n_*K^{-1}(\log p)^{-1}\}^{1/2}]
+\exp\{-C(n_*KL_n^{-2}\log p)^{\varphi/(4\varphi+2)}\}
+\exp\{-C(n_*K^{-1})^{\iota/(4\iota+4)}(\log p)^{1/2}\}+\exp\{-C(n_*K^{-1})^{\eta}\}$
with some sufficiently large $w>0$, where $w\rightarrow\infty$ as $\alpha\rightarrow\infty$.
Due to
$
\max_{i,j\in[p]}\mathbb{E}(|\hat{\sigma}_{u,i,j}^{{\rm jump}}-\sigma_{u,i,j}|^2)\lesssim n_*^{1/2}\leq  p^{2c_*}$ with $c_*=1/(4\kappa)$, where $\kappa$ is specified in the beginning of Section \ref{sec:theoryfU}, if $\log p=o[\min\{(n_*L_n^{-2}K)^{\varphi/(3\varphi+2)},(n_*K^{-1})^{\chi_2}\}]$ with $\chi_2=\min\{\iota\gamma/(2\iota+2\gamma+\iota\gamma),\iota/(2\iota+2),\gamma/(\gamma+4),1/3\}$, then
$
A_2^*\lesssim p^{2+c_*-w}+\exp[-C\{n_*KL_n^{-2}(\log p)^{-1}\}^{\varphi/(2\varphi+2)}]+\exp[-C\{n_*K^{-1}(\log p)^{-1}\}^{1/2}]+\exp\{-C(n_*KL_n^{-2}\log p)^{\varphi/(4\varphi+2)}\}+\exp\{-C(n_*K^{-1})^{\iota/(4\iota+4)}(\log p)^{1/2}\}+\exp\{-C(n_*K^{-1})^{\eta}\}
=o\{(n_*^{-1}K\log p)^4\}$.
Hence,
$$
\sup_{\mathcal{P}_3}\mathbb{E}(|\widehat{\bSigma}_{u}^{{\rm jump}}-\bSigma_{u}|_\infty)\lesssim(n_*^{-1}K\log p)^{1/2}$$ provided that $K^{-\varphi}L_n^{\varphi}\log\{n_*(K\log p)^{-1}\}=o(1)$, $\lambda_*^2(n_*^{-1}K)^{\iota/(2\iota+2)} (\log p)^{-1}=o(1)$, $K\gtrsim L_n$ and $K (\log n_*)^{\chi_1}=o(n_*)$ with $\chi_1=\max\{(2\iota+2\gamma+\iota\gamma)/(\iota\gamma),2+2/\iota,1+2/\gamma\}$, and  $\log p=o[\min\{(n_*L_n^{-2}K)^{\varphi/(3\varphi+2)},\\(n_*K^{-1})^{\chi_2}\}]$ with $\chi_2=\min\{\iota\gamma/(2\iota+2\gamma+\iota\gamma),\iota/(2\iota+2),\gamma/(\gamma+4),1/3\}$.

Now we consider the case with $\varphi=\infty$. As we discussed in Remark \ref{rek:1}(i), if $\{\bU_{t_k}\}$ is an independent sequence, we can select $L_n=1/2$. Due to $K\geq 1$, we have $K>L_n$ in this case. Without loss of generality, we can always assume $K>L_n$ when $\varphi=\infty$. Based on Remark \ref{rek:1}, it holds that $\{1+n_*(K+L_n)^{-1}\rho^{-1}v^2\}^{-\rho/2}+v^{-1}\exp(-Cn_{*}K^{-1}\rho^{-1}v)$ for any $v>0$ and $\rho\geq 1$ under either of the scenarios: (i) $\{\bU_{t_k}\}$ is an independent sequence, and (ii) $\{\bU_{t_k}\}$ is an $L_n$-dependent sequence. Repeating the arguments for $\varphi<\infty$, we have $
\sup_{\mathcal{P}_3}\mathbb{E}(|\widehat{\bSigma}_{u}^{{\rm jump}}-\bSigma_{u}|_\infty)\lesssim(n_*^{-1}K\log p)^{1/2}$ provided that  $\lambda_*^2(n_*^{-1}K)^{\iota/(2\iota+2)} (\log p)^{-1}=o(1)$, $K (\log n_*)^{\chi_1}=o(n_*)$ with $\chi_1=\max\{(2\iota+2\gamma+\iota\gamma)/(\iota\gamma),2+2/\iota,1+2/\gamma\}$ and $\log p=o\{(n_*K^{-1})^{\chi_2}\}$ with $\chi_2=\min\{\iota\gamma/(2\iota+2\gamma+\iota\gamma),\iota/(2\iota+2),\gamma/(\gamma+4),1/3\}$.   We complete the proof of part (i) in Theorem \ref{eq:the6}. $\hfill\Box$

%

\subsection{Proof of Proposition \ref{eq:proj4}}\label{sec:pro4}
Recall $\xi=\max_{i,j\in[p]}\max_{ k\in[n_{i,j}]}\max_{t_{i,j,\ell}\in S_{i,j,k}}|t_{i,j,\ell}-t_{i,j,k}|\asymp n_*^{-1}K$, $S_{i,j,k}=\{t_{i, j, \ell}:K\leq|\ell-k|\leq K+\Delta_K \}$  and $N_{i ,j , k}=|S_{i,j,k}|$.
In the sequel, we will bound the tail probabilities of
$\max_{i,j\in[p]}|\textrm{II}_5(i,j)|$,
$\max_{i,j\in[p]}|\textrm{II}_6(i,j)|$, $\max_{i,j\in[p]}|\textrm{II}_7(i,j)|$ and  $\max_{i,j\in[p]}|\textrm{II}_8(i,j)|$, respectively.

We  first bound the tail probabilities  of $\max_{i,j\in[p]}|\textrm{II}_5(i,j)|$ and $\max_{i,j\in[p]}|\textrm{II}_8(i,j)|$.
 Notice that for sufficiently large $n$, $\min_{k \in[n_{i,j}]}N_{i,j,k}=\Delta_K+1$ and $\max_{k\in[n_{i,j}]}N_{i,j,k}=2\Delta_K+2$. Since  $\Delta_K$ is a fixed constant, it holds that
\begin{align}\label{eq:IIJ(i,j)}
|{\rm II}_5(i,j)|\leq &~\frac{1}{2n_{i,j}}\sum_{k=1}^{n_{i,j}}\frac{1}{N_{i,j,k}}\sum_{K\leq |\ell-k|\leq K+\Delta_K}\bigg(\int_{t_{i,j,k}\wedge t_{i,j,\ell}}^{t_{i,j,k}\vee t_{i,j,\ell}}|\mu_{i,s}|\,{\rm d}s\bigg)\bigg(\int_{t_{i,j,k}\wedge t_{i,j,\ell}}^{t_{i,j,k}\vee t_{i,j,\ell}}|J_{j,s}|\,{\rm d}M_{j,s}\bigg)\notag\\
\leq&~\bigg(\max_{ k\in[n_{i,j}]}\max_{t_{i,j,\ell}\in S_{i,j,k}}\int_{t_{i,j,k}\wedge t_{i,j,\ell}}^{t_{i,j,k}\vee t_{i,j,\ell}}|\mu_{i,s}|\,{\rm d}s\bigg)\notag\\&~\times\frac{1}{2n_{i,j}}\sum_{\ell=1}^{n_{i,j}-1}\bigg[\sum_{k=\ell-K-\Delta_K+1}^{\ell-K}\frac{\min\{k+K+\Delta_K,n_{i,j}\}-\ell}{N_{i,j,k}}\notag\\
&~~~~~~~~~~~~~~~~~~~~+\sum_{k=\ell-K+1}^{\ell}\frac{|\min\{k+K+\Delta_K,n_{i,j}\}-k-K+1|_+}{N_{i,j,k}}\\
&~~~~~~~~~~~~~~~~~~~~+\sum_{k=\ell+1}^{\ell+K}\frac{|k-K+1-\max\{k-K-\Delta_K,1\}|_+}{N_{i,j,k}}\notag\\
&~~~~~~~~~~~~~~~~~~~~+\sum_{k=\ell+K+1}^{\ell+K+\Delta_K}\frac{\ell+1-\max\{k-K-\Delta_K,1\}}{N_{i,j,k}}\bigg]\int_{t_{i,j,\ell}}^{t_{i,j,\ell+1}}|J_{j,s}|\,{\rm d}M_{j,s}\notag\\
\lesssim&~\bigg(\max_{ k\in[n_{i,j}]}\max_{t_{i,j,\ell}\in S_{i,j,k}}\int_{t_{i,j,k}\wedge t_{i,j,\ell}}^{t_{i,j,k}\vee t_{i,j,\ell}}K|\mu_{i,s}|\,{\rm d}s\bigg)\times\frac{1}{2n_{i,j}}\sum_{\ell=1}^{n_{i,j}-1}\int_{t_{i,j,\ell}}^{t_{i,j,\ell+1}}|J_{j,s}|\,{\rm d}M_{j,s}
\notag\\:=&~Q_{i,j}^{\mu}\cdot \frac{1}{2n_{i,j}}\sum_{\ell=1}^{n_{i,j}-1}\int_{t_{i,j,\ell}}^{t_{i,j,\ell+1}}|J_{j,s}|\,{\rm d}M_{j,s}\notag\,,
\end{align}
where we adopt the convention $N_{i,j,k}=\infty$ if $k>n_{i,j}$ or $k<0$. For any constant $d_1\in (0,K^{-1}\xi^{-1}]$, define $\mathcal{E}_{j,d_1}=\{\sup_{0\leq s\leq  T}|J_{j,s}|> d_1\}$.
Recall $(\Delta M_{i,\cdot})_{i,j}=M_{i,t_{i,j,n_{i,j}}}-M_{i,t_{i,j,1}}$ for any $i,j\in[p]$.
By Assumption \ref{as:jumpsize}, for any $v> 0$, we have
\begin{align*}
	&\mathbb{P}\bigg(Q_{i,j}^{\mu}\cdot \frac{1}{2n_{i,j}}\sum_{\ell=1}^{n_{i,j}-1}\int_{t_{i,j,\ell}}^{t_{i,j,\ell+1}}|J_{j,s}|\,{\rm d}M_{j,s}\geq v\bigg)\\&~~~~\leq \mathbb{P}\bigg(Q_{i,j}^{\mu}\cdot \frac{1}{2n_{i,j}}\sum_{\ell=1}^{n_{i,j}-1}\int_{t_{i,j,\ell}}^{t_{i,j,\ell+1}}|J_{j,s}|\,{\rm d}M_{j,s}\geq v,\mathcal{E}_{j,d_1}^c\bigg)+\mathbb{P}(\mathcal{E}_{j,d_1})\\&~~~~\leq
	\mathbb{P}\{d_1Q_{i,j}^{\mu}\cdot (\Delta M_{j,\cdot})_{j,i} \geq 2n_{i,j}v\}+C\exp(-Cd_1^{\iota})\,.
\end{align*}
 By Bonferroni inequality,
$
\mathbb{P}(d_1Q_{i,j}^{\mu}\geq v)\lesssim n_*\max_{ k\in[n_{i,j}]}\max_{t_{i,j,\ell}\in S_{i,j,k}}\mathbb{P}(\int_{t_{i,j,k}\wedge t_{i,j,\ell}}^{t_{i,j,k}\vee t_{i,j,\ell}}Kd_1|\mu_{i,s}|\,{\rm d}s\geq v)
$
for any $v>0$.
For any $\theta\in (0,C_5d_1^{-1}K^{-1}\xi^{-1}]$, by Jensen's inequality, Assumption \ref{as:drifdiff} implies that
$
\mathbb{E}\{\exp(\theta\int_{t_{i,j,k}\wedge t_{i,j,\ell}}^{t_{i,j,k}\vee t_{i,j,\ell}}Kd_1|\mu_{i,s}|\,{\rm d}s)\}\leq {|t_{i,j,\ell}-t_{i,j,k}|^{-1}}\int_{t_{i,j,k}\wedge t_{i,j,\ell}}^{t_{i,j,k}\vee t_{i,j,\ell}}\mathbb{E}\{\exp(\theta Kd_1|t_{i,j,\ell}-t_{i,j,k}||\mu_{i,s}|)\}\,{\rm d}s
\leq  \sup_{0\leq s\leq  T}\mathbb{E}\{\exp(\theta Kd_1\xi|\mu_{i,s}|\,)\}
\lesssim \exp(CK^2d_1^2\xi^2\theta^2)$.
Selecting $\theta\asymp K^{-1}d_1^{-1}\xi^{-1}$, applying Markov's inequality, we have
$
\max_{ k\in[n_{i,j}]}\max_{t_{i,j,\ell}\in S_{i,j,k}}\mathbb{P}(\int_{t_{i,j,k}\wedge t_{i,j,\ell}}^{t_{i,j,k}\vee t_{i,j,\ell}}Kd_1|\mu_{i,s}|\,{\rm d}s\geq v) \lesssim
\exp(-CK^{-1}d_1^{-1}\xi^{-1}v)
$
for any $v>0$, which implies that
\begin{align}\label{eq:Qijmu2}
\max_{ i,j\in[p]}\mathbb{P}(d_1Q_{i,j}^{\mu}\geq v) \lesssim n_*\exp(-CK^{-1}d_1^{-1}\xi^{-1}v)
\end{align}
for any $v> 0$.  By Assumption \ref{as:jump}, applying Proposition 2.9 of \cite{Wainwright_2019}, we have
$
	\mathbb{P}[|(\Delta M_{j,\cdot})_{j,i}-\mathbb{E}\{(\Delta M_{j,\cdot})_{j,i}\}|\geq v]\leq 2\exp[-v^2/\{4\lambda_j(t_{j,i,n_{j,i}}-t_{j,i,1})\}]+2\exp\{-v/(2C_{10})\}
$
for any $v>0$.  By Assumption \ref{as:space} and \eqref{eq:nstar},  $(t_{j,i,n_{j,i}}-t_{j,i,1})\asymp 1$ and $n_{i,j}\asymp n_*$, which implies that
$
\mathbb{P}[|(\Delta M_{j,\cdot})_{j,i}-\mathbb{E}\{(\Delta M_{j,\cdot})_{j,i}\}|\geq n_{i,j}v]\lesssim \exp(-Cn_*^2v^2/\lambda_j)+\exp(-Cn_*v)
$
for any $v>0$. By Assumption \ref{as:jump}(ii), we have $\max_{ i,j\in[p]}|n_{i,j}^{-1}\mathbb{E}\{(\Delta M_{j,\cdot})_{j,i}\}|\lesssim n_*^{-1}\lambda_*$. Then for any $v\gg n_*^{-1}\lambda_*$,  it holds that
\begin{align}\label{eq:mijl}
\max_{ i,j\in[p]}	\mathbb{P}\{|(\Delta M_{j,\cdot})_{j,i}|\geq n_{i,j}v\}\leq&~\max_{ i,j\in[p]}\mathbb{P}[|(\Delta M_{j,\cdot})_{j,i}-\mathbb{E}\{(\Delta M_{j,\cdot})_{j,i}\}|\geq n_{i,j}v-|\mathbb{E}\{(\Delta M_{j,\cdot})_{j,i}\}|]\notag\\\lesssim&~
 \max_{ i,j\in[p]}\mathbb{P}[|(\Delta M_{j,\cdot})_{j,i}-\mathbb{E}\{(\Delta M_{j,\cdot})_{j,i}\}|\geq n_{i,j}v/2]
 \\\lesssim&~ \exp(-C\lambda_*^{-1}n_*^2v^2)+\exp(-Cn_*v)\lesssim \exp(-Cn_*v)\notag\,.
\end{align}
Combining \eqref{eq:Qijmu2} and \eqref{eq:mijl}, we have
\begin{align*}
\mathbb{P}\{d_1Q_{i,j}^{\mu} (\Delta M_{j,\cdot})_{j,i}\geq 2n_{i,j}v\}\leq&~ \mathbb{P}(d_1Q_{i,j}^{\mu}\geq K)+\mathbb{P}(|(\Delta M_{j,\cdot})_{j,i}|\geq 2K^{-1}n_{i,j}v)\\\lesssim&~
n_*\exp(-Cd_1^{-1}\xi^{-1})+\exp(-Cn_*K^{-1}v)
\end{align*}
for any $v\gg n_*^{-1}K\lambda_*$. Then
$
\mathbb{P}\{Q_{i,j}^{\mu}\cdot (2n_{i,j})^{-1}\sum_{\ell=1}^{n_{i,j}-1}\int_{t_{i,j,\ell}}^{t_{i,j,\ell+1}}|J_{j,s}|\,{\rm d}M_{j,s}\geq v\}\lesssim n_*\exp(-Cd_1^{-1}\xi^{-1})\notag+\exp(-Cn_*K^{-1}v)+\exp(-Cd_1^{\iota})
$
for any $v\gg n_*^{-1}K\lambda_*$. Notice that $\xi\asymp n_*^{-1}K$. Thus, by \eqref{eq:IIJ(i,j)},
$
\max_{ i,j\in[p]}\mathbb{P}\{|{\rm II}_5(i,j)|\geq v\}\lesssim \exp(-Cn_*K^{-1}v)+n_*\exp(-Cd_1^{-1}n_*K^{-1})+\exp(-Cd_1^{\iota})
$
for any $v\gg n_*^{-1}K\lambda_*$. In order to make $\max_{ i,j\in[p]}|{\rm II}_5(i,j)|=O_{\p}(n_*^{-1/2}K^{1/2}\log ^{1/2}p)$ and  $p$ diverge as fast as possible,
 we require $\log p=o(n_*K^{-1})$ and
 select $d_1\asymp (n_*K^{-1})^{1/(1+\iota)}$. Then if $K(\log n_*)^{(1+\iota)/\iota}=o(n_*)$, it holds that
\begin{align}\label{eq:II5ij}
\max_{ i,j\in[p]}\mathbb{P}\{|{\rm II}_5(i,j)|\geq v\}\lesssim \exp(-Cn_*K^{-1}v)+\exp\{-C(n_*K^{-1})^{\iota/(\iota+1)}\}
\end{align}
for any $v\gg n_*^{-1}K\lambda_*$.
Identically, if $K(\log n_*)^{(1+\iota)/\iota}=o(n_*)$, we also have
\begin{align}\label{eq:II8ij}
\max_{ i,j\in[p]}\mathbb{P}\{|{\rm II}_8(i,j)|\geq v\}\lesssim \exp(-Cn_*K^{-1}v)+\exp\{-C(n_*K^{-1})^{\iota/(\iota+1)}\}
\end{align}
for any $v\gg n_*^{-1}K\lambda_*$.

 Now we consider ${\rm II}_6(i,j)$ and ${\rm II}_7(i,j)$.
Analogously to \eqref{eq:IIJ(i,j)}, we have
\begin{align}\label{eq:IIJ(i,j)6}
|{\rm II}_6(i,j)|\leq &~\frac{1}{2n_{i,j}}\sum_{k=1}^{n_{i,j}}\frac{1}{N_{i,j,k}}\sum_{K\leq |\ell-k|\leq K+\Delta_K}\bigg|\int_{t_{i,j,k}\wedge t_{i,j,\ell}}^{t_{i,j,k}\vee t_{i,j,\ell}}\sigma_{i,s}\,{\rm d}B_{i,s}\bigg|\bigg(\int_{t_{i,j,k}\wedge t_{i,j,\ell}}^{t_{i,j,k}\vee t_{i,j,\ell}}|J_{j,s}|\,{\rm d}M_{j,s}\bigg)\notag\\
\lesssim&~\bigg(\max_{ k\in[n_{i,j}]}\max_{t_{i,j,\ell}\in S_{i,j,k}}K\bigg|\int_{t_{i,j,k}\wedge t_{i,j,\ell}}^{t_{i,j,k}\vee t_{i,j,\ell}}\sigma_{i,s}\,{\rm d}B_{i,s}\bigg|\bigg)\times\frac{1}{2n_{i,j}}\sum_{\ell=1}^{n_{i,j}-1}\int_{t_{i,j,\ell}}^{t_{i,j,\ell+1}}|J_{j,s}|\,{\rm d}M_{j,s}
\\:=&~Q_{i,j}^{\sigma}\cdot \frac{1}{2n_{i,j}}\sum_{\ell=1}^{n_{i,j}-1}\int_{t_{i,j,\ell}}^{t_{i,j,\ell+1}}|J_{j,s}|\,{\rm d}M_{j,s}\notag\,.
\end{align}
For any constant $d_2\in (0,K^{-1}\xi^{-1}]$, define $\mathcal{E}_{j,d_2}=\{\sup_{0\leq s\leq  T}|J_{j,s}|> d_2\}$.
By Assumption \ref{as:jumpsize},
\begin{align*}
&\mathbb{P}\bigg(Q_{i,j}^{\sigma}\cdot \frac{1}{2n_{i,j}}\sum_{\ell=1}^{n_{i,j}-1}\int_{t_{i,j,\ell}}^{t_{i,j,\ell+1}}|J_{j,s}|\,{\rm d}M_{j,s}\geq v\bigg)\\&~~~~\leq \mathbb{P}\bigg(Q_{i,j}^{\sigma}\cdot \frac{1}{2n_{i,j}}\sum_{\ell=1}^{n_{i,j}-1}\int_{t_{i,j,\ell}}^{t_{i,j,\ell+1}}|J_{j,s}|\,{\rm d}M_{j,s}\geq v,\mathcal{E}_{j,d_2}^c\bigg)+\mathbb{P}(\mathcal{E}_{j,d_2})\\&~~~~\leq
\mathbb{P}\{d_2Q_{i,j}^{\sigma}\cdot (\Delta M_{j,\cdot})_{j,i}\geq 2n_{i,j}v\}+C\exp(-Cd_2^{\iota})
\end{align*}
for any $v> 0$.   By Bonferroni inequality, we have
\begin{align}\label{eq:Qijsigma}
\mathbb{P}(d_2Q_{i,j}^{\sigma}\geq v)\leq&~ \sum_{k=1}^{n_{i,j}}\sum_{t_{i,j,\ell}\in S_{i,j,k}}\mathbb{P}\bigg(Kd_2\bigg|\int_{t_{i,j,k}\wedge t_{i,j,\ell}}^{t_{i,j,k}\vee t_{i,j,\ell}}\sigma_{i,s}\,{\rm d}B_{i,s}\bigg|\geq v\bigg)\notag\\\lesssim &~n_*\max_{ k\in[n_{i,j}]}\max_{t_{i,j,\ell}\in S_{i,j,k}}\mathbb{P}\bigg(d_2K\bigg|\int_{t_{i,j,k}\wedge t_{i,j,\ell}}^{t_{i,j,k}\vee t_{i,j,\ell}}\sigma_{i,s}\,{\rm d}B_{i,s}\bigg|\geq v\bigg)
\end{align}
for any $v>0$.
 For any constant $d\in(0,d_2^{-1/2}K^{-1/2}\xi^{-1/2}]$, define a stopping time $
\Gamma_{i,d}=T\wedge\inf\{t>0:\sup_{0\leq s\leq t}\sigma_{i,s}>d\}$.
Then it holds that
\begin{align}\label{eq:sigmaJ}
&\mathbb{P}\bigg(d_2K\bigg|\int_{t_{i,j,k}\wedge t_{i,j,\ell}}^{t_{i,j,k}\vee t_{i,j,\ell}}\sigma_{i,s}\,{\rm d}B_{i,s}\bigg|\geq v,\Gamma_{i,d}=T\bigg)\notag\\&~~~~~~\lesssim \exp(-C\theta v^2)\mathbb{E}\bigg[\exp\bigg\{\theta d_2^2K^2\bigg(\int_{t_{i,j,k}\wedge t_{i,j,\ell}}^{t_{i,j,k}\vee t_{i,j,\ell}}\sigma_{i,s}\,{\rm d}B_{i,s}\bigg)^2\bigg\}I(\Gamma_{i,d}=T)\bigg]
\end{align}
for any $\theta>0$.
Identical to (\ref{eq:re1}) and (\ref{eq:re2}),
$
\mathbb{E}[\exp\{\theta d_2^2K^2(\int_{t_{i,j,k}\wedge t_{i,j,\ell}}^{t_{i,j,k}\vee t_{i,j,\ell}}\sigma_{i,s}\,{\rm d}B_{i,s})^2\}{I}(\Gamma_{i,d}=T)]\lesssim 1$.
With $\theta=  d_2^{-2}d^{-2}K^{-2}\xi^{-1}/4$,  \eqref{eq:sigmaJ} implies  $$\mathbb{P}\bigg(d_2K\bigg|\int_{t_{i,j,k}\wedge t_{i,j,\ell}}^{t_{i,j,k}\vee t_{i,j,\ell}}\sigma_{i,s}\,{\rm d}B_{i,s}\bigg|\geq v,\Gamma_{i,d}=T\bigg)\lesssim \exp(-Cd_2^{-2}d^{-2}K^{-2}\xi^{-1}v^2)\,. $$
By Assumption \ref{as:diff},
$
\max_{ k\in[n_{i,j}]}\max_{t_{i,j,\ell}\in S_{i,j,k}}\mathbb{P}(d_2K|\int_{t_{i,j,k}\wedge t_{i,j,\ell}}^{t_{i,j,k}\vee t_{i,j,\ell}}\sigma_{i,s}\,{\rm d}B_{i,s}|\geq v)\lesssim
\exp(-Cd_2^{-2}d^{-2}K^{-2}\xi^{-1}v^2)+\exp(-Cd^{\gamma})
$
for any $v>0$. Then \eqref{eq:Qijsigma} implies that
$
\mathbb{P}(d_2Q_{i,j}^{\sigma}\geq v)\lesssim n_*\exp(-Cd_2^{-2}d^{-2}K^{-2}\xi^{-1}v^2)+n_*\exp(-Cd^{\gamma})
$
for any $v>0$. Together with \eqref{eq:mijl}, we have
$
\mathbb{P}\{d_2Q_{i,j}^{\sigma} (\Delta M_{j,\cdot})_{j,i}\geq 2n_{i,j}v\}\leq \mathbb{P}(d_2Q_{i,j}^{\sigma}\geq K)+\mathbb{P}\{|(\Delta M_{j,\cdot})_{j,i}|\geq 2K^{-1}n_{i,j}v\}\lesssim
n_*\exp(-Cd_2^{-2}d^{-2}\xi^{-1})+n_*\exp(-Cd^{\gamma})+\exp(-Cn_*K^{-1}v)$
for any $v\gg n_*^{-1}K\lambda_*$, which implies
$\max_{ i,j\in[p]}
\mathbb{P}\{|{\rm II}_6(i,j)|\geq v\}\lesssim n_*\exp(-Cd_2^{-2}d^{-2}\xi^{-1})+n_*\exp(-Cd^{\gamma})+\exp(-Cd_2^{\iota})+\exp(-Cn_*K^{-1}v)
$
for any $v\gg n_*^{-1}K\lambda_*$. To make $\max_{ i,j\in[p]}|{\rm II}_6(i,j)|=O_{\p}(n_*^{-1/2}K^{1/2}\log ^{1/2}p)$ and  $p$ diverge as fast as possible, we require $\log p=o(n_*K^{-1})$, and  select $d_2\asymp (n_*K^{-1})^{\gamma/(2\iota+2\gamma+\iota\gamma)}$ and $d\asymp (d_2^{-2}n_*K^{-1})^{1/(\gamma+2)}$. Then if $K (\log n_*)^{(2\iota+2\gamma+\iota\gamma)/(\iota\gamma)}=o(n_*)$,  we have
\begin{align}\label{eq:II6(ij)}
\max_{ i,j\in[p]}\mathbb{P}\{|{\rm II}_6(i,j)|\geq v\}\lesssim \exp(-Cn_*K^{-1}v) +
\exp\{-C(n_*K^{-1})^{\iota\gamma/(2\iota+2\gamma+\iota\gamma)}\}
\end{align}
for any $v\gg n_*^{-1}K\lambda_*$.
Analogously, if $K (\log n_*)^{(2\iota+2\gamma+\iota\gamma)/(\iota\gamma)}=o(n_*)$, it holds that
\begin{align}\label{eq:II7(ij)}
\max_{ i,j\in[p]}\mathbb{P}\{|{\rm II}_7(i,j)|\geq v\}\lesssim \exp(-Cn_*K^{-1}v) +
\exp\{-C(n_*K^{-1})^{\iota\gamma/(2\iota+2\gamma+\iota\gamma)}\}
\end{align}
for any $v\gg n_*^{-1}K\lambda_*$.
 Notice that  $\iota\gamma/(2\iota+2\gamma+\iota\gamma)\leq \iota/(1+\iota)$.  Combining  \eqref{eq:II5ij}, \eqref{eq:II8ij}, \eqref{eq:II6(ij)} and \eqref{eq:II7(ij)}, if $K (\log n_*)^{(2\iota+2\gamma+\iota\gamma)/(\iota\gamma)}=o(n_*)$, it holds that
$
\max_{i,j\in[p]}\mathbb{P}\{|{\rm II}''(i,j)|>v\}\lesssim~ \exp(-Cn_*K^{-1}v)+
\exp\{-C(n_*K^{-1})^{\iota\gamma/(2\iota+2\gamma+\iota\gamma)}\}
$ for any
$v\gg n_*^{-1}K\lambda_*$.

By \eqref{eq:IIJ(i,j)} and Cauchy-Schwarz inequality, for any positive integer $m$, we have  $\mathbb{E}\{|{\rm II}_5(i,j)|^{m}\}\lesssim [\mathbb{E}\{(K^{-1}Q_{i,j}^{\mu})^{2m}\}]^{1/2}\cdot [\mathbb{E}\{(n_{i,j}^{-1}K\sum_{\ell=1}^{n_{i,j}-1}\int_{t_{i,j,\ell}}^{t_{i,j,\ell+1}}|J_{j,s}|\,{\rm d}M_{j,s})^{2m}\}]^{1/2}=:E_{i,j,1}\cdot E_{i,j,2}$. By  Assumption \ref{as:drifdiff} and Jensen's inequality, it holds that
 \begin{align*}
 	E_{i,j,1}^2\lesssim&~ n_*\max_{ k\in[n_{i,j}]}\max_{t_{i,j,\ell}\in S_{i,j,k}} \mathbb{E} \bigg\{\bigg(\int_{t_{i,j,k}\wedge t_{i,j,\ell}}^{t_{i,j,k}\vee t_{i,j,\ell}}|\mu_{i,s}|\,{\rm d}s\bigg)^{2m}\bigg\}
 	\\\lesssim&~
 n_*\xi^{2m-1}\max_{ k\in[n_{i,j}]}\max_{t_{i,j,\ell}\in S_{i,j,k}} \int_{t_{i,j,k}\wedge t_{i,j,\ell}}^{t_{i,j,k}\vee t_{i,j,\ell}}\sup_{0\leq s\leq  T}\mathbb{E}(|\mu_{i,s}|^{2m})\,{\rm d}s\lesssim n_*\xi^{2m}\,.
 \end{align*}
 If $\lambda_*n_*^{-1}K=o(1)$, it follows from Assumptions \ref{as:jump} and \ref{as:jumpsize} that
\begin{align}\label{eq:D2}
E_{i,j,2}^2\leq &~\mathbb{E}\bigg[\bigg(\sup_{0\leq s\leq  T}|J_{i,s}|^{2m}\bigg)\{Kn_{i,j}^{-1}(\Delta M_{j,\cdot})_{j,i}\}^{2m}\bigg]\notag
\\\lesssim&~
\bigg\{\mathbb{E}\bigg(\sup_{0\leq s\leq  T}|J_{i,s}|^{4m}\bigg)\bigg\}^{1/2}\cdot (\mathbb{E}[\{Kn_{i,j}^{-1}(\Delta M_{j,\cdot})_{j,i}\}^{4m}])^{1/2}
\\\lesssim&~(\mathbb{E}[\exp\{Kn_{i,j}^{-1}(\Delta M_{j,\cdot})_{j,i}\}])^{1/2}\lesssim 1\notag\,.
\end{align}
Thus, $\max_{ i,j\in[p]}\mathbb{E}\{|{\rm II}_5(i,j)|^{m}\}\lesssim n_*^{1/2}\xi^m$. Analogously, we have  $\max_{ i,j\in[p]}\mathbb{E}\{|{\rm II}_8(i,j)|^{m}\}\lesssim n_*^{1/2}\xi^m$.
By \eqref{eq:IIJ(i,j)6} and Cauchy-Schwarz inequality, it holds that  $\mathbb{E}\{|{\rm II}_6(i,j)|^{m}\}\lesssim [\mathbb{E}\{(K^{-1}Q_{i,j}^{\sigma})^{2m}\}]^{1/2}\cdot [\mathbb{E}\{(n_{i,j}^{-1}K\sum_{\ell=1}^{n_{i,j}-1}\int_{t_{i,j,\ell}}^{t_{i,j,\ell+1}}|J_{j,s}|\,{\rm d}M_{j,s})^{2m}\}]^{1/2}=:E_{i,j,3}\cdot E_{i,j,2}$ for any positive integer $m$.
By Assumptions \ref{as:drifdiff} and  \ref{as:diff}, Burkholder-Davis-Gundy inequality implies
\begin{align*}
E_{i,j,3}^2\lesssim&~ n_*\max_{ k\in[n_{i,j}]}\max_{t_{i,j,\ell}\in S_{i,j,k}}\mathbb{E}\bigg\{\bigg(\int_{t_{i,j,k}\wedge t_{i,j,\ell}}^{t_{i,j,k}\vee t_{i,j,\ell}}\sigma_{i,s}\,{\rm d}B_{i,s}\bigg)^{2m}\bigg\}
\lesssim n_*\xi^{m}\,.
 \end{align*}
Together with \eqref{eq:D2}, $\max_{ i,j\in[p]}\mathbb{E}\{|{\rm II}_6(i,j)|^{m}\}\lesssim n_*^{1/2}\xi^{m/2}$. Analogously, $\max_{ i,j\in[p]}\mathbb{E}\{|{\rm II}_7(i,j)|^{m}\}\lesssim n_*^{1/2}\xi^{m/2}$. Hence, $\max_{ i,j\in[p]}\mathbb{E}\{|{\rm II}''(i,j)|^{m}\}\lesssim n_*^{1/2}\xi^{m/2}$
for any positive integer $m$ provided that $\lambda_*n_*^{-1}K=o(1)$.
 We complete the proof of Proposition \ref{eq:proj4}.  $\hfill\Box$

\subsection{Proof of Proposition \ref{eq:pfpro5} }\label{sec:pro5}
 Notice that for sufficiently large $n$, $\min_{k \in[n_{i,j}]}N_{i,j,k}=\Delta_K+1$ and $\max_{k\in[n_{i,j}]}N_{i,j,k}=2\Delta_K+2$. Analogously to \eqref{eq:IIJ(i,j)}, we have that
\begin{align}\label{eq:III2J(i,j)}
|{\rm III}''(i,j)|\lesssim &~\frac{1}{2n_{i,j}}\sum_{k=1}^{n_{i,j}}\sum_{K\leq |\ell-k|\leq K+\Delta_K}\bigg(\int_{t_{i,j,k}\wedge t_{i,j,\ell}}^{t_{i,j,k}\vee t_{i,j,\ell}}|J_{i,s}|\,{\rm d}M_{i,s}\bigg)|U_{j,t_{i,j,k}}|\notag
\\\lesssim&~\bigg(K\max_{ k\in[n_{i,j}]}|U_{j,t_{i,j,k}}|\bigg)\times\frac{1}{2n_{i,j}}\sum_{\ell=1}^{n_{i,j}-1}\int_{t_{i,j,\ell}}^{t_{i,j,\ell+1}}|J_{i,s}|\,{\rm d}M_{i,s}
\,.
\end{align}
For any constant $d_3>0$, define $\mathcal{E}_{i,d_3}=\{\sup_{0\leq s\leq  T}|J_{i,s}|> d_3\}$.
By Assumption \ref{as:jumpsize}, it holds that
\begin{align}\label{eq:UijJ1}
&\mathbb{P}\bigg\{\bigg(K\max_{ k\in[n_{i,j}]}|U_{j,t_{i,j,k}}|\bigg)\cdot \frac{1}{2n_{i,j}}\sum_{\ell=1}^{n_{i,j}-1}\int_{t_{i,j,\ell}}^{t_{i,j,\ell+1}}|J_{i,s}|\,{\rm d}M_{i,s}\geq v\bigg\}\notag\\&~~~~\leq \mathbb{P}\bigg\{\bigg(K\max_{ k\in[n_{i,j}]}|U_{j,t_{i,j,k}}|\bigg)\cdot \frac{1}{2n_{i,j}}\sum_{\ell=1}^{n_{i,j}-1}\int_{t_{i,j,\ell}}^{t_{i,j,\ell+1}}|J_{i,s}|\,{\rm d}M_{i,s}\geq v,\mathcal{E}_{i,d_3}^c\bigg\}+\mathbb{P}(\mathcal{E}_{i,d_3})\\&~~~~\leq
\mathbb{P}\bigg\{\bigg(d_3K\max_{ k\in[n_{i,j}]}|U_{j,t_{i,j,k}}|\bigg) (\Delta M_{i,\cdot})_{i,j}\geq 2n_{i,j}v\bigg\}+C\exp(-Cd_3^{\iota})\notag
\end{align}
for any $v> 0$.
By Bonferroni inequality and Assumption \ref{as:moment},
$
\mathbb{P}(d_3K\max_{ k\in[n_{i,j}]}|U_{j,t_{i,j,k}}|\geq v)\leq \sum_{k=1}^{n_{i,j}}\mathbb{P}(|U_{j,t_{i,j,k}}|\geq d_3^{-1}K^{-1}v)\lesssim
n_*\max_{ k\in[n_{i,j}]}\mathbb{P}(|U_{j,t_{i,j,k}}|\geq d_3^{-1}K^{-1}v)
\lesssim n_*\exp(-Cd_3^{-2}K^{-2}v^2)
$
for any $v>0$. Recall $\xi\asymp n_*^{-1}K$. Together with \eqref{eq:mijl}, for any  $\delta\in(0,1)$, we have that
$
\mathbb{P}\{(d_3K\max_{ k\in[n_{i,j}]}|U_{j,t_{i,j,k}}|)(\Delta M_{i,\cdot})_{i,j}\geq 2n_{i,j}v\}\leq  \mathbb{P}(d_3K\max_{ k\in[n_{i,j}]}|U_{j,t_{i,j,k}}|\geq K\xi^{-\delta})+\mathbb{P}\{ (\Delta M_{i,\cdot})_{i,j}\geq 2n_{i,j}K^{-1}\xi^\delta v\}\lesssim n_*\exp(-Cd_3^{-2}\xi^{-2\delta})+\exp(-C\xi^{\delta-1}v)
$
for any $v\gg \xi^{1-\delta}\lambda_*$.
By \eqref{eq:III2J(i,j)} and  \eqref{eq:UijJ1},
$
\max_{ i,j\in[p]}\mathbb{P}\{|{\rm III}''(i,j)|\geq v\}\lesssim n_*\exp(-Cd_3^{-2}\xi^{-2\delta})+\exp(-C\xi^{\delta-1}v)+\exp(-Cd_3^{\iota})
$
for any $v\gg \xi^{1-\delta}\lambda_*$.
To make $\max_{ i,j\in[p]}|{\rm III}''(i,j)|=O_{\p}(n_*^{-1/2}K^{1/2}\log ^{1/2}p)$, it suffices to require
$\log p=o[\min\{d_3^{-2}\xi^{-2\delta},d_3^\iota,\xi^{2\delta-1}\}]$, $\lambda_*^2\xi^{1-2\delta}(\log p)^{-1}=o(1)$ and $\log n_*=o(d_3^{-2}\xi^{-2\delta})$. To make $p$ diverge as fast as possible, we select $\delta$ such that $d_3^{-2}\xi^{-2\delta}\asymp d_3^\iota\asymp\xi^{2\delta-1}$. Then $\delta =(\iota+2)/(4\iota+4)$ and $d_3=\xi^{-1/(2\iota+2)}$.
 Hence, if $K(\log n_*)^{2(\iota+1)/\iota}=o(n_*)$,  we have
$
\max_{i,j\in[p]}\mathbb{P}\{|{\rm{III}}''(i,j)|\geq v\}\lesssim\exp\{-C(n_*K^{-1})^{(3\iota+2)/(4\iota+4)}v\}+\exp\{-C(n_*K^{-1})^{\iota/(2\iota+2)}\}
$
for any $v\gg (n_*^{-1}K)^{(3\iota+2)/(4\iota+4)}\lambda_*$.

By \eqref{eq:III2J(i,j)} and Cauchy-Schwarz inequality,   \begin{align*}
\mathbb{E}\{|{\rm III}''(i,j)|^{m}\}\lesssim &~ \bigg\{\mathbb{E}\bigg(\max_{ k\in[n_{i,j}]}|U_{j,t_{i,j,k}}|^{2m}\bigg)\bigg\}^{1/2}\cdot \bigg[\mathbb{E}\bigg\{\bigg(n_{i,j}^{-1}K\sum_{\ell=1}^{n_{i,j}-1}\int_{t_{i,j,\ell}}^{t_{i,j,\ell+1}}|J_{i,s}|\,{\rm d}M_{i,s}\bigg)^{2m}\bigg\}\bigg]^{1/2}\\=:&~F_{i,j,1}\cdot F_{i,j,2}
\end{align*} for any positive integer $m$.
It follows Assumption \ref{as:moment} that $F_{i,j,1}^2\lesssim n_*\max_{ k\in[n_{i,j}]}\mathbb{E}(U_{j,t_{i,j,k}}^{2m})\lesssim n_*$. Analogously to \eqref{eq:D2}, we have $F_{i,j,2}^2\lesssim 1$ if  $\lambda_*n_*^{-1}K=o(1)$. Hence,
$\max_{ i,j\in[p]}\mathbb{E}\{|{\rm III}''(i,j)|^{m}\}\lesssim n_*^{1/2}$ for any positive integer $m$ provided that $\lambda_*n_*^{-1}K=o(1)$.
 We complete the proof of Proposition \ref{eq:pfpro5}.  $\hfill\Box$

\subsection{Proof of Part (ii) in Theorem \ref{eq:the6}}

We first consider the case with $\varphi<\infty$. Write $\aleph=(n_*^{-1}K\log p)^{1/2}$. For each $i,j\in[p]$, we  define the event
$
A_{i,j}=\{|\hat{\sigma}_{u,i,j}^{{\rm jump,thre}}-\sigma_{u,i,j}|\leq 4\min(|\sigma_{u,i,j}|,\alpha\aleph)\}
$
with some constant $\alpha>0$, and
$
d_{i,j}=(\hat{\sigma}_{u,i,j}^{{\rm jump,thre}}-\sigma_{u,i,j}){I}(A_{i,j}^c)$.
Identical to \eqref{eq:speupp1}, we have
$
\mathbb{E}(\|\widehat{\bSigma}_{u}^{{\rm jump,  thre}}-\bSigma_{u}\|_2^2)\lesssim \mathbb{E}\{(\max_{i\in[p]}\sum_{j=1}^p|d_{i,j}|)^2\}+c_p^2\aleph^{2(1-q)}$.
 It holds that
\begin{align*}
\mathbb{E}\bigg\{\bigg(\max_{i\in[p]}\sum_{j=1}^p|d_{i,j}|\bigg)^2\bigg\}\leq&~p\sum_{i,j=1}^p\mathbb{E}\big\{|\hat{\sigma}_{u,i,j}^{{\rm jump,thre}}-\sigma_{u,i,j}|^2{I}(A_{i,j}^c)\big\}\\
=&~\underbrace{p\sum_{i,j=1}^p\mathbb{E}\big(|\hat{\sigma}_{u,i,j}^{{\rm jump,thre}}-\sigma_{u,i,j}|^2{I}[A_{i,j}^c\cap\{\hat{\sigma}_{u,i,j}^{{\rm jump,thre}}=0\}]\big)}_{\textrm{I}}\\
&+\underbrace{p\sum_{i,j=1}^p\mathbb{E}\big(|\hat{\sigma}_{u,i,j}^{{\rm jump,thre}}-\sigma_{u,i,j}|^2{I}[A_{i,j}^c\cap\{\hat{\sigma}_{u,i,j}^{{\rm jump,thre}}=\hat{\sigma}_{u,i,j}^{{\rm jump}}\}]\big)}_{\textrm{II}}\,.
\end{align*}
Recall $\hat{\sigma}_{u,i,j}^{{\rm jump,thre}}=\hat{\sigma}_{u,i,j}^{{\rm jump}}{I}(|\hat{\sigma}_{u,i,j}^{{\rm jump}}|\geq \beta\aleph)$ for any $i,j\in[p]$. Identical to \eqref{eq:th3I},
$
{\rm I}
\leq p\sum_{i,j=1}^p\sigma_{u,i,j}^2\mathbb{P}\{|\hat{\sigma}_{u,i,j}^{{\rm jump}}-\sigma_{u,i,j}|\geq(4\alpha-\beta)\aleph\}$.
Selecting $\alpha=\beta/2$ and $\beta$ being sufficiently large,  identical to the arguments used in Section \ref{sec:prothe6} for bounding the convergence rate of $A_2^*$, we know
$
{\rm I}\lesssim  o(\aleph^4)$
provided that $K^{-\varphi}L_n^{\varphi}\log\{n_*(K\log p)^{-1}\}=o(1)$, $\lambda_*^2(n_*^{-1}K)^{\iota/(2\iota+2)} (\log p)^{-1}=o(1)$, $K\gtrsim L_n$, $K (\log n_*)^{\chi_1}=o(n_*)$ with $\chi_1=\max\{(2\iota+2\gamma+\iota\gamma)/(\iota\gamma),2+2/\iota,1+2/\gamma\}$, and $\log p=o[\min\{(n_*L_n^{-2}K)^{\varphi/(3\varphi+2)},(n_*K^{-1})^{\chi_2}\}]$ with $\chi_2=\min\{\iota\gamma/(2\iota+2\gamma+\iota\gamma),\iota/(2\iota+2),\gamma/(\gamma+4),1/3\}$. Also, identical to \eqref{eq:IIthm3},
$
{\rm II}
\leq p\sum_{i,j=1}^p\{\mathbb{E}(|\hat{\sigma}_{u,i,j}^{{\rm jump}}-\sigma_{u,i,j}|^4)\}^{1/2}\{\mathbb{P}(|\hat{\sigma}_{u,i,j}^{{\rm jump}}-\sigma_{u,i,j}|>4\alpha\aleph)\}^{1/2}+p\sum_{i,j=1}^p\{\mathbb{E}(|\hat{\sigma}_{u,i,j}^{{\rm jump}}-\sigma_{u,i,j}|^4)\}^{1/2}[\mathbb{P}\{|\hat{\sigma}_{u,i,j}^{{\rm jump}}-\sigma_{u,i,j}|>(\beta-\alpha)\aleph\}]^{1/2}$.
Notice that
$
\max_{i,j\in[p]}\mathbb{E}(|\hat{\sigma}_{u,i,j}^{{\rm jump}}-\sigma_{u,i,j}|^4)\lesssim n_*^{1/2}\leq  p^{2c_*}$ with $c_*=1/(4\kappa)$, where $\kappa$ is specified in the beginning of Section \ref{sec:theoryfU}.
Since $\alpha=\beta/2$, repeating the arguments used in Section \ref{sec:prothe6} for bounding the convergence rate of $A_2^*$ again, it holds that ${\rm II}\lesssim p^{3+c_*}[\max_{i,j\in[p]}\mathbb{P}\{|\hat{\sigma}_{u,i,j}^{{\rm jump}}-\sigma_{u,i,j}|>(\beta-\alpha)\aleph\}]^{1/2}\lesssim o(\aleph^2)$
if  $K^{-\varphi}L_n^{\varphi}\log\{n_*(K\log p)^{-1}\}=o(1)$, $\lambda_*^2(n_*^{-1}K)^{\iota/(2\iota+2)} (\log p)^{-1}=o(1)$, $K\gtrsim L_n$, $K (\log n_*)^{\chi_1}=o(n_*)$ with $\chi_1=\max\{(2\iota+2\gamma+\iota\gamma)/(\iota\gamma),2+2/\iota,1+2/\gamma\}$, $\log p=o[\min\{(n_*L_n^{-2}K)^{\varphi/(3\varphi+2)},(n_*K^{-1})^{\chi_2}\}]$ with $\chi_2=\min\{\iota\gamma/(2\iota+2\gamma+\iota\gamma),\iota/(2\iota+2),\gamma/(\gamma+4),1/3\}$.
Therefore,
$
\mathbb{E}\{(\max_{i\in[p]}\sum_{j=1}^p|d_{i,j}|)^2\}\leq{\rm I}+{\rm II}\lesssim o(\aleph
^2)$, which implies
$
\sup_{\mathcal{P}_4}\mathbb{E}(\|\widehat{\bSigma}_{u}^{{\rm jump,  thre}}-{\bf\Sigma}_{u}\|_2^2)\lesssim c_p^2(n_*^{-1}K\log p)^{1-q}
$. Analogously, in the case with $\varphi=\infty$, we have $\sup_{\mathcal{P}_4}\mathbb{E}(\|\widehat{\bSigma}_{u}^{{\rm jump,  thre}}-{\bf\Sigma}_{u}\|_2^2)\lesssim c_p^2(n_*^{-1}K\log p)^{1-q}$ provided that $K>L_n$, $\lambda_*^2(n_*^{-1}K)^{\iota/(2\iota+2)} (\log p)^{-1}=o(1)$, $K (\log n_*)^{\chi_1}=o(n_*)$ with $\chi_1=\max\{(2\iota+2\gamma+\iota\gamma)/(\iota\gamma),2+2/\iota,1+2/\gamma\}$ and $\log p=o\{(n_*K^{-1})^{\chi_2}\}$ with $\chi_2=\min\{\iota\gamma/(2\iota+2\gamma+\iota\gamma),\iota/(2\iota+2),\gamma/(\gamma+4),1/3\}$.
We complete the proof of Part (ii) in Theorem \ref{eq:the6}. $\hfill\Box$

\subsection{Proof of Theorem \ref{thwij}}
 Notice that for sufficiently large $n$, $\min_{k \in[n_{i,j}]}N_{i,j,k}=\Delta_K+1$ and $\max_{k\in[n_{i,j}]}N_{i,j,k}=2\Delta_K+2$. Analogously to \eqref{eq:IIJ(i,j)}, we have that
\begin{align}\label{eq:w(i,j)}
|\varpi_{i,j}|\leq  &~\frac{1}{2n_{i,j}}\sum_{k=1}^{n_{i,j}}\sum_{K\leq |\ell-k|\leq K+\Delta_K}\bigg(\int_{t_{i,j,k}\wedge t_{i,j,\ell}}^{t_{i,j,k}\vee t_{i,j,\ell}}|J_{i,s}|\,{\rm d}M_{i,s}\bigg)\bigg(\int_{t_{i,j,k}\wedge t_{i,j,\ell}}^{t_{i,j,k}\vee t_{i,j,\ell}}|J_{j,s}|\,{\rm d}M_{j,s}\bigg)\notag
\\\lesssim&~\bigg(K\max_{ k\in[n_{i,j}]}\max_{t_{i,j,\ell}\in S_{i,j,k}}\int_{t_{i,j,k}\wedge t_{i,j,\ell}}^{t_{i,j,k}\vee t_{i,j,\ell}}|J_{j,s}|\,{\rm d}M_{j,s}\bigg)\times\frac{1}{2n_{i,j}}\sum_{\ell=1}^{n_{i,j}-1}\int_{t_{i,j,\ell}}^{t_{i,j,\ell+1}}|J_{i,s}|\,{\rm d}M_{i,s}
\\:=&~Q_{i,j}^{J}\cdot\frac{1}{2n_{i,j}}\sum_{\ell=1}^{n_{i,j}-1}\int_{t_{i,j,\ell}}^{t_{i,j,\ell+1}}|J_{i,s}|\,{\rm d}M_{i,s}\notag
\,.
\end{align}
For any constant $d\geq 1$ and $i\in[p]$, define $\mathcal{E}_{i,d}=\{\sup_{0\leq s\leq  T}|J_{i,s}|> d\}$.	Let $(\Delta M_{j,\cdot})_{j,i}^{(k,\ell)}=\int_{t_{i,j,k}\wedge t_{i,j,\ell}}^{t_{i,j,k}\vee t_{i,j,\ell}}\,{\rm d}M_{j,s}$ for any $k\in[n_{i,j}]$ and $t_{i,j,\ell}\in S_{i,j,k}$. Recall $(\Delta M_{i,\cdot})_{i,j}=\int_{t_{i,j,1}}^{t_{i,j,n_{i,j}}}\,{\rm d}M_{i,s}$.
By Assumption \ref{as:jumpsize}, it holds that
\begin{align}\label{eq:qij}
&\mathbb{P}\bigg(Q_{i,j}^{J}\cdot \frac{1}{2n_{i,j}}\sum_{\ell=1}^{n_{i,j}-1}\int_{t_{i,j,\ell}}^{t_{i,j,\ell+1}}|J_{i,s}|\,{\rm d}M_{i,s}\geq v\bigg)\notag\\&~~~~\leq \mathbb{P}\bigg(Q_{i,j}^{J}\cdot \frac{1}{2n_{i,j}}\sum_{\ell=1}^{n_{i,j}-1}\int_{t_{i,j,\ell}}^{t_{i,j,\ell+1}}|J_{i,s}|\,{\rm d}M_{i,s}\geq v\,,\mathcal{E}_{i,d}^c\,,\mathcal{E}_{j,d}^c\bigg)+\mathbb{P}(\mathcal{E}_{i,d})+\mathbb{P}(\mathcal{E}_{j,d})\\&~~~~\lesssim
\mathbb{P}\bigg\{d^2K\max_{ k\in[n_{i,j}]}\max_{t_{i,j,\ell}\in S_{i,j,k}}(\Delta M_{j,\cdot})_{j,i}^{(k,\ell)}\cdot (\Delta M_{i,\cdot})_{i,j}\geq 2n_{i,j}v\bigg\}+\exp(-Cd^{\iota})\notag
\end{align}
for any $v> 0$.
 By  Bonferroni inequality, we have that
\begin{align}\label{eq:tijk}
\mathbb{P}\bigg\{\max_{ k\in[n_{i,j}]}\max_{t_{i,j,\ell}\in S_{i,j,k}}(\Delta M_{j,\cdot})_{j,i}^{(k,\ell)}\geq v\bigg\}\lesssim&~ n_*\max_{ k\in[n_{i,j}]}\max_{t_{i,j,\ell}\in S_{i,j,k}}\mathbb{P}\{(\Delta M_{j,\cdot})_{j,i}^{(k,\ell)}\geq v\}
\end{align}
for any $v>0$. Recall $\xi=\max_{i,j\in[p]}\max_{ k\in[n_{i,j}]}\max_{t_{i,j,\ell}\in S_{i,j,k}}|t_{i,j,\ell}-t_{i,j,k}|\asymp n_*^{-1}K$. By Assumption \ref{as:jump1}(ii), $\max_{ k\in[n_{i,j}]}\max_{t_{i,j,\ell}\in S_{i,j,k}}\mathbb{E}\{(\Delta M_{j,\cdot})_{j,i}^{(k,\ell)}\}\lesssim \lambda_*n_*^{-1}K$. If $\lambda_*n_*^{-1}K=o(1)$,   applying Markov's inequality, Assumption \ref{as:jump1} implies that
\begin{align*}
\mathbb{P}\{(\Delta M_{j,\cdot})_{j,i}^{(k,\ell)}\geq v\}=&~\mathbb{P}[(\Delta M_{j,\cdot})_{j,i}^{(k,\ell)}-\mathbb{E}\{(\Delta M_{j,\cdot})_{j,i}^{(k,\ell)}\}\geq v-\mathbb{E}\{(\Delta M_{j,\cdot})_{j,i}^{(k,\ell)}\}]
\\\lesssim&~
\mathbb{P}[(\Delta M_{j,\cdot})_{j,i}^{(k,\ell)}-\mathbb{E}\{(\Delta M_{j,\cdot})_{j,i}^{(k,\ell)}\}\geq v/2]
\\\lesssim&~
\exp(-Cv)\exp(C\lambda_*n_*^{-1}K)\lesssim \exp(-Cv)
\end{align*}
for any $v\gg \lambda_*n_*^{-1}K$. By \eqref{eq:tijk}, we have
$\mathbb{P}\{\max_{ k\in[n_{i,j}]}\max_{t_{i,j,\ell}\in S_{i,j,k}}(\Delta M_{j,\cdot})_{j,i}^{(k,\ell)}\geq v\}\lesssim n_*\exp(-Cv)$ for any $v\gg \lambda_*n_*^{-1}K$.
Together with \eqref{eq:mijl}, for any  $\delta\in(0,1)$, it holds that
\begin{align*}
&\mathbb{P}\bigg\{d^2K\max_{ k\in[n_{i,j}]}\max_{t_{i,j,\ell}\in S_{i,j,k}}(\Delta M_{j,\cdot})_{j,i}^{(k,\ell)}\cdot (\Delta M_{i,\cdot})_{i,j}\geq 2n_{i,j}v\bigg\}
\\&~~~~~\lesssim
\mathbb{P}\bigg\{\max_{ k\in[n_{i,j}]}\max_{t_{i,j,\ell}\in S_{i,j,k}}(\Delta M_{j,\cdot})_{j,i}^{(k,\ell)}\geq \xi^{-\delta}\}+\mathbb{P}\{d^2K (\Delta M_{i,\cdot})_{i,j}\geq 2n_{i,j}\xi^{\delta}v\}
\\&~~~~~\lesssim
n_*\exp(-C\xi^{-\delta})+\exp(-Cd^{-2}\xi^{\delta-1}v)
\end{align*}
for any $v\gg\xi^{1-\delta}d^2\lambda_*$. Together with \eqref{eq:w(i,j)} and \eqref{eq:qij},
$
\max_{ i,j\in[p]}\mathbb{P}(|\varpi_{i,j}|\geq v)\lesssim n_*\exp(-C\xi^{-\delta})+\exp(-Cd^{-2}\xi^{\delta-1}v)+\exp(-Cd^{\iota})
$
for any $v\gg\xi^{1-\delta}d^2\lambda_*$. In order to make $\max_{ i,j\in[p]}|\varpi_{i,j}|=O_{\p}(\xi^{1/2}\log^{1/2}p)$, it suffices to require  $\log p=o[\min\{\xi^{-\delta},d^{-4}\xi^{2\delta-1},d^{\iota}\}]$, $\lambda_*^2\xi^{1-2\delta}d^4(\log p)^{-1}=o(1)$ and $\log n_*=o(\xi^{-\delta})$. To make $p$ diverge as fast as possible, we select
$\delta=\iota/(3\iota+4)$ and $d^{\iota}\asymp \xi^{-\delta}$. Then if $K(\log n_*)^{(3\iota+4)/\iota}=o(n_*)$,
$\max_{ i,j\in[p]}\mathbb{P}(|\varpi_{i,j}|\geq v)\lesssim \exp\{-C(n_*K^{-1})^{(2\iota+2)/(3\iota+4)}v\}+ \exp\{-C(n_*K^{-1})^{\iota/(3\iota+4)}\}$ for any $v\gg (n_*^{-1}K)^{(2\iota+2)/(3\iota+4)}\lambda_*$.

Analogous to \eqref{eq:D2}, by Assumptions \ref{as:jumpsize} and \ref{as:jump1},  it holds that
\begin{align*}
	\max_{ i,j\in[p]}\mathbb{E}(|\varpi_{i,j}|^m)\lesssim \max_{ i,j\in[p]}\max_{ k\in[n_{i,j}]}\max_{t_{i,j,\ell}\in S_{i,j,k}}\mathbb{E}\bigg\{\bigg(\int_{t_{i,j,k}\wedge t_{i,j,\ell}}^{t_{i,j,k}\vee t_{i,j,\ell}}|J_{i,s}|\,{\rm d}M_{i,s}\bigg)^{2m}\bigg\}\lesssim 1
\end{align*}
for any fixed positive integer $m$.
We complete the proof of Theorem \ref{thwij}.
$\hfill\Box$

\subsection{Proof of Theorem \ref{eq:jjj}}
Write $\aleph=(n_*^{-1}{K\log p})^{1/2}$. Recall $\hat{\sigma}_{u,i,j}^{{\rm jump}}=\hat{\sigma}_{u,i,j}-\varpi_{i,j}$. By Theorem \ref{thwij} and the proof of Section \ref{sec:prothe6}, for some sufficiently large constant $\alpha_*>0$, if $\lambda_*^2(n_*^{-1}K)^{\iota/(3\iota+4)}(\log p)^{-1}=o(1)$, we have
$
\mathbb{P}(\max_{i,j\in[p]}|\hat{\sigma}_{u,i,j}-\sigma_{u,i,j}|\geq \alpha_*\aleph)
\lesssim p^2\max_{i,j\in[p]}\mathbb{P}(|\hat{\sigma}_{u,i,j}^{{\rm jump}}-\sigma_{u,i,j}|> \alpha_*\aleph/2)+p^2	\max_{ i,j\in[p]}\mathbb{P}(|\varpi_{i,j}|\geq \alpha_*\aleph/2)=
o(\aleph^4)
$
provided that $\log p=o[\min\{(n_*L_n^{-2}K)^{\varphi/(3\varphi+2)},(n_*K^{-1})^{\chi^*}\}]$, $K^{-\varphi}L_n^{\varphi}\log\{n_*(K\log p)^{-1}\}=o(1)$ and $K\gtrsim L_n$, where $\chi^*=\min\{\iota\gamma/(2\iota+2\gamma+\iota\gamma),\iota/(3\iota+4),\gamma/(\gamma+4)\}$. Hence, Part (i) holds. For any $C>0$, by Markov's inequality, we have
$
\mathbb{P}(\|\widehat{\bSigma}_u^{{\rm thre}}-\bSigma_u\|_2\geq C c_p\aleph^{1-q})\leq
C^{-2}c_p^{-2}\aleph^{2(q-1)}\mathbb{E}(\|\widehat{\bSigma}_u^{{\rm thre}}-\bSigma_u\|_2^2)$.
In the sequel, we will show $\mathbb{E}(\|\widehat{\bSigma}_u^{{\rm thre}}-\bSigma_u\|_2^2)\lesssim c_p^2\aleph^{2(1-q)}$. Based on this result, we know Part (ii) holds. For each $i,j\in[p]$, we  define the event
$
A_{i,j}=\{|\hat{\sigma}_{u,i,j}^{\rm  thre}-\sigma_{u,i,j}|\leq 4\min(|\sigma_{u,i,j}|,\alpha\aleph)\}
$
with some $\alpha>0$, and
$
d_{i,j}=(\hat{\sigma}_{u,i,j}^{\rm thre}-\sigma_{u,i,j}){I}(A_{i,j}^c)$. As shown in \eqref{eq:speupp1} and \eqref{eq:Edij}, $
\mathbb{E}(\|\widehat{\bSigma}_{u}^{\rm  thre}-\bSigma_{u}\|_2^2)\lesssim \mathbb{E}\{(\max_{i\in[p]}\sum_{j=1}^p|d_{i,j}|)^2\}+c_p^2\aleph^{2(1-q)}$ and
\begin{align*}
\mathbb{E}\bigg\{\bigg(\max_{i\in[p]}\sum_{j=1}^p|d_{i,j}|\bigg)^2\bigg\}\leq&~\underbrace{p\sum_{i,j=1}^p\mathbb{E}\big(|\hat{\sigma}_{u,i,j}^{\rm  thre}-\sigma_{u,i,j}|^2{I}[A_{i,j}^c\cap\{\hat{\sigma}_{u,i,j}^{\rm  thre}=0\}]\big)}_{\textrm{I}}\\
&+\underbrace{p\sum_{i,j=1}^p\mathbb{E}\big(|\hat{\sigma}_{u,i,j}^{\rm  thre}-\sigma_{u,i,j}|^2{I}[A_{i,j}^c\cap\{\hat{\sigma}_{u,i,j}^{\rm  thre}=\hat{\sigma}_{u,i,j}^{\rm }\}]\big)}_{\textrm{II}}\,.
\end{align*}
Recall $\hat{\sigma}_{u,i,j}^{\rm  thre}=\hat{\sigma}_{u,i,j}{I}(|\hat{\sigma}_{u,i,j}|\geq \beta\aleph)$. By \eqref{eq:th3I},
$
{\rm I}
\leq p\sum_{i,j=1}^p\sigma_{u,i,j}^2\mathbb{P}\{|\hat{\sigma}_{u,i,j}-\sigma_{u,i,j}|\geq(4\alpha-{\beta})\aleph\}$.
Notice that
$
\max_{i,j\in[p]}\mathbb{E}\{|\hat{\sigma}_{u,i,j}-\sigma_{u,i,j}|^4\}
\lesssim n_*^{1/2}$. Identical to \eqref{eq:IIthm3},
$
{\rm II}\lesssim pn_*^{1/4}\sum_{i,j=1}^p\{\mathbb{P}(|\hat{\sigma}_{u,i,j}-\sigma_{u,i,j}|>4\alpha\aleph)\}^{1/2}+pn_*^{1/4}\sum_{i,j=1}^p[\mathbb{P}\{|\hat{\sigma}_{u,i,j}-\sigma_{u,i,j}|>({\beta}-\alpha)\aleph\}]^{1/2}$. Selecting $\alpha={\beta}/2$ for some sufficiently large $\beta>0$, we have
$
{\rm I}+{\rm II}\lesssim o(\aleph
^2)$ provided that  $\lambda_*^2(n_*^{-1}K)^{\iota/(3\iota+4)}(\log p)^{-1}=o(1)$, $\log p=o[\min\{(n_*L_n^{-2}K)^{\varphi/(3\varphi+2)},(n_*K^{-1})^{\chi^*}\}]$,  $K^{-\varphi}L_n^{\varphi}\log\{n_*(K\log p)^{-1}\}=o(1)$ and $K\gtrsim L_n$, where $\chi^*=\min\{\iota\gamma/(2\iota+2\gamma+\iota\gamma),\iota/(3\iota+4),\gamma/(\gamma+4)\}$. Thus
$
\mathbb{E}(\|\widehat{\bf\Sigma}_{u}^{\rm  thre}-{\bf\Sigma}_{u}\|_2^2)\lesssim c_p^2\aleph^{2(1-q)}
$.
We complete the proof of Theorem \ref{eq:jjj}. $\hfill\Box$

\subsection{Proof of Lemma \ref{bern.inq}}\label{sec:pfberninq}

Write $S_k=\sum_{t=1}^{k}z_t$. Note that $\max_{t\in[\tilde{n}]}{\rm{Var}}(z_t)<\infty$. We will apply the Fuk-Nagaev inequality \cite[Theorem 6.2]{Rio_2017} to bound the tail probability of $\max_{k\in[\tilde{n}]}|S_k|$.
Define
$
\alpha^{-1}(u)=\sum_{k=1}^{\infty}I\{u<\alpha(k)\}$.
Since $\alpha(k)\leq  a_1\exp(-a_2\tilde{L}_{\tilde{n}}^{-\varphi}|k-m|_+^{\varphi})$, then $
\alpha^{-1}(u)\leq m+a_2^{-1/\varphi}\tilde{L}_{\tilde{n}}\log^{1/\varphi}(a_1u^{-1})$ for any $u>0$.
Define $Q(u)=\sup_{t\in[\tilde{n}]}Q_t(u)$ with $Q_t(u)=\inf\{x>0:\mathbb{P}(|z_t|>x)\leq u\}$. Since $\mathbb{P}(|z_t|>x)\leq b_1\exp(-b_2x^{r})$, then $Q(u)\leq b_2^{-1/r}\log^{1/r}(b_1u^{-1})$. Define $R(u)=\alpha^{-1}(u)Q(u)$. We have
$
R(u)\leq
c_1\tilde{L}_{\tilde{n}}\log^{1/r_*}(c_2u^{-1})+c_1m\log^{1/r}(c_2u^{-1})$ with $c_1=b_2^{-1/r}\max(a_2^{-1/\varphi},1)$, $c_2=\max(a_1,b_1)$ and $r_*=r\varphi/(r+\varphi)$.
Since  $R(u)$ is a right-continuous and non-increasing function, then its inverse function
$
\notag H(x)=R^{-1}(x)=\inf\{u:R(u)\leq x\}
\leq\inf\{u:c_1\tilde{L}_{\tilde{n}}\log^{1/r_*}(c_2u^{-1})+c_1m\log^{1/r}(c_2u^{-1})\leq x\}
\leq
\tilde{c}_1\exp(-\tilde{c}_2\tilde{L}_{\tilde{n}}^{-r_*}x^{r_*})+\tilde{c}_1\exp(-\tilde{c}_2m^{-r}x^{r})$
for any $x>0$ with $\tilde{c}_1=c_2$ and $\tilde{c}_2=\min\{(2c_1)^{-r_*},(2c_1)^{-r}\}
$. Write $\tilde{u}=\tilde{c}_1\exp(-\tilde{c}_2\tilde{L}_{\tilde{n}}^{-r_*}x^{r_*})+\tilde{c}_1\exp(-\tilde{c}_2m^{-r}x^{r})$. Therefore,
\begin{align*}
\int_{0}^{H(x)}Q(u)\,\mathrm{d}u\leq&\,\int_{0}^{\tilde{u}}b_2^{-1/r}\log^{1/r}(b_1u^{-1})\,\mathrm{d}u\\
\lesssim&\, \int_{\log(b_1\tilde{u}^{-1})}^{\infty}y^{1/r}e^{-y}\,\mathrm{d}y\lesssim
b_1^{-1}\tilde{u}\log^{1/r}(b_1\tilde{u}^{-1})\lesssim \tilde{u}^{1/2}\{\tilde{u}^{r/2}\log(b_1\tilde{u}^{-1})\}^{1/r}\,.
\end{align*}
As $x\to +\infty$, we have $\tilde{u}\to 0^{+}$, which implies $\tilde{u}^{r/2}\log(b_1\tilde{u}^{-1})\to 0$. Hence there exists a universal constant $\varepsilon>0$ such that $\tilde{u}^{r/2}\log(b_1\tilde{u}^{-1})<\varepsilon$ for any $\tilde{u}\in(0,2\tilde{c}_1]$. By the definition of $\tilde{u}$, we have
$
\int_{0}^{H(x)}Q(u)\,\mathrm{d}u\lesssim
\exp(-C\tilde{L}_{\tilde{n}}^{-r_*}x^{r_*})+\exp(-Cm^{-r}x^{r})$
for any $x>0$. Recall $s_{\tilde{n}}^2=\sum_{t_1,t_2=1}^{\tilde{n}}|\mathrm{Cov}(z_{t_1},z_{t_2})|$. By the Fuk-Nagaev inequality,
$
\mathbb{P}(\max_{k\in[\tilde{n}]}|S_k|\geq 4\tilde{n}\lambda)\leq
4(1+\tilde{n}^2\rho^{-1}s_{\tilde{n}}^{-2}\lambda^2)^{-\rho/2}+4\lambda^{-1}\int_{0}^{H(\tilde{n}\lambda/\rho)}Q(u)\,\mathrm{d}u\lesssim
(1+\tilde{n}^2\rho^{-1}s_{\tilde{n}}^{-2}\lambda^2)^{-\rho/2}+\lambda^{-1}\{\exp(-C\tilde{n}^{r_*}\lambda^{r_*}\rho^{-r_*}\tilde{L}_{\tilde{n}}^{-r_*})+\exp(-Cm^{-r}\tilde{n}^{r}\lambda^{r}\rho^{-r})\}$
for any $\lambda>0$ and $\rho\geq 1$. We complete the proof of Lemma \ref{bern.inq}. $\hfill\Box$

\subsection{Proof of Lemma \ref{eq:mixingbound}}\label{se:pflemma2}

We first consider the case with $\varphi<\infty$. Define $\Lambda_{\tilde{n}}(\alpha,u)=\max\{1,\max_{s\in[\tilde{n}]}\sum_{t=s}^{\tilde{n}}\alpha^{1/u}(t-s)\}$. Notice that $\alpha(k)\leq a_1\exp\{-a_2(\tilde{L}_{\tilde{n}}^{-1}k)^{\varphi}\} $ for any integer $k\geq 1$. For any $k\geq 2$,  it holds that
\begin{align*}
\Lambda_{\tilde{n}}\{\alpha,2(k-1)\}\leq&~ 1+\max_{s\in[\tilde{n}]}\sum_{t=s}^{\tilde{n}}\{\alpha(t-s)\}^{1/(2k-2)}\leq1+\sum_{m=0}^{\tilde{n}}\{\alpha(m)\}^{1/(2k-2)}\\
\leq &~ 1+a_1+a_1\sum_{m=1}^{\tilde{n}}\int_{m-1}^{m}\exp\{-a_2(2k-2)^{-1}\tilde{L}_{\tilde{n}}^{-\varphi} x^{\varphi}\}\,\mathrm{d}x	\\
\leq&~ 1+a_1+a_1\int_{0}^{\infty}\exp\{-a_2(2k-2)^{-1}\tilde{L}_{\tilde{n}}^{-\varphi}x^{\varphi}\}\,\mathrm{d}x\leq C_*(k-1)^{1/\varphi}\tilde{L}_{\tilde{n}}
\end{align*}
for some constant $C_*>0$ independent of $k$. Due to $k^k\leq k!e^k$ for any integer $k\geq1$, we then have
$
\Lambda_{\tilde{n}}^{k-1}\{\alpha,2(k-1)\}\leq C_*^{k-1}\{(k-1)^{k-1}\}^{1/\varphi}{\tilde{L}}_{\tilde{n}}^{k-1}\leq (C_*e^{2/\varphi}{\tilde{L}}_{\tilde{n}})(C_*e^{1/\varphi}{\tilde{L}}_{\tilde{n}})^{k-2}(k!)^{1/\varphi}
$.
Let $\Gamma_k(x)$ be the $k$-th order cumulant of the random variable $x$. By Theorem 4.17 of \cite{SS_1991} with $\delta=1$, we have
$
|\Gamma_k(S_{\tilde{n}})|\leq (k!)^{2+\tilde{r}+1/\varphi}(\bar{C}{\tilde{L}}_{\tilde{n}}H_{\tilde{n}}^2)(\bar{C}{\tilde{L}}_{\tilde{n}}H_{\tilde{n}})^{k-2}\tilde{n}$ with $\bar{C}=2^{7+2\tilde{r}}C_*e^{2/\varphi}$, which implies $|\Gamma_k(S_{\tilde{n}}/\tilde{n})|=\tilde{n}^{-k}|\Gamma_k(S_{\tilde{n}})|\leq (k!)^{2+\tilde{r}+1/\varphi}(\bar{C}{\tilde{L}}_{\tilde{n}}H_{\tilde{n}}^2\tilde{n}^{-1})(\bar{C}{\tilde{L}}_{\tilde{n}}H_{\tilde{n}}\tilde{n}^{-1})^{k-2}$.
By Lemma 2.4 of \cite{SS_1991},  it holds that
\begin{align*}
\mathbb{P}(|S_{\tilde{n}}|\geq \tilde{n}x)\lesssim \exp(-C\tilde{n}{\tilde{L}}_{\tilde{n}}^{-1}H_{\tilde{n}}^{-2}x^2)+\exp\{-C(\tilde{n}{\tilde{L}}_{\tilde{n}}^{-1}H_{\tilde{n}}^{-1}x)^{1/(1+\check{r})}\}
\end{align*}
for any $x>0$, where $\check{r}=1+\tilde{r}+\varphi^{-1}$. Analogously, when $\varphi =\infty$, we can also show
\begin{align*}
\mathbb{P}(|S_{\tilde{n}}|\geq \tilde{n} x)\lesssim \exp(-C\tilde{n}{\tilde{L}}_{\tilde{n}}^{-1}H_{\tilde{n}}^{-2}x^2)+\exp\{-C(\tilde{n}{\tilde{L}}_{\tilde{n}}^{-1}H_{\tilde{n}}^{-1}x)^{1/(1+\check{r})}\}
\end{align*}
for any $x>0$, where $ \check{r}=1+\tilde{r}$. We complete the proof of Lemma \ref{eq:mixingbound}.  $\hfill\Box$

\newpage
\normalsize

\begin{center}
	{\bf  \Large
		Supplementary Material for ``Optimal Covariance Matrix Estimation for  High-dimensional  Noise in  High-frequency Data" by Jinyuan Chang, Qiao Hu, Cheng Liu and Cheng Yong Tang}  \\
\end{center}

\setcounter{page}{1}

\renewcommand{\theequation}{S.\arabic{section}.\arabic{equation}}
	
	\renewcommand{\thesection}{S\arabic{section}}
	
	\setcounter{equation}{0}
	
	\setcounter{section}{0}
	
	
\section{Simulation results for independent measurement errors $\{\bU_{\tilde{t}_k}\}$}	

When the measurement errors $\{\bU_{\tilde{t}_k}\}$ are independent and identically distributed, we can compare our proposed estimators $\widehat{\bSigma}_u^{{\rm thre}}$ and $\widehat{\bSigma}_u^{{\rm bc,thre}}$ given in \eqref{eq:praest} with $\widehat{\bSigma}_{u}^{*,{\rm thre}}$, which is developed based on thresholding the estimator proposed in \cite{Christensen2013}:
\begin{align}\label{eq:ac1est}
\widehat{\bSigma}_{u}^{*,{\rm thre}}&=(\hat{\sigma}_{u,i, j}^{*,\rm thre})_{p\times p}=\big\{\hat{\sigma}_{u,i, j}^*{I}\big(|\hat{\sigma}^*_{u, i, j}|\geq\varpi_{i,j}^*\big)\big\}_{p\times p}\,,
\end{align}
where
\[
\begin{split}
\hat{\sigma}_{u,i, j}^*=-\frac{1}{2(n_{i, j}-2)}\sum_{k=2}^{n_{i, j}-1}&\big\{(Y_{i,t_{i, j, k+1}}-Y_{i, t_{i, j, k}})(Y_{j,t_{i, j, k}}-Y_{j, t_{i, j, k-1}})\\
&~~~~~~~~~~+(Y_{j,t_{i, j, k+1}}-Y_{j, t_{i, j, k}})(Y_{i,t_{i, j, k}}-Y_{i, t_{i, j, k-1}})\big\}\,.
\end{split}
\]
To determine the thresholding level $\varpi_{i,j}^*$, we defined $\zeta_{i,j,k}^*=-2^{-1}\{(Y_{i,t_{i, j, k+1}}-Y_{i, t_{i, j, k}})(Y_{j,t_{i, j, k}}-Y_{j, t_{i, j, k-1}})+(Y_{j,t_{i, j, k+1}}-Y_{j, t_{i, j, k}})(Y_{i,t_{i, j, k}}-Y_{i, t_{i, j, k-1}})\}$ for each $k=2,\ldots,n_{i,j}-1$, and used the same procedure  as that for choosing the threshold level in  our proposed estimators $\widehat{\bSigma}_{u}^{{\rm thre}}$ and $\widehat{\bSigma}_{u}^{{\rm bc, thre}}$, 
	upon replacing $\zeta_{i,j,k}$ by $\zeta_{i,j,k}^*$. We kept the same data generating procedure for $\bX_t$ and generated the measurement errors independently from $N(\bzero,\sigma_e^2\bR)$ with $\sigma_e^2=0.005^2$ and $0.001^2$. The settings for $\bR$ are identical to those in the main paper. Table \ref{Table-IndepN} reports the averages of relative estimation errors for different estimators based on 1000 repetitions.
 We can find that our estimator $\widehat{\bSigma}_{u}^{{\rm bc, thre}}$ outperforms  $\widehat{\bSigma}_{u}^{*,{\rm thre}}$ in all cases, and $\widehat{\bSigma}_{u}^{{\rm thre}}$ outperforms  $\widehat{\bSigma}_{u}^{*,{\rm thre}}$ for the cases with high signal-to-noise ratio $252\sigma_e^2/\bar{\sigma}^2=0.063$.
\begin{table}[htbp]
	\tiny
	\centering
	\caption{Averages of the relative estimation errors ($\times100$) for the proposed estimators and $\widehat{\bSigma}_{u}^{*,{\rm thre}}$ with i.i.d. measurement errors when jumps exist and do not exist (in parentheses) based on 1000 repetitions.}
	\begin{tabular}{ccccccccccccc}
		\hline
		\hline
		\multicolumn{3}{c}{Synchronous Data}& &  \multicolumn{3}{c}{Model 1}  &  \multicolumn{3}{c}{Model 2}  &  \multicolumn{3}{c}{Model 3}\\
		$252\sigma_e^2/\bar{\sigma}^2$ &   $p $    & Estimators  & $\Delta_K$ & $\Delta=3$     & $\Delta=2$     & $\Delta=1$     & $\Delta=3$     & $\Delta=2$     & $\Delta=1$     & $\Delta=3$     & $\Delta=2$     & $\Delta=1$ \\
		\hline

    0.063 & 50    & $\widehat{\bSigma}_{u}^{{\rm thre}}$       & 1     & 3.5(3.5) & 2.9(2.9) & 2.0(2.0)  & 4.0(4.1) & 3.3(3.4) & 2.3(2.4) & 6.0(6.0)  & 5.0(5.0)  & 3.6(3.6) \\
          &       &       & 2     & 3.4(3.4) & 2.8(2.8) & 2.0(1.9) & 3.9(4.0) & 3.2(3.3) & 2.2(2.3) & 5.8(5.8) & 4.8(4.8) & 3.4(3.4) \\
          &       &       & 3     & 3.4(3.4) & 2.7(2.7) & 1.9(1.9) & 3.9(4.0) & 3.2(3.2) & 2.2(2.3) & 5.6(5.6) & 4.6(4.6) & 3.3(3.3) \\
          &       &   $\widehat{\bSigma}_{u}^{{\rm bc, thre}}$     & 1     & 3.6(3.6) & 2.9(2.9) & 2.1(2.1) & 4.1(4.2) & 3.3(3.4) & 2.4(2.4) & 6.1(6.1) & 5.0(5.0)  & 3.6(3.6) \\
          &       &       & 2     & 3.5(3.5) & 2.8(2.8) & 2.0(2.0)  & 4.0(4.1) & 3.2(3.3) & 2.3(2.3) & 5.8(5.8) & 4.8(4.8) & 3.4(3.4) \\
          &       &       & 3     & 3.4(3.4) & 2.8(2.8) & 2.0(1.9) & 4.0(4.1) & 3.2(3.3) & 2.3(2.3) & 5.6(5.6) & 4.6(4.6) & 3.4(3.3) \\
          &       &  $\widehat{\bSigma}_{u}^{*, {\rm thre}}$      &       & 5.9(5.8) & 4.8(4.8) & 3.4(3.4) & 7.8(7.5) & 5.6(5.6) & 3.8(4.0) & 9.9(9.9) & 8.2(8.2) & 5.9(5.8) \\
          & 100   &  $\widehat{\bSigma}_{u}^{{\rm thre}}$      & 1     & 3.9(3.9) & 3.2(3.2) & 2.2(2.2) & 4.2(4.2) & 3.4(3.5) & 2.4(2.4) & 6.8(6.8) & 5.6(5.7) & 4.0(4.0) \\
          &       &       & 2     & 3.8(3.8) & 3.1(3.1) & 2.2(2.2) & 4.1(4.1) & 3.3(3.4) & 2.3(2.4) & 6.5(6.5) & 5.4(5.4) & 3.8(3.9) \\
          &       &       & 3     & 3.7(3.7) & 3.0(3.0)  & 2.1(2.1) & 4.0(4.1) & 3.3(3.3) & 2.3(2.3) & 6.2(6.2) & 5.2(5.2) & 3.7(3.8) \\
          &       &   $\widehat{\bSigma}_{u}^{{\rm bc, thre}}$     & 1     & 3.9(3.9) & 3.2(3.2) & 2.3(2.3) & 4.2(4.3) & 3.5(3.5) & 2.4(2.5) & 6.8(6.8) & 5.6(5.7) & 4.0(4.0) \\
          &       &       & 2     & 3.8(3.8) & 3.1(3.1) & 2.2(2.2) & 4.1(4.2) & 3.4(3.4) & 2.4(2.4) & 6.5(6.5) & 5.4(5.5) & 3.8(3.9) \\
          &       &       & 3     & 3.8(3.8) & 3.1(3.1) & 2.2(2.2) & 4.1(4.1) & 3.3(3.4) & 2.3(2.4) & 6.3(6.3) & 5.2(5.3) & 3.7(3.8) \\
          &       &   $\widehat{\bSigma}_{u}^{*, {\rm thre}}$     &       & 6.5(6.5) & 5.3(5.3) & 3.7(3.7) & 9.2(8.8) & 6.0(5.9) & 4.0(4.1) & 11.1(11.3) & 9.3(9.3) & 6.6(6.6) \\
          & 200   &  $\widehat{\bSigma}_{u}^{{\rm thre}}$      & 1     & 4.2(4.2) & 3.4(3.5) & 2.4(2.4) & 4.0(4.0)  & 3.3(3.2) & 2.3(2.3) & 7.5(7.5) & 6.2(6.2) & 4.4(4.4) \\
          &       &       & 2     & 4.1(4.1) & 3.3(3.3) & 2.3(2.4) & 3.9(3.9) & 3.2(3.1) & 2.2(2.2) & 7.2(7.1) & 5.9(6.0) & 4.2(4.2) \\
          &       &       & 3     & 4.0(4.0)  & 3.3(3.3) & 2.3(2.3) & 3.8(3.8) & 3.1(3.1) & 2.2(2.2) & 6.9(6.9) & 5.7(5.8) & 4.1(4.1) \\
          &       &   $\widehat{\bSigma}_{u}^{{\rm bc, thre}}$     & 1     & 4.3(4.3) & 3.5(3.5) & 2.5(2.5) & 4.0(4.0)  & 3.3(3.3) & 2.3(2.3) & 7.5(7.5) & 6.2(6.3) & 4.4(4.4) \\
          &       &       & 2     & 4.1(4.1) & 3.4(3.4) & 2.4(2.4) & 3.9(3.9) & 3.2(3.2) & 2.3(2.2) & 7.2(7.2) & 6.0(6.0)  & 4.2(4.3) \\
          &       &       & 3     & 4.1(4.1) & 3.3(3.3) & 2.3(2.3) & 3.9(3.8) & 3.2(3.1) & 2.2(2.2) & 6.9(6.9) & 5.8(5.8) & 4.1(4.1) \\
          &       &  $\widehat{\bSigma}_{u}^{*, {\rm thre}}$      &       & 7.1(7.1) & 5.7(5.8) & 4.0(4.1) & 10.0(9.2) & 6.3(6.0) & 3.8(3.8) & 12.5(12.4) & 10.3(10.4) & 7.3(7.3) \\
 \hline
    0.00252 & 50    & $\widehat{\bSigma}_{u}^{{\rm thre}}$       & 1     & 6.4(6.7) & 4.3(4.5) & 2.4(2.4) & 9.9(10.2) & 6.6(6.5) & 3.5(3.4) & 7.0(7.2) & 5.1(5.2) & 3.6(3.6) \\
          &       &       & 2     & 8.7(9.0) & 5.7(6.0) & 2.9(3.0) & 14.8(15.5) & 9.2(9.5) & 4.4(4.4) & 9.0(9.4) & 6.2(6.4) & 3.6(3.7) \\
          &       &       & 3     & 11.5(11.7) & 7.4(7.6) & 3.6(3.7) & 21.1(21.0) & 12.9(13.3) & 5.7(5.9) & 11.3(11.8) & 7.6(7.9) & 4.1(4.2) \\
          &       & $\widehat{\bSigma}_{u}^{{\rm bc, thre}}$       & 1     & 3.6(3.6) & 2.9(2.9) & 2.1(2.1) & 4.2(4.3) & 3.4(3.5) & 2.4(2.4) & 6.1(6.1) & 5.1(5.1) & 3.6(3.6) \\
          &       &       & 2     & 3.5(3.5) & 2.8(2.8) & 2.0(2.0)  & 4.2(4.2) & 3.3(3.4) & 2.3(2.4) & 5.9(5.9) & 4.9(4.9) & 3.5(3.5) \\
          &       &       & 3     & 3.5(3.5) & 2.8(2.8) & 2.0(2.0)  & 4.3(4.3) & 3.4(3.4) & 2.3(2.4) & 5.8(5.8) & 4.7(4.7) & 3.4(3.4) \\
          &       & $\widehat{\bSigma}_{u}^{*, {\rm thre}}$       &       & 5.9(5.9) & 4.8(4.8) & 3.4(3.4) & 7.9(7.7) & 5.6(5.7) & 3.9(4.0) & 10.0(10.0) & 8.3(8.3) & 5.9(5.8) \\
          & 100   & $\widehat{\bSigma}_{u}^{{\rm thre}}$       & 1     & 7.9(7.4) & 5.3(4.9) & 2.7(2.6) & 10.1(9.6) & 6.6(6.2) & 3.5(3.3) & 8.5(8.4) & 5.8(5.9) & 4.0(4.0) \\
          &       &       & 2     & 10.7(10.4) & 7.1(6.7) & 3.5(3.3) & 15.0(14.4) & 9.4(8.9) & 4.4(4.1) & 11.3(11.1) & 7.5(7.4) & 4.1(4.1) \\
          &       &       & 3     & 14(13.6) & 9.0(8.8) & 4.4(4.2) & 20.9(19.8) & 12.9(12.4) & 5.8(5.5) & 14.0(13.9) & 9.3(9.1) & 4.9(4.9) \\
          &       & $\widehat{\bSigma}_{u}^{{\rm bc, thre}}$       & 1     & 4.0(4.0)  & 3.2(3.2) & 2.3(2.3) & 4.4(4.4) & 3.5(3.6) & 2.5(2.5) & 6.9(6.9) & 5.7(5.8) & 4.0(4.1) \\
          &       &       & 2     & 3.8(3.8) & 3.1(3.1) & 2.2(2.2) & 4.3(4.3) & 3.5(3.5) & 2.4(2.4) & 6.6(6.7) & 5.5(5.5) & 3.9(3.9) \\
          &       &       & 3     & 3.8(3.8) & 3.1(3.1) & 2.2(2.2) & 4.4(4.4) & 3.5(3.5) & 2.4(2.4) & 6.4(6.5) & 5.3(5.4) & 3.8(3.8) \\
          &       & $\widehat{\bSigma}_{u}^{*, {\rm thre}}$       &       & 6.5(6.5) & 5.3(5.3) & 3.7(3.7) & 9.4(9.0) & 6.1(6.0) & 4.1(4.1) & 11.3(11.4) & 9.4(9.4) & 6.6(6.6) \\
          & 200   &  $\widehat{\bSigma}_{u}^{{\rm thre}}$      & 1     & 8.0(7.6) & 5.3(5.1) & 2.8(2.7) & 7.8(7.4) & 5.2(5.0) & 2.9(2.9) & 8.8(8.8) & 6.2(6.2) & 4.3(4.3) \\
          &       &       & 2     & 10.7(10.4) & 7.1(6.8) & 3.6(3.4) & 11.2(10.7) & 7.0(6.8) & 3.5(3.3) & 11.5(11.5) & 7.6(7.5) & 4.3(4.3) \\
          &       &       & 3     & 14(13.7) & 9.0(8.8) & 4.4(4.2) & 15.4(14.8) & 9.6(9.2) & 4.4(4.2) & 14.2(14.3) & 9.4(9.3) & 4.9(5.0) \\
          &       &  $\widehat{\bSigma}_{u}^{{\rm bc, thre}}$      & 1     & 4.3(4.3) & 3.5(3.5) & 2.5(2.5) & 4.2(4.1) & 3.4(3.3) & 2.4(2.3) & 7.7(7.7) & 6.3(6.3) & 4.4(4.5) \\
          &       &       & 2     & 4.2(4.2) & 3.4(3.4) & 2.4(2.4) & 4.1(4.1) & 3.3(3.2) & 2.3(2.3) & 7.4(7.4) & 6.1(6.1) & 4.3(4.3) \\
          &       &       & 3     & 4.1(4.1) & 3.4(3.4) & 2.4(2.4) & 4.1(4.1) & 3.3(3.2) & 2.3(2.2) & 7.1(7.1) & 5.9(5.9) & 4.1(4.2) \\
          &       &  $\widehat{\bSigma}_{u}^{*, {\rm thre}}$      &       & 7.1(7.1) & 5.8(5.8) & 4.0(4.1) & 10.3(9.5) & 6.4(6.1) & 3.9(3.8) & 12.6(12.6) & 10.4(10.4) & 7.3(7.3) \\

		\hline
		\multicolumn{3}{c}{Asynchronous Data} & & \multicolumn{3}{c}{Model 1}  &  \multicolumn{3}{c}{Model 2}  &  \multicolumn{3}{c}{Model 3}\\
		$252\sigma_e^2/\bar{\sigma}^2$ &   $p $    & Estimators  & $\Delta_K$ & $\lambda=3$     & $\lambda=2$     & $\lambda=1$     & $\lambda=3$     & $\lambda=2$     & $\lambda=1$     & $\lambda=3$     & $\lambda=2$     & $\lambda=1$\\
		\hline

    0.063 & 50    &   $\widehat{\bSigma}_{u}^{{\rm thre}}$     & 1     & 4.9(4.9) & 3.8(3.8) & 2.6(2.7) & 10(10.3) & 5.2(5.2) & 3.3(3.4) & 10.2(10.3) & 7.7(7.7) & 5.1(5.1) \\
          &       &       & 2     & 4.8(4.8) & 3.7(3.7) & 2.6(2.6) & 9.1(9.1) & 5.0(5.0)  & 3.2(3.3) & 9.6(9.7) & 7.3(7.3) & 4.9(4.8) \\
          &       &       & 3     & 4.8(4.8) & 3.7(3.7) & 2.6(2.6) & 8.5(8.4) & 5.0(5.0)  & 3.3(3.3) & 9.1(9.2) & 7.0(7.0)  & 4.7(4.7) \\
          &       &  $\widehat{\bSigma}_{u}^{{\rm bc, thre}}$      & 1     & 4.9(4.9) & 3.8(3.8) & 2.7(2.7) & 10.2(10.4) & 5.2(5.3) & 3.4(3.4) & 10.3(10.3) & 7.7(7.7) & 5.1(5.1) \\
          &       &       & 2     & 4.8(4.9) & 3.7(3.7) & 2.6(2.6) & 9.3(9.3) & 5.1(5.1) & 3.3(3.4) & 9.7(9.7) & 7.3(7.3) & 4.9(4.9) \\
          &       &       & 3     & 4.9(4.9) & 3.7(3.7) & 2.6(2.6) & 8.6(8.6) & 5.1(5.1) & 3.3(3.4) & 9.1(9.1) & 6.9(6.9) & 4.7(4.7) \\
          &       & $\widehat{\bSigma}_{u}^{*, {\rm thre}}$       &       & 8.0(8.0)  & 5.8(5.8) & 4.1(4.1) & 24.4(24.1) & 14.3(14.8) & 5.8(5.9) & 17.1(17.1) & 12.5(12.5) & 8.1(8.1) \\
          & 100   &   $\widehat{\bSigma}_{u}^{{\rm  thre}}$     & 1     & 5.3(5.3) & 4.1(4.1) & 2.9(2.9) & 12.6(12.4) & 6.2(5.9) & 3.5(3.5) & 11.6(11.7) & 8.7(8.8) & 5.7(5.7) \\
          &       &       & 2     & 5.2(5.2) & 4.0(4.0)  & 2.8(2.8) & 11.3(11.0) & 5.6(5.5) & 3.4(3.4) & 10.9(11.0) & 8.3(8.4) & 5.5(5.5) \\
          &       &       & 3     & 5.2(5.2) & 3.9(3.9) & 2.8(2.8) & 10.2(10.0) & 5.4(5.3) & 3.4(3.4) & 10.3(10.4) & 7.9(8.0) & 5.3(5.3) \\
          &       &   $\widehat{\bSigma}_{u}^{{\rm bc, thre}}$     & 1     & 5.3(5.3) & 4.1(4.1) & 2.9(2.9) & 12.7(12.5) & 6.3(6.0) & 3.5(3.6) & 11.7(11.8) & 8.7(8.7) & 5.7(5.8) \\
          &       &       & 2     & 5.2(5.2) & 4.0(4.0)  & 2.9(2.9) & 11.3(11.1) & 5.7(5.5) & 3.4(3.5) & 11.0(11.1) & 8.2(8.3) & 5.5(5.5) \\
          &       &       & 3     & 5.2(5.3) & 4.0(4.0)  & 2.9(2.8) & 10.3(10.1) & 5.5(5.4) & 3.4(3.5) & 10.3(10.4) & 7.8(7.9) & 5.3(5.3) \\
          &       &  $\widehat{\bSigma}_{u}^{*, {\rm thre}}$      &       & 9.5(9.4) & 6.2(6.2) & 4.4(4.4) & 31.5(31.1) & 18.5(18.3) & 7.1(6.7) & 19.6(19.5) & 14.2(14.2) & 9.1(9.1) \\
          & 200   &  $\widehat{\bSigma}_{u}^{{\rm thre}}$      & 1     & 5.6(5.6) & 4.3(4.3) & 3.1(3.1) & 14.7(13.6) & 7.0(6.5) & 3.4(3.3) & 13.0(13.0) & 9.6(9.6) & 6.3(6.3) \\
          &       &       & 2     & 5.5(5.5) & 4.2(4.2) & 3.1(3.1) & 12.9(11.9) & 6.3(5.8) & 3.3(3.3) & 12.2(12.2) & 9.2(9.2) & 6.0(6.1) \\
          &       &       & 3     & 5.6(5.5) & 4.2(4.2) & 3.0(3.0)  & 11.5(10.6) & 5.8(5.4) & 3.3(3.2) & 11.4(11.5) & 8.8(8.8) & 5.8(5.8) \\
          &       &  $\widehat{\bSigma}_{u}^{{\rm bc, thre}}$      & 1     & 5.7(5.7) & 4.4(4.4) & 3.2(3.2) & 14.8(13.7) & 7.0(6.5) & 3.4(3.4) & 13.1(13.1) & 9.5(9.6) & 6.3(6.3) \\
          &       &       & 2     & 5.6(5.6) & 4.3(4.3) & 3.1(3.1) & 13.0(12.0) & 6.4(5.9) & 3.3(3.3) & 12.3(12.3) & 9.1(9.1) & 6.1(6.1) \\
          &       &       & 3     & 5.6(5.6) & 4.3(4.3) & 3.1(3.1) & 11.6(10.7) & 5.8(5.5) & 3.3(3.3) & 11.5(11.6) & 8.7(8.7) & 5.9(5.9) \\
          &       & $\widehat{\bSigma}_{u}^{*, {\rm thre}}$       &       & 12.1(12.0) & 6.6(6.6) & 4.7(4.7) & 42.7(40.2) & 22.2(20.9) & 7.8(7.2) & 21.8(21.7) & 15.5(15.5) & 10.0(10.0) \\
\hline
    0.00252 & 50    &  $\widehat{\bSigma}_{u}^{{\rm  thre}}$      & 1     & 23.0(23.0) & 11.7(11.9) & 4.8(5.0) & 41.4(41.1) & 18.4(19.4) & 6.8(7.0) & 23.5(24.1) & 13.1(13.5) & 5.6(5.8) \\
          &       &       & 2     & 33.6(33.9) & 16.8(16.8) & 6.4(6.6) & 63.8(62.4) & 30.6(29.7) & 10.2(10.7) & 31.0(32.0) & 17.1(17.7) & 7.2(7.4) \\
          &       &       & 3     & 45.2(46.3) & 22.7(22.9) & 8.4(8.5) & 86.5(86.4) & 42.9(41.9) & 14.8(15.0) & 38.3(39.6) & 21.2(21.9) & 8.9(9.2) \\
          &       &$\widehat{\bSigma}_{u}^{{\rm bc, thre}}$        & 1     & 5.3(5.3) & 3.9(3.9) & 2.7(2.7) & 12.6(12.4) & 5.8(5.8) & 3.4(3.5) & 11.3(11.4) & 8.0(8.0)  & 5.1(5.1) \\
          &       &       & 2     & 5.4(5.4) & 3.9(3.9) & 2.7(2.7) & 12.7(12.1) & 5.8(5.6) & 3.4(3.5) & 11.1(11.2) & 7.7(7.7) & 5.0(5.0) \\
          &       &       & 3     & 5.6(5.5) & 3.9(3.9) & 2.6(2.6) & 13.3(12.2) & 5.9(5.7) & 3.5(3.5) & 11.0(11.1) & 7.5(7.5) & 4.8(4.8) \\
          &       &  $\widehat{\bSigma}_{u}^{*, {\rm thre}}$      &       & 8.9(8.8) & 5.9(5.9) & 4.1(4.1) & 26.3(25.5) & 15.2(15.4) & 5.9(5.9) & 18.1(18.0) & 12.8(12.8) & 8.1(8.1) \\
          & 100   &  $\widehat{\bSigma}_{u}^{{\rm thre}}$      & 1     & 27.3(26.9) & 14.0(13.7) & 5.8(5.4) & 41.1(39.8) & 18.7(18.0) & 6.9(6.6) & 28.7(28.0) & 15.5(15.3) & 6.7(6.7) \\
          &       &       & 2     & 41.3(39.7) & 20.1(19.7) & 7.8(7.6) & 63.0(61.0) & 29.6(28.5) & 10.3(9.9) & 37.6(36.3) & 20.3(20.4) & 8.8(8.7) \\
          &       &       & 3     & 57.1(54.6) & 27.8(27.0) & 10.2(10.0) & 86.9(84.0) & 42.5(40.6) & 14.7(14.0) & 46.3(44.3) & 25.3(25.5) & 10.9(10.8) \\
          &       &  $\widehat{\bSigma}_{u}^{{\rm bc, thre}}$      & 1     & 5.8(5.8) & 4.2(4.2) & 3.0(3.0)  & 16.2(15.5) & 7.3(6.9) & 3.6(3.7) & 13.0(13.1) & 9.0(9.1) & 5.8(5.8) \\
          &       &       & 2     & 5.9(5.9) & 4.2(4.2) & 2.9(2.9) & 16.3(15.2) & 6.9(6.7) & 3.6(3.6) & 12.9(12.9) & 8.8(8.8) & 5.6(5.6) \\
          &       &       & 3     & 6.1(6.1) & 4.2(4.2) & 2.9(2.9) & 17.0(15.4) & 6.8(6.6) & 3.6(3.6) & 12.9(13.0) & 8.5(8.5) & 5.4(5.5) \\
          &       &  $\widehat{\bSigma}_{u}^{*, {\rm thre}}$      &       & 11.5(11.3) & 6.3(6.3) & 4.4(4.4) & 33.9(33.5) & 19.6(19.2) & 7.3(7.0) & 20.7(20.6) & 14.5(14.6) & 9.2(9.2) \\
          & 200   &  $\widehat{\bSigma}_{u}^{{\rm thre}}$      & 1     & 27.6(27.2) & 14.2(13.9) & 5.9(5.6) & 33.7(31.1) & 14.7(13.7) & 5.4(5.2) & 29.4(29.4) & 15.6(15.7) & 7.0(7.0) \\
          &       &       & 2     & 41.7(40.6) & 20.3(20.0) & 7.9(7.6) & 48.6(46.8) & 22.5(21.3) & 7.7(7.4) & 38.1(37.9) & 20.6(20.8) & 9.0(9.0) \\
          &       &       & 3     & 57.1(56.8) & 27.8(27.4) & 10.2(10.0) & 65.4(64.4) & 31.6(30.7) & 10.9(10.4) & 46.5(46.5) & 25.8(26.0) & 11.2(11.2) \\
          &       &  $\widehat{\bSigma}_{u}^{{\rm bc, thre}}$      & 1     & 6.2(6.1) & 4.5(4.5) & 3.2(3.2) & 19.0(17.7) & 8.1(7.4) & 3.5(3.5) & 14.5(14.5) & 10.0(10.0) & 6.4(6.4) \\
          &       &       & 2     & 6.4(6.2) & 4.5(4.5) & 3.1(3.1) & 18.9(17.4) & 7.8(7.1) & 3.4(3.4) & 14.4(14.3) & 9.7(9.7) & 6.2(6.2) \\
          &       &       & 3     & 6.8(6.5) & 4.5(4.5) & 3.1(3.1) & 19.4(17.7) & 7.6(7.0) & 3.4(3.4) & 14.4(14.4) & 9.4(9.4) & 6.0(6.0) \\
          &       &  $\widehat{\bSigma}_{u}^{*, {\rm thre}}$      &       & 13.8(13.8) & 6.7(6.7) & 4.7(4.8) & 46.1(43.5) & 23.6(22.2) & 8.1(7.4) & 22.9(22.8) & 15.9(16.0) & 10.1(10.1) \\
		\hline		\hline
	\end{tabular}\label{Table-IndepN}
\end{table}%

\end{document}